\documentclass[11pt,reqno,twoside]{article}

\usepackage[hang]{footmisc}
\usepackage{lipsum}
\usepackage{amsmath}

\usepackage{cmap} 

\usepackage[T1]{fontenc}
\usepackage[utf8]{inputenc}
\usepackage{graphicx}
\usepackage{placeins}
\usepackage{enumerate}
\usepackage{algcompatible}
\usepackage{subcaption}
\usepackage{booktabs}
\usepackage{verbatim}
\newcommand{\comments}[1]{}

\usepackage{soul}

\usepackage{setspace}

\let\counterwithin\relax  
\usepackage{lmodern} 
\usepackage[scale=0.88]{tgheros} 

\usepackage{bm} 

\usepackage{bbold}

\usepackage{pgfplots}

\usepackage{amsmath,amsbsy,amsgen,amscd,amsthm,amsfonts,amssymb} 
\allowdisplaybreaks[4]

\usepackage[centering,top=1in,bottom=1in,left=1in,right=1in]{geometry}

\usepackage{titling}
\graphicspath{ {./images/} }

\usepackage[sf,bf,compact]{titlesec}

\definecolor{dark-gray}{gray}{0.3}
\definecolor{dkgray}{rgb}{.4,.4,.4}
\definecolor{dkblue}{rgb}{0,0,.5}
\definecolor{medblue}{rgb}{0,0,.75}
\definecolor{rust}{rgb}{0.5,0.1,0.1}

\usepackage{url}
\usepackage[colorlinks=true]{hyperref}
\hypersetup{linkcolor=dkblue}    
\hypersetup{citecolor=rust}      
\hypersetup{urlcolor=rust}     

\usepackage[final]{microtype} 

\newtheoremstyle{myThm} 
    {\topsep}                    
    {\topsep}                    
    {\itshape}                   
    {}                           
    {\sffamily\bfseries}                   
    {.}                          
    {.5em}                       
    {}  

\newtheoremstyle{myRem} 
    {\topsep}                    
    {\topsep}                    
    {}                   
    {}                           
    {\sffamily}                   
    {.}                          
    {.5em}                       
    {}  

\newtheoremstyle{myDef} 
    {\topsep}                    
    {\topsep}                    
    {}                   
    {}                           
    {\sffamily\bfseries}                   
    {.}                          
    {.5em}                       
    {}  

\theoremstyle{myThm}
\newtheorem{theorem}{Theorem}[section]
\newtheorem{lemma}[theorem]{Lemma}
\newtheorem{proposition}[theorem]{Proposition}
\newtheorem{corollary}[theorem]{Corollary}

\newtheorem{dataassumption}[theorem]{Data Assumption}

\theoremstyle{myRem}
\newtheorem{remarkx}[theorem]{Remark}
 \newenvironment{remark}
  {\pushQED{\qed}\remarkx}
  {\popQED\endremarkx}

\theoremstyle{myDef}
\newtheorem{definition}[theorem]{Definition}

\usepackage{fancyhdr}
\let\originalleft\left
\let\originalright\right
\renewcommand{\left}{\mathopen{}\mathclose\bgroup\originalleft}
\renewcommand{\right}{\aftergroup\egroup\originalright}

\usepackage{mathtools}
\mathtoolsset{centercolon}  

\usepackage[normalem]{ulem}

\definecolor{mygreen}{rgb}{0.1,0.75,0.2}

\newcommand{\nc}{\normalcolor}

\providecommand{\mathbbm}{\mathbb} 

\newcommand{\R}{\mathbbm{R}}

\newcommand{\N}{\mathbbm{N}}

\newcommand{\G}{\mathcal{G}}
\newcommand{\F}{\mathcal{F}}

\renewcommand{\L}{\mathcal{L}}

\renewcommand{\phi}{\varphi}
\newcommand{\eps}{\varepsilon}

\usepackage[font = small, margin=30pt]{caption}

\usepackage[]{algorithm}
\usepackage{algpseudocode}
\usepackage{enumerate}

\usepackage{authblk}
\usepackage[square,numbers]{natbib}
\usepackage{chngcntr}
\usepackage{mathrsfs}
\counterwithin{table}{section}
\counterwithin{algorithm}{section}

\newcommand{\iid}{\stackrel{\mathsf{i.i.d.}}{\sim}}

\usepackage{multirow}
\usepackage{enumitem}

\usepackage[scr = esstix, cal = cm, frak=euler]{mathalfa}

\definecolor{mygreen}{rgb}{0.1,0.75,0.2}

\newcommand{\black}{\color{black}}

\usepackage[scr = esstix, cal = cm, frak=euler]{mathalfa}

\usepackage{authblk}
\title{\huge Elliptic Bayesian Inverse Problems on Metric Graphs}
\author[1]{\Large D. Bolin}
\author[,2]{\Large W. Li\footnote{Corresponding author: wenwenli@uchicago.edu}}
\author[3]{\Large D. Sanz-Alonso}
\affil[1]{\normalsize CEMSE Division, King Abdullah University of Science and Technology, Saudi Arabia}
\affil[2]{\normalsize Committee on Computational and Applied Mathematics, University of Chicago, USA}
\affil[3]{\normalsize Department of Statistics, University of Chicago, USA}
\date{}

\makeatletter\@addtoreset{section}{part}\makeatother%
\numberwithin{equation}{section}

\newcommand{\upperRomannumeral}[1]{\uppercase\expandafter{\romannumeral#1}}

\renewcommand{\hat}{\widehat}

\providecommand{\keywords}[1]{\textbf{{Keywords:}} #1}

\begin{document}
\maketitle 


\abstract{This paper studies the formulation, well-posedness, and numerical solution of Bayesian inverse problems on metric graphs, in which the edges represent one-dimensional wires connecting vertices.  
We focus on the inverse problem of recovering the diffusion coefficient of a (fractional) elliptic equation on a metric graph from noisy measurements of the solution. Well-posedness hinges on both stability of the forward model and an appropriate choice of prior. 
We establish the stability of elliptic and fractional elliptic forward models using recent regularity theory for differential equations on metric graphs. For the prior,  we leverage modern Gaussian Whittle--Matérn process models on metric graphs with sufficiently smooth sample paths. 
Numerical results demonstrate accurate reconstruction and effective uncertainty quantification.}

\bigskip


\bigskip

\section{Introduction}
Differential equations on metric graphs are  
used to model
diverse phenomena ranging from 
 traffic flow in road networks \cite{bolin2025log} to wave propagation in thin structures arising in mesoscopic systems, photonic crystals, and nanotechnology \cite{kuchment2002graph}. These and other applications  motivate the study of inverse problems for differential equations on metric graphs.  In this work, we focus on the inverse problem of recovering the diffusion coefficient of a (fractional) elliptic partial differential equation (PDE) from noisy measurements of the solution. Importantly, the unknown to-be-reconstructed diffusion coefficient represents a function defined along the edges of the metric graph.  

 We will adopt the Bayesian approach to inverse problems \cite{kaipio2006statistical,stuart2010inverse,calvetti2007introduction,sanzstuarttaeb}, placing a prior on the unknown parameter, and conditioning on the observed data to obtain its posterior distribution. In the Bayesian approach, the posterior not only yields point estimators for the unknown parameter, but also enables to quantify the uncertainty in the reconstruction. 
 Furthermore, the Bayesian formulation can result in a form of \emph{well-posedness}, by which small perturbations in the data, the prior, or the forward model yield small perturbations in the posterior \cite{stuart2010inverse,latz2020well}.
In this paper, we are concerned only with well-posedness with respect to perturbations in the data. This notion hinges on stability of the forward map together with an appropriate choice of prior. Ensuring well-posedness is particularly delicate when the to-be-reconstructed parameter is a function, as is the case in this work.

In Euclidean domains, well-posedness of Bayesian inversion in function space has been extensively studied, see e.g. \cite{stuart2010inverse,dashti2013bayesian,latz2020well} and references therein. The papers \cite{dashti2011uncertainty,trillos2016bayesian} established the formulation and well-posedness of elliptic and fractional elliptic inverse problems that have become standard test-beds for theory and methods for Bayesian inversion and uncertainty quantification. More recently, the works \cite{harlim2020kernel,harlim2022graph,kim2024optimization} studied Bayesian inversion and optimization for differential equations on manifolds and point clouds. Relatedly, \cite{trillos2017consistency,garcia2018continuum,trillos1718mathematical} formulate semi-supervised learning as a Bayesian inverse problem on a graph, where vertices represent labeled and unlabeled features and edges encode similarities between them. In that setting, the goal is to assign labels to all the features, which involves specifying a prior for functions defined on the \emph{vertices}  \cite{sanz2022spde}, and conditioning on the labeled data to find the posterior.

To our knowledge, this is the first paper to study the formulation and well-posedness of Bayesian inversion on metric graphs. Our goal is to recover a function defined on the \emph{edges} that represents an unknown parameter of a differential equation on a metric graph.  Extending the theory from the Euclidean, manifold, and graph settings to the metric graph setting is far from straightforward: the piecewise one-dimensional geometry together with the vertex coupling conditions introduce new regularity and stability challenges. 
\black In particular, the presence of vertices breaks the smooth-manifold structure, so one cannot directly invoke standard global elliptic regularity results as in Euclidean or manifold settings. 
The analysis must instead combine edgewise (one-dimensional) estimates with the vertex coupling conditions, and track how these estimates enter the stability bounds for the forward map. \nc
On the other hand, since metric graphs are locally one-dimensional, some classical Euclidean estimates can be improved. 
We analyze the stability of forward models for elliptic and fractional elliptic inverse problems building on recent regularity theory for differential equations on metric graphs developed in \cite{bolin2024regularity}. To ensure well-posedness, we leverage Gaussian Whittle--Matérn prior models on metric graphs with sufficiently smooth sample paths, introduced and analyzed in \cite{bolin2024gaussian}.

The rest of this paper is organized as follows. Section~\ref{sec: background} contains background on metric graphs. Section~\ref{sec:problemformulation} introduces the formulation of Bayesian inversion on metric graphs. Section~\ref{sec:well-posedness} contains our main technical contributions; we establish continuous dependence (in Hellinger distance) of the posterior distribution on the data by developing new stability theory for forward models on compact metric graphs. Section~\ref{sec:numerics} reports numerical experiments for elliptic and fractional elliptic inverse problems on a letter-shaped metric graph. Our numerical results demonstrate both accurate reconstruction and successful uncertainty quantification. The code is publicly available at: \url{https://github.com/WenwenLi2002/Bayesian-Inverse-Problems-on-Metric-Graphs}. Finally, Section~\ref{sec:conclusions} closes with conclusions and open directions for future research.

\section{Background}\label{sec: background}
This section introduces necessary background on metric graphs, as well as the function spaces in which we will formulate and analyze Bayesian inversion on metric graphs. 

\subsection{Metric Graphs and Function Spaces}
Let $\Gamma$ be a graph with vertices $V = \{v_i\}$ and edges $E = \{e_j\}$. We are interested in applications where the edges do not represent abstract relationships between vertices, but rather physical one-dimensional wires connecting them. 
To model this scenario, we assign to each edge $e\in E$ a positive length $\ell_e>0$, and then we orient each edge arbitrarily and identify it with the interval $[0,\ell_e]$ via a coordinate $t_e$. A graph $\Gamma$ supplemented with this structure is called a \emph{metric graph}, where the metric is naturally given by the shortest path distance \cite{berkolaiko2013introduction}. A generic point $x$ on a metric graph $\Gamma$ can be represented as $x = (e,t_e)$ for some $e \in E$ and $t_e \in [0, \ell_e].$  Thus, points on a metric graph include not only the vertices (for which $t_e \in \{0, \ell_e \}$), but also intermediate points on the edges (for which $t_e \in (0, \ell_e)$).   
In this paper, we work with \emph{compact metric graphs} comprising finitely many vertices and edges, with every edge of finite length. It can be shown that a metric graph equipped with its natural metric defines a compact metric space \cite{berkolaiko2013introduction}.

A function \(f\) on a metric graph \(\Gamma\) can be represented as a collection \(\{f_e\}_{e \in E}\), where \(f_e\) is a function defined on the edge \(e\cong [0,\ell_{e}]\). We write 
$$ \int_\Gamma f \,  {\rm dx} : = \sum_{e \in E} \int_0^{\ell_e}  f_e(t)  \, {\rm dt}. $$
Function spaces on $\Gamma$ can be defined in terms of function spaces on the edges. For instance, we let $L^2(e)$ be the space of Lebesgue square-integrable functions on $e\cong [0,\ell_{e}]$ and define
\[
  L^2(\Gamma)
  := \bigoplus_{e\in E}L^2(e)
  = \Bigl\{\,f = (f_e)_{e\in E}\;\Bigm|\;f_e\in L^2(e),\quad \forall e\in E\Bigr\}, \quad \|f\|_{L^2(\Gamma)}^2
  := \sum_{e\in E}\|f_e\|_{L^2(e)}^2.
\]
In many applications, we are interested in global functions with additional regularity in each edge. Furthermore, appropriate matching conditions at the vertices ---such as continuity or Kirchhoff-type conditions--- may be imposed to ensure that the edgewise functions are appropriately connected, resulting in a function on \(\Gamma\) with desired global behavior.
Table~\ref{tab:FractionalSobolev} summarizes several function spaces defined on \(\Gamma\) that will be used in this paper. These include spaces defined on each edge and spaces defined on the entire graph, along with their associated norms and dual spaces.  Together, these function spaces provide a rigorous framework for studying differential equations  on metric graphs.

\begin{table}[h]
\centering
\caption{Function spaces on \(\Gamma\) and on edges $e\in E$.}
\label{tab:FractionalSobolev}
\resizebox{.875\textwidth}{!}{%
\begin{tabular}{l l l}
\toprule
\textbf{Notation} & \textbf{Definition} & \textbf{Norm} \\ 
\midrule
\(H^1(e),\, e\in E\) & See \cite[p.~77]{mclean2000strongly} & \(\|\cdot\|_{H^1(e)}\) (see \cite[p.~77]{mclean2000strongly}) \\
\(H^s(e),\, e\in E\) & See \cite[p.~77]{mclean2000strongly} & \(\|\cdot\|_{H^s(e)}\) (see \cite[p.~77]{mclean2000strongly}) \\
\(H^s(\Gamma),\, s \in (0,1) \) & \((L^2(\Gamma),H^1(\Gamma))_s\) & \(\|\cdot\|_{H^s(\Gamma)}:=\|\cdot\|_{(L^2(\Gamma),H^1(\Gamma))_s}\) \\
\(L^\infty(\Gamma)\) & \(\displaystyle \Bigl\{f:\Gamma\to\mathbb{R}: \|f\|_{L^\infty(\Gamma)}=\operatorname*{ess\,sup}_{x\in\Gamma}|f(x)|<\infty\Bigr\}\) & \(\displaystyle \|f\|_{L^\infty(\Gamma)}=\operatorname*{ess\,sup}_{x\in\Gamma}|f(x)|\) \\
\(C(\Gamma)\) & \(\displaystyle \{f\in L^2(\Gamma):\; f \text{ is continuous on } \Gamma\}\) & \(\displaystyle \|f\|_{C(\Gamma)}=\sup_{x\in\Gamma}|f(x)|\) \\
\( C^{0,s}(\Gamma), \, s \in (0,1]\) & See \cite[Section 2]{bolin2024gaussian} & See \cite[Section 2]{bolin2024gaussian} \\
\(H^1(\Gamma)\) & \(\displaystyle \Bigl\{ f\in \bigoplus_{e\in E} H^1(e):\; f \text{ is continuous on } \Gamma\Bigr\}\) & \(\displaystyle \|f\|_{H^1(\Gamma)}=\Bigl(\sum_{e\in E}\|f\|_{H^1(e)}^2\Bigr)^{1/2}\) \\
\(H^{-1}(\Gamma)\) & Dual space of \(H^1(\Gamma)\) & Dual norm \\
\(H^s(\Gamma),\, s \in (1,2)\) &  \((H^1(\Gamma),\widetilde{H}_C^2(\Gamma))_{s-1}\)  & \(\|\cdot\|_{H^s(\Gamma)}:=\|\cdot\|_{(H^1(\Gamma),\widetilde{H}_C^2(\Gamma))_{ s-1 }}\) \\
\(\widetilde{H}^s(\Gamma),\, s \in (0,1)\) & \((L^2(\Gamma),\widetilde{H}^1(\Gamma))_s\) & \(\|\cdot\|_{\widetilde{H}^s(\Gamma)}:=\|\cdot\|_{(L^2(\Gamma),\widetilde{H}^1(\Gamma))_s}\) \\
\(\widetilde{H}^s(\Gamma),\, s \in (1,2)\)  &  \((H^1(\Gamma),\widetilde{H}^2(\Gamma))_{s-1}\)  & \(\|\cdot\|_{\widetilde{H}^s(\Gamma)}:=\|\cdot\|_{(\widetilde{H}^1(\Gamma),\widetilde{H}^2(\Gamma))_{ s-1}}\) \\
\black \(\widetilde{H}_C^2(\Gamma)\) \nc
& \black \(\displaystyle \widetilde{H}^2(\Gamma)\cap C(\Gamma)
   = \Bigl\{f\in \bigoplus_{e\in E} H^2(e):\ f\ \text{is continuous on }\Gamma\Bigr\}\) \nc
& \black \(\displaystyle \|f\|_{\widetilde{H}_C^2(\Gamma)}:=\|f\|_{\widetilde{H}^2(\Gamma)}
   =\Bigl(\sum_{e\in E}\|f\|_{H^2(e)}^2\Bigr)^{1/2}\) \nc
\\
\bottomrule
\end{tabular}%
}
\end{table}

Here, for two Hilbert spaces $\mathcal{H}_1 \subset \mathcal{H}_0,$ we denote by $(\mathcal{H}_0,\mathcal{H}_1)_s$ the real interpolation of order $s \in (0,1)$ between $\mathcal{H}_0$ and $\mathcal{H}_1$ using the $K$-method ---see \cite[Appendix A]{bolin2024gaussian}. 
Furthermore, $\widetilde{H}^r(\Gamma) = \bigoplus_{e\in E} H^r(e)$ for $r\in\mathbb{N}$ is a decoupled Sobolev space defined as the direct sum of the corresponding Sobolev spaces on the edges, and 
$\widetilde{H}_C^r(\Gamma)$ means $\widetilde{H}_C^r(\Gamma) 
= \widetilde{H}^r(\Gamma) \,\cap\, C(\Gamma)$ for  $r \in \mathbb{N}.$ We note in particular that the space \(\widetilde{H}^1(\Gamma)\) is defined as the decoupled direct sum of \(H^1(e)\) spaces on each edge whereas the space \(H^1(\Gamma)\) requires global continuity over the entire metric graph \(\Gamma\). This distinction is crucial when studying differential equations on \(\Gamma\), as the imposition of continuity conditions at the vertices influences the analytical properties of the solutions.

\subsection{Quantum Graphs}\label{ssec:quantum_graphs}
A metric graph is called a \emph{quantum graph} when equipped with a differential operator and appropriate \emph{vertex conditions}. \black The term quantum graph, which is standard  in the literature (see \cite{berkolaiko2013introduction,kuchment2008quantum}), originates from mathematical physics, where operators on metric graphs serve as effective Schr\"odinger-type Hamiltonians modeling wave or quantum-particle propagation on thin network-like structures, e.g.\ wire or waveguide networks, with vertex conditions encoding the coupling at junctions.
\nc
We first introduce a general second-order elliptic differential operator \( \mathcal{L}\) on a compact metric graph \(\Gamma\). Consider a sufficiently smooth function \(p\) defined on \(\Gamma\). Locally, on each edge \(e \in E\), the operator \(\mathcal{L}\) acts on $p$ by
\begin{equation}
\label{eq:edge_ODE}
- \frac{d}{dt}\Bigl(a_e(t)\,\frac{d}{dt} p_e(t)\Bigr) + \kappa_e^2(t)\,p_e(t), \quad t \in (0,\ell_e),
\end{equation}
where \(p_e\) is the restriction of \(p\) to the edge \(e\), \(a_e(t)\) is a \black continuous \nc and strictly positive function bounded away from zero, and \(\kappa_e(t)\) is bounded and has a strictly positive lower bound. 

To formulate a global problem $\mathcal{L}p = f$ on \(\Gamma\) for some $f\in L^2(\Gamma)$, we impose coupling conditions at each vertex. Specifically, we require the solution \(p\) to be continuous at all vertices, and we enforce generalized Kirchhoff-type conditions at vertices given explicitly by 
\begin{equation}
\label{eq:kirchhoff}
p \text{ is continuous on } \Gamma, \quad \text{and} \quad 
\forall v \in V:\quad  \sum_{e\in E_v} a_e(v)\partial_e p(v) = \theta p(v),
\end{equation}
where \(E_v\) is the set of edges incident to vertex \(v\), \black \(a:=\bigoplus_{e\in E} a_e \in C(\Gamma)\)\nc, \(\partial_e p(v)\)  denotes the outward-directed derivative of \(p\) at vertex \(v\) along edge \(e\), and \(\theta \in \mathbb{R}\) is a given parameter. The standard Kirchhoff conditions correspond to \(\theta = 0\), which ensures that flux is conserved at every vertex. 
\black Under these conditions, we will show in Proposition~\ref{prop:forward_under_H1} below \nc
that the operator \(\mathcal{L}\)  is self-adjoint on \(L^2(\Gamma)\), strictly positive, and has compact inverse. Consequently, \(\mathcal{L}\) possesses a discrete spectrum consisting of eigenvalues \(\{\lambda_j\}_{j=1}^{\infty}\) in non-decreasing order, with corresponding eigenfunctions \(\{\phi_j\}_{j=1}^{\infty}\) forming a complete orthonormal basis of \(L^2(\Gamma)\). \nc
Therefore, for \(\beta>0\), we can define the fractional power \(\L^\beta\) of the operator \(\mathcal{L}\)  via its spectral decomposition as follows:
\begin{equation*}
\mathcal{L}^\beta p = \sum_{j=1}^{\infty}\lambda_j^\beta \langle p,\phi_j\rangle_{L^2(\Gamma)}\phi_j, \quad p\in\mathcal{D}(\mathcal{L}^\beta),
\end{equation*}
where $\mathcal{D}(\mathcal{L}^\beta)$ is the natural domain of $\mathcal{L}^\beta$, and $\langle f, g \rangle_{L^2(\Gamma)} = \int_{\Gamma} f g \, {\rm dx}$ is the $L_2(\Gamma)$ inner product. This leads to the fractional elliptic differential equation
\begin{equation*}
\mathcal{L}^\beta p = f \quad \text{on } \Gamma,
\end{equation*}
with given \(f \in L^2(\Gamma)\), which we will study in detail below.

\section{Inverse Problem Formulation}\label{sec:problemformulation}
Let $\Gamma$ be, here and throughout, a compact metric graph. For given $\beta\geq 1$, consider the (fractional) elliptic differential equation: 
\begin{equation}
\label{equa:p}
\begin{aligned}
    & \mathcal{L}_u^\beta \,p := \Bigl(\kappa^2  - \nabla \cdot \bigl(e^u\, \nabla\bigr)\Bigr)^\beta p = f \quad \text{on } \Gamma,
\end{aligned}
\end{equation}
where $\mathcal{L}_u$ is equipped with the Kirchhoff vertex conditions \eqref{eq:kirchhoff}, and $f\in L^2(\Gamma)$ represents a given source term. 
\black In this work, we adopt a log-parameterization of the diffusion coefficient and set
\(
a := e^{u}
\)
in Equation~\ref{equa:p}, so that the unknown field \(u\) enters the model through the strictly positive diffusion coefficient \(a\).
We will endow \(u\) with a Gaussian prior and choose its support so that \(u\) (and hence \(a=e^{u}\)) is sufficiently regular for the assumptions on \(a\) in Section~\ref{ssec:quantum_graphs} to hold.
\nc
For convenience, we only consider the standard Kirchhoff vertex condition (\(\theta=0\) in Equation~\ref{eq:kirchhoff}) and  constant coefficient \(\kappa>0\). \black This restriction is made to align with the available regularity theory for fractional elliptic problems on compact metric graphs, which is developed for the standard Kirchhoff case (see \cite{bolin2024regularity}). \nc We will refer to $\beta =1$ as the \emph{elliptic case} and to $\beta >1$ as the \emph{fractional case}. Our goal is to find the unknown function $u$ specifying the diffusion coefficient of $\mathcal{L}_u$ from noisy measurements of the solution $p.$

\black
\begin{proposition}[Properties of $\L_u$]
\label{prop:forward_under_H1}
Let \(u\in H^1(\Gamma)\) and set \(a:=e^{u}\).
Assume \(\kappa\in L^\infty(\Gamma)\) with \(\operatorname*{ess\,inf}_{x\in\Gamma}\kappa(x)\ge \kappa_{\min}>0\).
Let \(\L_u:=\kappa^2-\nabla\!\cdot(a\nabla)\) be endowed with standard Kirchhoff vertex conditions.
Then the following hold:
\begin{enumerate}
\item \textbf{Self-adjointness and positivity.}
The symmetric bilinear form on \(H^1(\Gamma)\) given by
\[
\mathfrak a_u(\phi,\psi)
:=\int_\Gamma a\,\nabla\phi\cdot\nabla\psi\, {\rm dx} +\int_\Gamma \kappa^2\,\phi\psi\, {\rm dx},
\qquad \phi,\psi\in H^1(\Gamma),
\]
is closed and coercive. The associated Kirchhoff realization \(\L_u\) is self-adjoint and strictly positive on \(L^2(\Gamma)\).

\item \textbf{Compact resolvent.}
\(\L_u^{-1}:L^2(\Gamma)\to L^2(\Gamma)\) is compact; in particular \(\L_u\) has purely discrete spectrum and
\(\L_u^\beta\) is well-defined by spectral calculus for every \(\beta>0\).

\item \textbf{Fractional equation and \(H^1\)-solution.}
For any \(\beta\ge 1\) and \(f\in L^2(\Gamma)\), the fractional equation
\[
\L_u^\beta p=f
\]
admits a unique solution \(p=\L_u^{-\beta}f\in D(\L_u)\subset H^1(\Gamma)\).

\item \textbf{Kirchhoff condition in the weak sense.}
The solution \(p\) satisfies the standard Kirchhoff vertex condition (cf. Equation~\ref{eq:kirchhoff}) in the weak sense:
the vertex boundary terms vanish in the edgewise integration-by-parts identity associated with \(\mathfrak a_u\).
\end{enumerate}
\end{proposition}
\begin{proof}
On a compact metric graph $\Gamma$ we have the continuous embedding $H^1(\Gamma)\hookrightarrow \black C(\Gamma) \nc$.
Since $u\in H^1(\Gamma)$, it is continuous and bounded on $\Gamma$.
Therefore $a:=e^u\in \black C(\Gamma) \nc$ and there exist constants $0<a_{\min}\le a_{\max}<\infty$ such that
\[
a_{\min}\le a(x)\le a_{\max}\qquad \forall\,x\in\Gamma .
\]

\smallskip
\noindent\emph{(1).}
By $a\le a_{\max}$ and $\kappa\in L^\infty(\Gamma)$,
\[
|\mathfrak a_u(\phi,\psi)|
\le a_{\max}\|\nabla\phi\|_{L^2(\Gamma)}\|\nabla\psi\|_{L^2(\Gamma)}
+\|\kappa\|_{L^\infty(\Gamma)}^2\|\phi\|_{L^2(\Gamma)}\|\psi\|_{L^2(\Gamma)}
\lesssim \|\phi\|_{H^1(\Gamma)}\|\psi\|_{H^1(\Gamma)} .
\]
Moreover, coercivity follows from $a_{\min}>0$ and $\kappa_{\min}>0$:
\[
\mathfrak a_u(\phi,\phi)
\ge a_{\min}\|\nabla\phi\|_{L^2(\Gamma)}^2+\kappa_{\min}^2\|\phi\|_{L^2(\Gamma)}^2
\gtrsim \|\phi\|_{H^1(\Gamma)}^2 .
\]
Hence the form norm induced by $\mathfrak a_u$ is equivalent to the $H^1(\Gamma)$-norm, so $\mathfrak a_u$ is closed and coercive on $H^1(\Gamma)$.
By the representation theorem for closed, densely-defined, coercive symmetric forms
(see, e.g., \cite[Ch.~VI]{kato1995perturbation} or \cite[Theorem~5]{kuchment2008quantum}),
there exists a unique self-adjoint operator $\L_u$ on $L^2(\Gamma)$ associated with $\mathfrak a_u$.
Strict positivity follows from coercivity.

\smallskip
\noindent\emph{(2).}
For each $f\in L^2(\Gamma)$, by Lax--Milgram, there exists a unique $p\in H^1(\Gamma)$ such that
\[
\mathfrak a_u(p,\varphi)=(f,\varphi)_{L^2(\Gamma)}\qquad \forall\,\varphi\in H^1(\Gamma),
\]
and $\|p\|_{H^1(\Gamma)}\lesssim \|f\|_{L^2(\Gamma)}$.
Thus $\L_u^{-1}:L^2(\Gamma)\to H^1(\Gamma)$ is bounded.
Since the embedding $H^1(\Gamma)\hookrightarrow L^2(\Gamma)$ is compact (compactness of $\Gamma$),
the composition $\L_u^{-1}:L^2(\Gamma)\to L^2(\Gamma)$ is compact.
Hence $\L_u$ has compact resolvent and spectral calculus defines $\L_u^\beta$ for all $\beta>0$.

\smallskip
\noindent\emph{(3).}
Let $\beta\ge 1$ and $f\in L^2(\Gamma)$.
Define $g:=\L_u^{-(\beta-1)}f\in L^2(\Gamma)$ and $p:=\L_u^{-1}g=\L_u^{-\beta}f$.
Then $p\in D(\L_u)$ by construction, hence $p\in H^1(\Gamma)$.
Uniqueness follows from strict positivity: if $\L_u^\beta p=0$, then
\[
0=(\L_u^\beta p,p)_{L^2(\Gamma)}=\|\L_u^{\beta/2}p\|_{L^2(\Gamma)}^2,
\]
so $p=0$.

\smallskip
\noindent\emph{(4).}
Given that $g:=\L_u^{-(\beta-1)}f\in L^2(\Gamma)$, the fractional equation $\L_u^\beta p=f$ can be rewritten as
\(\L_u p=g,\)
which is equivalent to
\[
\mathfrak a_u(p,\varphi)=(g,\varphi)_{L^2(\Gamma)}\qquad \forall\,\varphi\in H^1(\Gamma).
\]
On each edge $e$ we have
\[
-(a\,\partial_e p)' = g-\kappa^2 p \in L^2(e),
\]
so $a\,\partial_e p\in H^1(e)$. Since $H^1(e)\hookrightarrow \black C(\bar e)\nc$ in one dimension,
the flux $a\,\partial_e p$ is continuous on $\bar e$ and its endpoint values are well-defined.
Performing integration by parts edgewise gives
\[
\sum_{e}\int_e a\,\partial_e p\,\partial_e\varphi\,ds
= -\sum_{e}\int_e (a\,\partial_e p)'\,\varphi\,ds
+ \sum_{v}\Big(\sum_{e\in E_v}(a\,\partial_e p)(v)\Big)\,\varphi(v).
\]
Since $\L_u p=g$ implies $-\nabla\!\cdot(a\nabla p)+\kappa^2 p=g$ in $L^2(\Gamma)$, the edgewise integration-by-parts identity reduces the weak formulation to
\[
\sum_{v}\Big(\sum_{e\in E_v}(a\,\partial_e p)(v)\Big)\,\varphi(v)=0
\qquad \forall\,\varphi\in H^1(\Gamma).
\]
Choosing test functions $\varphi$ that satisfy $\varphi(v_0)=1$ and $\varphi(v)=0$ for all other vertices $v\neq v_0$
yields, for each vertex $v_0$,
\[
\sum_{e\in E_{v_0}}(a\,\partial_e p)(v_0)=0,
\]
which is the Kirchhoff condition in the weak sense.
For the standard Kirchhoff condition in the unweighted form $\sum_{e\in E_v}\partial_e p(v)=0$,
then since $a$ is continuous on $\Gamma$ and $a(v)>0$, we have $(a\,\partial_e p)(v)=a(v)\partial_e p(v)$ and the two formulations are equivalent.
\end{proof}

\nc
 Thanks to \black Proposition~\ref{prop:forward_under_H1}, \nc we can define the \emph{forward map} from parameter $u$ to solution $p$:
\begin{align}\label{equa:forward map}
\mathcal{F}: \black {H}^{1}\nc(\Gamma) & \to H^1(\Gamma),\\
 u & \mapsto p. \nonumber
\end{align}
We assume to have access to $m$ noisy measurements of $p$ of the form 
\begin{equation}
\label{eq:obs_data}
y_j = l_j(p) + \eta_j, \quad j = 1, \ldots, m,
\end{equation}
where $\{\eta_j\}_{j=1}^m$ represent measurement errors, not necessarily independent. 
\black
We assume that the data are obtained by applying finitely many bounded linear functionals to the state:

\begin{dataassumption}\label{DA1}
The observation functionals \(\{l_j\}_{j=1}^m\) are bounded linear functionals on \(H^1(\Gamma)\).
\end{dataassumption}

A prototypical special case covered by Assumption~\ref{DA1} is pointwise observation, \(l_j(p)=p(x_j)\) for \(x_j\in\Gamma\). Indeed, since \(\Gamma\) is locally one-dimensional, the Sobolev embedding \(H^1(\Gamma)\hookrightarrow C(\Gamma)\subset L^\infty(\Gamma)\) holds, and hence point evaluation \(p\mapsto p(x)\) is a bounded linear functional on \(H^1(\Gamma)\). This is in contrast with higher-dimensional settings, where point evaluation is not continuous on \(H^1\) in general.

\nc
\black Given \(\{l_j\}_{j=1}^m\) as in Data Assumption~\ref{DA1}, we define the \emph{observation map} \nc
\begin{align*}
    \mathcal{Q}: \black H^1 \nc(\Gamma) & \to \mathbb{R}^m, \\
    p &\mapsto \bigl(l_1(p), \dots, l_m(p)\bigr)^\top,
\end{align*}
and we refer to the operator $\G := \mathcal{Q} \circ \mathcal{F}: \black H^1 \nc(\Gamma) \to \R^m$ obtained by composing the observation and forward maps as the \emph{forward model}.
Concatenating the observations, we then have
\begin{equation}\label{eq:data}
   y = \mathcal{G}(u) + \eta, 
\end{equation}
where  $ y = (y_1, \ldots, y_m)^\top$ and $\eta = (\eta_1, \ldots, \eta_m)^\top.$ 
We assume that the measurement noise is centered Gaussian with positive definite covariance $\Sigma,$ written $\eta \sim \mathcal{N}(0, \Sigma).$

Our goal is to recover the unknown parameter \black\(u \in  H^1 (\Gamma)\), \nc from the data \(y \in \R^m\) using the Bayesian approach to inverse problems \cite{kaipio2006statistical,stuart2010inverse,dashti2013bayesian,calvetti2007introduction,sanzstuarttaeb}. To that end, we will
specify a \emph{prior distribution} on $u$ and combine it with the \emph{likelihood function} determined by \eqref{eq:data} to obtain the \emph{posterior distribution} of $u$ given $y.$ In the Bayesian framework, the posterior distribution is used to obtain point estimates for $u$ and quantify the uncertainty in the reconstruction \cite{sanzstuarttaeb}.
However, in our metric graph setting a careful choice of prior distribution will be essential to rigorously justify that the posterior distribution is a well-defined probability measure, which can be written as a change of measure with respect to the prior. We introduce and motivate our choice of prior in Subsection~\ref{ssec:Gaussian Processes on Metric Graphs}. Next, we define the likelihood function in Subsection~\ref{ssec:likelihood}. Finally, we formally derive the posterior distribution in Subsection~\ref{ssec:posterior}. In Section~\ref{sec:well-posedness}, we will rigorously show that the posterior distribution is well defined and stable with respect to perturbations in the data.

\subsection{Prior Distribution} 
\label{ssec:Gaussian Processes on Metric Graphs}
We will take the prior distribution to be the law of a (generalized) Whittle-Matérn Gaussian process on the metric graph \cite{bolin2024gaussian,bolin2024regularity}. 
Specifically, we set the prior to be $\mu_0 = \text{Law}(u),$ where $u$ is the solution of the fractional stochastic PDE 
\begin{equation}
\label{equa:prior gaussian}
     \Bigl(\black \kappa_{\rm ref}^2 \nc - \nabla \cdot (\black a_{\rm ref}\nc \nabla)\Bigr)^\alpha u = \mathcal{W},
\end{equation}
\black
equipped with standard Kirchhoff vertex conditions. 
Here $\kappa_{\rm ref}\in L^\infty(\Gamma)$ satisfies $\operatorname{ess\,inf}_{x\in\Gamma}\kappa_{\rm ref}(x)\ge c>0$ for some $c>0$, and $ a_{\rm ref}:\Gamma\to\R$ is positive and Lipschitz. 
Gaussian white noise $\mathcal W$ on $\Gamma$ is understood as a centered Gaussian random generalized function with covariance given by the $L^2(\Gamma)$ inner product (cf.\ \black Proposition~\ref{prop:forward_under_H1}\nc). 
Equivalently, for any orthonormal basis $\{e_k\}_{k\in\N}$ of $L^2(\Gamma)$, it admits the series representation (in a weak sense)
\[
\mathcal W=\sum_{k\ge1}\xi_k e_k,
\qquad \xi_k\stackrel{\mathrm{i.i.d.}}{\sim}\mathcal N(0,1).
\]
The fractional power $\bigl(\black \kappa_{\rm ref}^2 \nc -\nabla\cdot(\black a_{\rm ref}\nc\nabla)\bigr)^\alpha$ is defined spectrally, in the same way as in Section~\ref{ssec:quantum_graphs}: let $\{(\mu_j,\psi_j)\}_{j\ge1}$ denote the eigenpairs of $\black \kappa_{\rm ref}^2\nc -\nabla\cdot(\black a_{\rm ref}\nc\nabla)$, where $\{\psi_j\}_{j\ge1}$ forms an orthonormal basis of $L^2(\Gamma)$. 
\black Using this eigenbasis, $\mathcal W$ admits the series representation (in a weak sense)
\[
\mathcal W=\sum_{j\ge1}\tilde\xi_j \psi_j,
\qquad \tilde\xi_j\stackrel{\mathrm{i.i.d.}}{\sim}\mathcal N(0,1),
\]
and the prior admits the Karhunen--Lo\`eve representation
\[
u=\sum_{j\ge1}\mu_j^{-\alpha}\,\tilde\xi_j\,\psi_j.
\]
We emphasize that \eqref{equa:prior gaussian} and \eqref{equa:p} correspond to two distinct modeling steps. 
Equation~\ref{equa:prior gaussian} is used solely to \emph{define} the Whittle--Mat\'ern Gaussian prior on $u$, via a fixed reference operator and fixed parameters $({\kappa_{\rm ref}, a_{\rm ref}},\alpha)$. 
By contrast, Equation~\ref{equa:p} defines the physical forward model mapping $u$ to the state $p$, through a fractional elliptic operator whose coefficient depends on $u$ and which, in general, involves different parameters (such as $\kappa$ and $\beta$). 
Thus, while both operators are of Laplace type and adapted to the geometry of $\Gamma$, they play different roles in the Bayesian formulation (prior specification versus data model).
\nc

\begin{proposition}[Regularity of Prior Gaussian Process (See \cite{bolin2024regularity})]
\label{prop:regularity of priori}
Let \(u\) be the solution of \eqref{equa:prior gaussian}, where the operator is equipped with Kirchhoff vertex conditions and the coefficients $\black a_{\rm ref}\nc$ and $\black \kappa_{\rm ref}\nc$ satisfy the assumptions above. Then, \(u \in H^1(\Gamma)\) \(\mathbb{P}\)-a.s.\ if and only if \(\alpha > \tfrac{3}{4}\). Moreover, for any \(\varepsilon>0\) and any \(\frac{1}{4}\leq\alpha\leq 2+\frac{1}{4}\), the sample paths of \(u\) belong to $\widetilde{H}^{2\alpha-\frac{1}{2}-\varepsilon}\black(\Gamma)\nc$ \(\mathbb{P}\)-a.s.
\end{proposition}

\black
\begin{corollary}[Forward Map on Prior Support]
\label{cor:alpha_gt_3over4_forward}
Assume \(\alpha>\tfrac34\) in \eqref{equa:prior gaussian}. Then \(u\in H^1(\Gamma)\) \(\mathbb P\)-a.s.,
and Theorem~\ref{prop:forward_under_H1} applies \(\mathbb P\)-a.s.\ to the random conductivity \(a=\exp(u)\).
In particular, the forward map \(u\mapsto p=\L_u^{-\beta}f\) is well-defined for every \(\beta\ge1\) and \(f\in L^2(\Gamma)\).
\end{corollary}

\begin{proof}
This is an immediate consequence of Propositions~\ref{prop:forward_under_H1} and \ref{prop:regularity of priori}.
\end{proof}
\nc

\subsection{Likelihood Function}\label{ssec:likelihood}
Building on the data generating process described above, we now define the likelihood function and introduce the corresponding potential function. Since the data $y$ take values in $\mathbb{R}^m$ and the noise $\eta$ is modeled as a mean-zero Gaussian random vector with covariance matrix $\Sigma\in\mathbb{R}^{m\times m}$, we define the potential function by
\begin{equation}
\label{equa:Phi}
  \Phi(u;y)\;:=\;\frac{1}{2}\,\bigl(y - \mathcal{G}(u)\bigr)^\top \,\Sigma^{-1}\,\bigl(y - \mathcal{G}(u)\bigr) =: \frac12\|y-\G(u)\|_{\Sigma^{-1}}^{2},
\end{equation}
where, recall, $\mathcal{G}$ denotes the forward model.
Then, the likelihood of $y$ given $u$ can be expressed as
\[
  \pi(y\mid u)\;\propto\;\exp\bigl(-\Phi(u;y)\bigr).
\]

\black
\begin{remark}\label{rem:infdata}
Our formulation can also accommodate \emph{functional} data (see \cite{dashti2011uncertainty,stuart2010inverse}), but we do not pursue this extension here. 
Concretely, one may take the data space \(Y\) to be an infinite-dimensional Banach or Hilbert space and model observations as
\[
y=\mathcal{O}\circ\mathcal{F}(u) +\eta,\qquad y\in Y,
\]
where \(\mathcal{O}\) is a bounded linear observation operator and \(\eta\) is Gaussian noise in \(Y\).
In this setting, well-definedness of the posterior and continuity with respect to perturbations in \(y\) can be studied using the general Banach-space theory of Bayesian inverse problems in \cite{dashti2011uncertainty}, which is formulated for general Banach parameter spaces \(X\) and data spaces \(Y\). In this paper, we restrict attention to finite-dimensional data, as this is the most common setting in the applications we consider.
\end{remark}
\nc

\subsection{Posterior Distribution}\label{ssec:posterior}
Combining via Bayes' formula the Gaussian Whittle–Matérn prior distribution $\mu_0$ introduced in Subsection \ref{ssec:Gaussian Processes on Metric Graphs} with the likelihood function introduced in Subsection~\ref{ssec:likelihood}, we can formally derive the posterior distribution $\mu^y$ of $u$ given $y,$ which can be expressed as 
    \begin{equation}\label{eq:posteriorintro}
   \frac{d\mu^y}{d\mu_0}(u) \propto \exp\bigl(-\Phi(u;y)\bigr). 
\end{equation} 
The fact that the posterior distribution $\mu^y$ is given by a change of measure with respect to the prior as expressed in \eqref{eq:posteriorintro} will be rigorously established in the next section, where we also show the continuous dependence of the posterior on the data in Hellinger distance.

\section{Well-posedness Theory}\label{sec:well-posedness}
This section establishes the well-posedness of the Bayesian inverse problem introduced in Section~\ref{sec:problemformulation}. The main result is Theorem~\ref{thm: posterior distribution def}, which establishes the existence and uniqueness of the posterior distribution, as well as its stability with respect to perturbations in the data. The proof of Theorem~\ref{thm: posterior distribution def} relies on a careful study of the stability of the forward map (Theorem~\ref{thm:boundedness_solution fractional operator}) and the forward model (Theorem~\ref{thm:observations_lipschitz fractional}).

\subsection{Main Results}
\label{sec: main results}
The main technical contribution of this paper is the following stability result for the forward map of elliptic $( \beta = 1)$ and fractional elliptic $(\beta>1)$ inverse problems. 
We recall that $\alpha$ controls the regularity of samples from the prior $\mu_0$ introduced in Subsection~\ref{ssec:Gaussian Processes on Metric Graphs} (cf.\ Proposition~\ref{prop:regularity of priori}), and that $f$ denotes the right-hand side of the differential equation \eqref{equa:p}. 
\black We first state in Theorem~\ref{thm:boundedness_solution fractional operator}  deterministic stability estimates for functions $u,u_1,u_2\in L^\infty(\Gamma)$. 
Then in Corollary~\ref{cor:stability_mu0}, we relate these deterministic estimates to draws from our Gaussian prior $\mu_0$ defined in \eqref{equa:prior gaussian}; notice that the only property we use from $\mu_0$ is that it is supported on $L^\infty(\Gamma),$ so that the theory readily extends to other priors satisfying this condition. \nc
The proof of Theorem~\ref{thm:boundedness_solution fractional operator} is deferred to Subsection~\ref{ssec:stabilityforwardmap}.

\begin{theorem}[Stability of the Forward Map]
\label{thm:boundedness_solution fractional operator}
Let $\beta\ge 1$ be a fixed constant and $f\in L^{2}(\Gamma)$. \black Assume $u,u_1,u_2\in L^\infty(\Gamma)$,\nc and denote by $p=\F(u)$ and $p_i=\F(u_i)$
the weak solutions to \eqref{equa:p} corresponding to $u$ and $u_i$, respectively.
Then there exists a constant
\[
c(\beta)=c\bigl(\beta,\|f\|_{L^{2}(\Gamma)},\Gamma\bigr),
\]
decreasing with respect to $\beta$, such that
\begin{align}
        \|\F(u)\|_V 
        &\leq c(\beta)\exp\,\!\bigl(\beta\|u\|_{L^\infty(\Gamma)}\bigr), 
        \label{eq:thmbound1}\\[1mm]
        \|\F(u_1)-\F(u_2)\|_V 
        &\leq c(\beta)\exp\,\!\Bigl((\beta+2)\max\{\|u_1\|_{L^\infty(\Gamma)},\|u_2\|_{L^\infty(\Gamma)}\}\Bigr)
        \|u_1-u_2\|_{L^\infty(\Gamma)}.
        \label{eq:thmbound2}
\end{align}
Here $\|\cdot\|_V=\|\cdot\|_{H^1(\Gamma)}$ or $\|\cdot\|_V=\|\cdot\|_{L^\infty(\Gamma)}$.
\end{theorem}

\black
\begin{corollary}[Application to the Gaussian Prior $\mu_0$]
\label{cor:stability_mu0}
Assume $\alpha>3/4$ in \eqref{equa:prior gaussian}. Then $u\in L^\infty(\Gamma)$ $\mu_0$-almost surely, and hence
the bounds \eqref{eq:thmbound1}--\eqref{eq:thmbound2} hold for $\mu_0$-almost all $u,u_1,u_2\sim\mu_0$.
\end{corollary}

\begin{proof}
By Proposition~\ref{prop:regularity of priori}, for any $\alpha>3/4$ we have 
$u\in H^1(\Gamma)$ $\mu_0$-almost surely, and hence $u\in L^\infty(\Gamma)$ $\mu_0$-almost surely, by Sobolev embedding on compact metric graphs.
The conclusion follows by applying Theorem~\ref{thm:boundedness_solution fractional operator}.
\end{proof}
\nc

\begin{remark}\label{rem:Linfty-stability}
Since compact metric graphs are locally one-dimensional, it suffices to establish the bounds \eqref{eq:thmbound1} and \eqref{eq:thmbound2} for the $H^1(\Gamma)$ norm. To see this, notice that each edge \(e\in E\) is isometric to a real interval, so the classical one‑dimensional Sobolev embedding
\[
  H^1(e)\;\hookrightarrow\;C(e),
  \qquad
  \|p\|_{L^\infty(e)} \le c(e)\,\|p\|_{H^1(e)}
\]
applies edge‑wise.  Enforcing continuity at the vertices then yields the global embedding
\[
  H^1(\Gamma)\;\hookrightarrow\;L^\infty(\Gamma),
  \qquad
  \|p\|_{L^\infty(\Gamma)}\le c(\Gamma)\,\|p\|_{H^1(\Gamma)}.
\]
As a result, bounds in the \(H^1(\Gamma)\) norm transfer immediately to the  \(L^\infty(\Gamma)\) norm.  This control of the $L^\infty(\Gamma)$ norm by the $H^1(\Gamma)$ norm holds even when the metric graph \(\Gamma\) is embedded in \(\mathbb R^d\), $d \ge 2,$ while the conclusion fails in high-dimensional Euclidean domains.
\end{remark}

\black
\begin{remark}[Implications and Limitations of  Theorem~\ref{thm:boundedness_solution fractional operator}]
\label{rmk:implications and limitations}
Theorem~\ref{thm:boundedness_solution fractional operator} ensures stability of the forward map $u\mapsto \F(u)$ and is a key input for data-to-posterior well-posedness. 
It does not by itself imply identifiability or accurate reconstruction of $u$ from finitely many noisy observations, which depend on the observation regime. 
The estimate \eqref{eq:thmbound2} is the metric-graph analogue of classical stability bounds for elliptic inverse problems in Euclidean domains; see, e.g., Proposition~3.4 in \cite{dashti2011uncertainty} and related refinements in \cite{Vollmer_2013,nickl2019convergenceratespenalisedsquares}.
In the Euclidean setting, such stability estimates have been used as key ingredients in posterior consistency analyses (cf.\ \cite{Vollmer_2013}) and in establishing convergence guarantees for MAP-type estimators (cf.\ \cite{nickl2019convergenceratespenalisedsquares}).
We view Theorem~\ref{thm:boundedness_solution fractional operator} as a natural starting point for analogous consistency analyses on compact metric graphs, which we leave to future work.
\end{remark}

\nc

The stability of the forward map, along with our assumptions on the observation map encoded in \black Data Assumption \ref{DA1} imply \nc the following stability result for the forward model. 
\begin{theorem}[Stability of the Forward Model]
\label{thm:observations_lipschitz fractional}
In the setting of Theorem~\ref{thm:boundedness_solution fractional operator} and \black under Data Assumption \ref{DA1}\nc, there exists a constant 
\[
c(\beta) = c\bigl(\beta, \|f\|_{L^{2}(\Gamma)}, \Gamma, \max_j \|l_j\|_{H^{-1}(\Gamma)}\bigr),
\]
decreasing with respect to $\beta$,
 such that, for almost all $u, u_1, u_2 \sim \mu_0$:
\begin{align}
        \|\mathcal{G}(u)\| &\leq c(\beta) \exp\bigl(\beta\|u\|_{L^\infty(\Gamma)}\bigr),\label{eq:thmGbound1}\\[1mm]
        \|\mathcal{G}(u_1) - \mathcal{G}(u_2)\| &\leq c(\beta) \exp\Bigl((\beta+2)\, \max\bigl\{\|u_1\|_{L^\infty(\Gamma)}, \|u_2\|_{L^\infty(\Gamma)}\bigr\}\Bigr) \|u_1 - u_2\|_{L^\infty(\Gamma)}.\label{eq:thmGbound2}
\end{align}
\end{theorem}
\begin{proof}
By Theorem~\ref{thm:boundedness_solution fractional operator}, there exists a constant \( \tilde{c}(\beta) \), independent of \( u \), such that:
\begin{equation}
\label{equa:difference}
\left\{
    \begin{aligned}
    &\|p\|_{H^1(\Gamma)}\leq \tilde{c}(\beta)\exp(\beta \|u\|_{L^\infty(\Gamma)}),\\
    &\|p_1 - p_2\|_{H^1(\Gamma)} \leq \tilde{c}(\beta)\, \exp\Bigl((\beta+2) \max\bigl\{\|u_1\|_{L^\infty(\Gamma)}, \|u_2\|_{L^\infty(\Gamma)}\bigr\}\Bigr)\, \|u_1 - u_2\|_{L^\infty(\Gamma)},
\end{aligned}
\right.
\end{equation}
where $p_1$ and $p_2$ denote the solutions to \eqref{equa:p} corresponding to $u_1$ and $u_2$, respectively.
 
Since each $l_j$ is a bounded linear functional on $H^1(\Gamma)$, the Riesz representation theorem implies that
\( 
    |l_j(p)| \leq \|l_j\|_{H^{-1}}\, \|p\|_{H^1}
\) 
for $j=1,\ldots,m.$
Applying the observation map to the solution $p$ and the difference $p_1 - p_2$ respectively, we obtain that 
\begin{equation}
\left\{
\begin{aligned}
    &\|\mathcal{G}(u)\|\leq \|l\|_{H^{-1}(\Gamma)}\|p\|_{H^1{(\Gamma)}},\\
    &\|\mathcal{G}(u_1) - \mathcal{G}(u_2)\| = \sqrt{\sum_{j=1}^m |l_j(p_1 - p_2)|^2} \leq \sqrt{\sum_{j=1}^m \|l_j\|_{H^{-1}}^2}\, \|p_1 - p_2\|_{H^1(\Gamma)}.
\end{aligned}
\right.
\end{equation}

Hence, it follows from \eqref{equa:difference} that there exists a constant  \( c(\beta):=c(\beta,\|f\|_{L^2(\Gamma)},\Gamma,\underset{j}{\max} \|l_j\|_{H^{-1}(\Gamma)})\) such that:
\begin{equation*}
\left\{
\begin{aligned}
     & \|\mathcal{G}(u)\|_{H^1(\Gamma)} \leq c(\beta)\exp(\beta \|u\|_{L^\infty(\Gamma)}),\\
     &\|\mathcal{G}(u_1) - \mathcal{G}(u_2)\| \leq c(\beta)\, \exp\Bigl((\beta+2) \max\bigl\{\|u_1\|_{L^\infty(\Gamma)}, \|u_2\|_{L^\infty(\Gamma)}\bigr\}\Bigr)\, \|u_1 - u_2\|_{L^\infty(\Gamma)},\\
\end{aligned}
\right.
\end{equation*}
establishing the boundedness and Lipschitz continuity of the forward model. 
\end{proof}

We will use the Hellinger distance to quantify the proximity between posterior distributions $\mu^y$ and $\mu^{y'}$ corresponding to different realizations \(y\) and \(y'\) of the data.

\begin{definition}[Hellinger Distance]
\label{def:hellinger_distance}
For two probability measures \(\nu_1\) and \(\nu_2\) on the same measurable space \((X, \mathcal{B})\), their \emph{Hellinger distance} is defined as
\[
d_H(\nu_1, \nu_2)
  \;:=\;
  \left(\frac{1}{2}\int_X \biggl(\sqrt{\frac{d\nu_1}{d\nu}} - \sqrt{\frac{d\nu_2}{d\nu}}\biggr)^2\,d\nu\right)^{1/2},
\]
where \(\nu\) is a reference measure that dominates both \(\nu_1\) and \(\nu_2\).
\end{definition}

The following is the main result of this section.
\begin{theorem}[Well‑posedness of the Bayesian Inverse Problem]
\label{thm: posterior distribution def}

Let \(\mu_{0}\) be the Gaussian Whittle–Matérn prior defined in
Subsection~\ref{ssec:Gaussian Processes on Metric Graphs} with regularity
parameter \(\alpha>1\).
For every \(\beta\ge 1\) and measurement \(y\in \R^{m}\), define the
posterior probability measure \(\mu^{y}\) on \black \(H^{1}(\Gamma)\) \nc by
\begin{equation}\label{equa:posterior distribution}
  \frac{d\mu^{y}}{d\mu_{0}}(u)
  \;=\;
  \frac{1}{Z(y)}\,
  \exp\bigl(-\Phi(u;y)\bigr),
  \qquad 
  Z(y):=\int_{\black H^1 \nc(\Gamma)}\exp\bigl(-\Phi(u;y)\bigr)\,d\mu_{0}(u),
\end{equation}
where the potential \(\Phi\colon \black H^1 \nc(\Gamma) \times \R^{m} \to\R\) is given by \eqref{equa:Phi}. Then, the following holds:

\medskip
\noindent
\textbf{(i) Well‑definedness.}
\(\mu^{y}\) is a well‑defined probability measure on \(\black H^1 \nc(\Gamma)\).

\smallskip
\noindent
\textbf{(ii) Stability with respect to the data.}
For every radius \(r>0\), there exists a constant \(c=c(r)>0\) such that,  
for all \(y,y'\in \R^{m}\) with \(\max \bigl\{\|y\|,\|y'\| \bigr\}\le r\),
\begin{equation}\label{equa: lipschitz of mu^y}
  d_{H} \bigl(\mu^{y},\mu^{y'}\bigr)
  \;\le\;
  c(r)\,\|y-y'\| .
\end{equation}
\end{theorem}

\begin{proof}
\label{proof: thm posterior}
 According to \cite{dashti2011uncertainty}, if the potential $\Phi(u,y)$ satisfies the following properties: 
\begin{itemize}
    \item[(I)] \textit{For every $\epsilon > 0$ and $r > 0$, there exists $M = M(\epsilon, r) \in \mathbb{R}$ such that, for all $u \in \black H^1 \nc(\Gamma)$ and for all $y \in \mathbb{R}^m$ with $\|y\| < r$,}
    \[
    \Phi(u,y) \geq M - \epsilon \|u\|^2_{\black H^1 \nc(\Gamma)};
    \]
    
    \item[(II)] \textit{For every $r > 0$, there exists $K = K(r) > 0$ such that, for all $u \in \black H^1 \nc(\Gamma)$ and for all $y \in \mathbb{R}^m$ with $\max
    \bigl\{ \|u\|_{\black H^1 \nc(\Gamma)}, \|y\| \bigr\} < r$,}
    \[
    \Phi(u,y) \leq K;
    \]
    
    \item[(III)] \textit{For every $r > 0$, there exists $L = L(r) > 0$ such that, for all $u_1, u_2 \in \black H^1 \nc(\Gamma)$ and for all $y \in \mathbb{R}^m$ with $\max
    \bigl\{\|u_1\|_{\black H^1 \nc(\Gamma)}, \|u_2\|_{\black H^1 \nc(\Gamma)}, \|y\| \bigr\} < r$,}
    \[
    |\Phi(u_1,y) - \Phi(u_2,y)| \leq L \|u_1 - u_2\|_{\black H^1 \nc(\Gamma)};
    \]
    
    \item[(IV)] \textit{For every $\epsilon > 0$ and $r > 0$, there exists $c = c(\epsilon, r) \in \mathbb{R}$ such that, for all $y_1, y_2 \in \mathbb{R}^m$ with $\max
    \bigl\{\|y_1\|, \|y_2\| \bigr\} < r$ and for every $u \in \black H^1 \nc(\Gamma)$,}
    \[
    |\Phi(u,y_1) - \Phi(u,y_2)| \leq \exp\Bigl(\epsilon \|u\|_{\black H^1 \nc(\Gamma)}^2 + c\Bigr) \,  \|y_1 - y_2\|,
    \]
\end{itemize}
and if \(\mu_0\big(\black H^1 \nc(\Gamma)\big)=1\), then the posterior measure $\mu^y$ defined in \eqref{equa:posterior distribution} is well defined and depends Lipschitz continuously on the data $y$, as expressed in \eqref{equa: lipschitz of mu^y}. We now verify that these properties hold in our setting. Throughout the proof, we regard the metric graph \(\Gamma\), the covariance matrix~\(\Sigma\) in \eqref{eq:data}, and the parameter~\(\beta\) in \eqref{equa:p}  as fixed. Accordingly, we allow the constants to implicitly depend on these parameters.

 First, since \(\Phi(u,y)\ge0\) for all \(u,y\), we may simply set $M(\epsilon,r)=0,$ to verify property~(I). Next, 
 \black under Data Assumption \ref{DA1}, where $\G(u)= \mathcal{Q}(p)=\bigl(l_1(p),\ldots,l_m(p)\bigr)^\top$ with each $l_j$ being a bounded linear functional, \nc
 the following inequality holds for any $(u,y)\in \black H^1 \nc(\Gamma) \times \R^m$ satisfying $\max \bigl\{\|u\|_{\black H^1 \nc(\Gamma)}, \|y\| \bigr\} < r:$
\begin{align}
  \Phi(u,y)
  &= \tfrac12\|y-\G(u)\|_{\Sigma^{-1}}^{2}\notag
  \le \tfrac12\bigl(\|y\|_{\Sigma^{-1}}+\|\G(u)\|_{\Sigma^{-1}}\bigr)^{2}\notag\\
  &\le \tfrac12\bigl(\|\Sigma^{-1}\|_{\mathrm{op}}\,r
          +\|\Sigma^{-1}\|_{\mathrm{op}}\|\G(u)\|\bigr)^{2}\notag\\
  &\overset{\hspace{-0.1cm}\text{(i)}}{\le}\;
      \tfrac12\Bigl(\|\Sigma^{-1}\|_{\mathrm{op}}\,r
          +\|\Sigma^{-1}\|_{\mathrm{op}}\,C_1
            e^{\beta\|u\|_{L^{\infty}(\Gamma)}}\Bigr)^{2}\notag\\
  &\overset{\hspace{-0.1cm}\text{(ii)}}{\le}\;
      \tfrac12\Bigl(\|\Sigma^{-1}\|_{\mathrm{op}}\,r
          + C_1\,\|\Sigma^{-1}\|_{\rm op}\,
            e^{\beta C_2\|u\|_{\black H^1 \nc(\Gamma)}}\Bigr)^{2}\notag\\
  &\le \tfrac12\Bigl(\|\Sigma^{-1}\|_{\mathrm{op}}\,r
          +C_1\,\|\Sigma^{-1}\|_{\rm op}\,e^{\beta C_2 r}\Bigr)^{2}
       \;=:\;K(r).
  \label{eq:Phi-upper-bound}
\end{align}
Here, inequality (i) follows from
Theorem~\ref{thm:observations_lipschitz fractional},
while inequality (ii) uses the Sobolev embedding
$\black H^1 \nc(\Gamma)\hookrightarrow L^{\infty}(\Gamma)$.
Hence, we have shown that property~(II) holds. 

We now seek to verify property ~(III).  Let $y\in\mathbb{R}^{m}$ be fixed and given. For $u\in \black H^1 \nc(\Gamma)$, we denote
\[
  \black \rho(u)\nc \;:=\;\bigl\|y-\mathcal{G}(u)\bigr\|_{\Sigma^{-1}} .
\]
Then, for any given $u_{1},u_{2}\in \black H^1 \nc(\Gamma)$ the elementary identity
$a^{2}-b^{2}=(a-b)(a+b)$ implies
\begin{equation}\label{eq:phi:diff}
\begin{aligned}
  \bigl|\Phi(u_{1},y)-\Phi(u_{2},y)\bigr|
  &=\frac12\,\bigl|\rho(u_{1})^{2}-\rho(u_{2})^{2}\bigr|
   =\frac12\,\bigl|\rho(u_{1})-\rho(u_{2})\bigr|\,  \bigl(\rho(u_{1})+\rho(u_{2})\bigr)\\[4pt]
  &\le \frac12\,\bigl\|\mathcal{G}(u_{1})-\mathcal{G}(u_{2})\bigr\|_{\Sigma^{-1}}
  \black \,\bigl(\rho(u_{1})+\rho(u_{2})\bigr) \nc\\
  &\le \frac12\,\bigl\|\Sigma^{-1}\bigr\|_{\rm op}\,
  \bigl\|\mathcal{G}(u_{1})-\mathcal{G}(u_{2})\bigr\|
  \,\bigl(\rho(u_{1})+\rho(u_{2})\bigr),
\end{aligned}
\end{equation}
where we have used the reverse triangle inequality.
To verify property~(III), suppose that $\max \bigl\{\|u_{i}\|_{\black H^1 \nc(\Gamma)},\|y\| \bigr\}\le r$ $(i=1,2).$ 
According to property~(II), there exists a constant
$c(r)>0$ such that
\begin{equation}\label{eq:phi:sum}
 \black \rho(u_i)\nc =\|y-\mathcal{G}(u_{i})\|_{\Sigma^{-1}}\;\le\;c(r),
  \qquad (i=1,2).
\end{equation}
Applying
Theorem~\ref{thm:observations_lipschitz fractional} yields a constant \(L_1(r)>0\) such that
\begin{equation}\label{eq:phi:diffTerm}
\begin{aligned}
    \bigl\|\mathcal{G}(u_{1})-\mathcal{G}(u_{2})\bigr\|
  \;&\overset{\hspace{-0.1cm}\text{(iii)}}{\le}\; C_3 \exp\Big((\beta+2) \|u_1-u_2\|_{L^\infty(\Gamma)}\Big)\\
  \;&\overset{\hspace{-0.1cm}\text{(iv)}}{\le}\; C_3\exp\big((\beta+2) C_4\|u_1-u_2\|_{\black H^1 \nc(\Gamma)}\big) \\
  \;& \leq\; L_1(r)\,\bigl\|u_{1}-u_{2}\bigr\|_{\black H^1 \nc(\Gamma)}, 
\end{aligned}
\end{equation}
where (iii) follows from Theorem~\ref{thm:observations_lipschitz fractional}, and (iv) can be obtained by applying the Sobolev embedding theorem on $L^\infty(\Gamma)$. Inserting~\eqref{eq:phi:sum} and~\eqref{eq:phi:diffTerm}
into~\eqref{eq:phi:diff} gives
\begin{align*}
     \bigl|\Phi(u_{1},y)-\Phi(u_{2},y)\bigr|
  \;\le\;
  \|\Sigma^{-1}\|_{\rm op}\,c(r)\,L_1(r)\,
  \bigl\|u_{1}-u_{2}\bigr\|_{\black H^1 \nc(\Gamma)}
  \;=:\;L(r)\,\bigl\|u_{1}-u_{2}\bigr\|_{\black H^1 \nc(\Gamma)},
\end{align*}
where \(L(r):=c(r)L_{1}(r)\), implying that property~(III) holds.

To verify property~(IV), fix \(u\in \black H^1 \nc(\Gamma)\) and let \(y_{1},y_{2}\in\R^{m}\) satisfy \(\max
\bigl\{\|y_{1}\|,\|y_{2}\|\bigr\}\le r\).  A direct expansion gives
\begin{equation}\label{eq:phi:y:diff}
\begin{aligned}
  \bigl|\Phi(u,y_{1})-\Phi(u,y_{2})\bigr|
  &= \tfrac12
     \Bigl|(y_{1}-y_{2})^{\top}\Sigma^{-1}\bigl(y_{1}+y_{2}-2\G(u)\bigr)\Bigr| \\[4pt]
  &\le \|y_{1}-y_{2}\|\;\|\Sigma^{-1}\|_{\mathrm{op}}
        \Bigl(\tfrac12\|y_{1}+y_{2}\| + \|\G(u)\|\Bigr) \\[4pt]
  &\le \|y_{1}-y_{2}\|\;\|\Sigma^{-1}\|_{\mathrm{op}}
        \bigl(r + \|\G(u)\|\bigr) \\[4pt]
  &\le \|y_{1}-y_{2}\|\;\|\Sigma^{-1}\|_{\mathrm{op}}
        \bigl(r + A\,e^{\beta\|u\|_{L^\infty(\Gamma)}}\bigr),
\end{aligned}
\end{equation}
where in the third line we used that \(\tfrac12\|y_{1}+y_{2}\|\le r\) by hypothesis, and in the last line we applied Theorem~\ref{thm:observations_lipschitz fractional} to bound \(\|\G(u)\|\le A\,e^{\beta\|u\|_{L^\infty(\Gamma)}}\).

For any given $\varepsilon>0$,  the inequality
\(
  e^{t s}\le 
  \exp\bigl(\varepsilon s^{2}+t^{2}/(4\varepsilon)\bigr)
\)
holds for any $s \ge 0$ and $t>0.$  Hence, by the Sobolev embedding theorem, 
\[
  \|\mathcal G(u)\|
  \;\le\;
  A\,\exp \Bigl(
     \beta C_5\|u\|_{\black H^1 \nc(\Gamma)}
    \Bigr)
  \;\le\;
  A\,\exp \Bigl(
      \varepsilon\|u\|_{\black H^1 \nc(\Gamma)}^{2}
      +\tfrac{(\beta C_5)^2}{4\varepsilon}
    \Bigr).
\]
Inserting the previous estimates into~\eqref{eq:phi:y:diff}, we obtain
\[
\begin{aligned}
  \bigl|\Phi(u,y_{1})-\Phi(u,y_{2})\bigr|
  &\le
    \|y_{1}-y_{2}\|\;\|\Sigma^{-1}\|_{\mathrm{op}}\,
    \Bigl(
      r
      +A\,\exp \bigl(\varepsilon\|u\|_{\black H^1 \nc(\Gamma)}^{2}
                      +\tfrac{(\beta C_5)^2}{4\varepsilon}\bigr)
    \Bigr)                                           \\[6pt]
  &\le
    \exp\Bigl(\varepsilon\|u\|_{\black H^1 \nc(\Gamma)}^{2}+c(\varepsilon,r)\Bigr)\,
    \|y_{1}-y_{2}\|,
\end{aligned}
\]
where the constant 
\begin{align*}
    c(\varepsilon,r)=\log\Bigl(
      \|\Sigma^{-1}\|_{\mathrm{op}}\,\bigl(r+A\,e^{(C_5\beta)^{2}/(4\varepsilon)} \bigr)
    \Bigr)
\end{align*}
depends only on $(\varepsilon,r)$.  
Thus, properties (I)–(IV) are satisfied in our setting. Moreover, it is straightforward to verify that \(\mu_0\big(\black H^1 \nc(\Gamma)\big) = 1\), which directly follows from Proposition~\ref{prop:regularity of priori}. Therefore, by applying Theorems 2.2 and 2.3 from \cite{dashti2011uncertainty}, the proof is complete. 
\end{proof}

\subsection{Stability of the Forward Map}\label{ssec:stabilityforwardmap}
This subsection contains the proof of Theorem~\ref{thm:boundedness_solution fractional operator}. In Subsection~\ref{ssec:elliptic} we study the case $\beta = 1,$ where the forward map concerns an elliptic problem.
Next, in Subsection~\ref{ssec:fractional} we establish the result for $\beta>1,$ where the forward map concerns a fractional elliptic problem. The proof for the fractional elliptic problem makes use of the result for the elliptic problem.

\subsubsection{Elliptic Problem}\label{ssec:elliptic}

The proof of Theorem~\ref{thm:boundedness_solution fractional operator} in the case $\beta = 1$ follows a similar structure to that in \cite{dashti2011uncertainty}, which shows stability of the forward map of an elliptic problem in Euclidean space. Some technical modifications are required to carry out the analysis in our metric graph setting. 

\begin{proof}[Proof of Theorem~\ref{thm:boundedness_solution fractional operator} for $\beta =1$]
    \black Fix $u\in L^\infty(\Gamma)$. \nc
    Let $p:=\F(u)$ be the \black weak \nc solution to \eqref{equa:p} with $\beta =1$ and input parameter $u$.
    \black Testing the weak formulation of \eqref{equa:p} with $\varphi=p$ gives \nc
    \begin{equation}\label{eq:auxeqproof}
        \int_\Gamma \kappa^2 p^2\,{\rm dx} + \int_\Gamma e^u\,|\nabla p|^2\,{\rm dx} = \int_\Gamma fp\,{\rm dx}.
    \end{equation}
    \black This identity can be viewed as the standard energy balance obtained from edgewise integration by parts, where the vertex boundary terms vanish due to the Kirchhoff condition in the weak sense. \nc
    On the other hand, we can bound the right-hand side of  \eqref{eq:auxeqproof} using  Cauchy-Schwarz inequality:
    \begin{align}
    \label{equa:rhs}
        \int_\Gamma f p \,{\rm dx} \leq \|f\|_{L^{2}(\Gamma)} \|p\|_{H^1(\Gamma)}.
    \end{align}
    Combining \eqref{eq:auxeqproof} and \eqref{equa:rhs}, and using $e^u\ge e^{-\|u\|_{L^\infty(\Gamma)}}$, we deduce that
    \begin{align*}
       \min\{\kappa^2, e^{-\|u\|_{L^\infty(\Gamma)}}\}\|p\|_{H^1(\Gamma)}^2 
       &\leq \kappa^2\|p\|_{L^2(\Gamma)}^2+e^{-\|u\|_{L^\infty(\Gamma)}}\|\nabla p\|_{L^2(\Gamma)}^2\\
       &\leq \kappa^2\|p\|_{L^2(\Gamma)}^2+\int_{\Gamma} e^u |\nabla p|^2\,{\rm dx}
       \leq \int_{\Gamma}fp\,{\rm dx} \leq  \|f\|_{L^{2}(\Gamma)} \|p\|_{H^1(\Gamma)},
    \end{align*}
    which implies that
    \begin{align}
        \label{equa:boundedness}
        \|p\|_{H^1(\Gamma)}
        &\leq \frac{1}{\min\{\kappa^2,e^{-\|u\|_{L^\infty(\Gamma)}}\}} \|f\|_{L^{2}(\Gamma)} 
        \leq \frac{e^{\|u\|_{L^\infty(\Gamma)}}}{\min\{\kappa^2,1\}} \|f\|_{L^{2}(\Gamma)} 
        \leq C\,e^{\|u\|_{L^\infty(\Gamma)}}\|f\|_{L^{2}(\Gamma)},
    \end{align}
    where $C$ is a constant independent of $u$. This completes the proof of the bound \eqref{eq:thmbound1} in the $H^1(\Gamma)$ norm for $\beta =1$. The bound in the $L^\infty(\Gamma)$ norm then follows from Remark~\ref{rem:Linfty-stability}.

    Next, we prove the bound \eqref{eq:thmbound2}. Let $p_1:= \F(u_1)$ and $p_2:= \F(u_2)$ be the solutions to \eqref{equa:p} with $\beta = 1$ and parameters $u_1$ and $u_2$, respectively. 
    \black We use the weak formulation: for $i=1,2$ and all $\varphi\in H^1(\Gamma)$,
    \[
      \int_\Gamma \kappa^2 p_i\,\varphi\,{\rm dx}+\int_\Gamma \exp(u_i)\,\nabla p_i\cdot \nabla\varphi\,{\rm dx}
      =\int_\Gamma f\,\varphi\,{\rm dx}.
    \]
    \nc
    Defining $q := p_1 - p_2$ and subtracting the identities for $i=1$ and $i=2$, we obtain \black for all $\varphi\in H^1(\Gamma)$ \nc
    \[
    \int_\Gamma \kappa^2 q\,\varphi\,{\rm dx}+\int_\Gamma \exp(u_1)\,\nabla q\cdot \nabla\varphi\,{\rm dx}
    =
    \int_\Gamma \bigl(\exp(u_2)-\exp(u_1)\bigr)\,\nabla p_2\cdot \nabla\varphi\,{\rm dx}.
    \]
    \black Taking $\varphi=q$ yields \nc
    \begin{equation}
    \label{equa:simplification 2}
        \kappa^2\int_\Gamma q^2 \,{\rm dx}+\int_\Gamma \exp(u_1)\lvert\nabla q\rvert^2\,{\rm dx}
        =
        \int_\Gamma \bigl(\exp(u_2)-\exp(u_1)\bigr)\nabla p_2\cdot \nabla q\,{\rm dx}.
    \end{equation}
    \black Since $u_1,u_2\in L^\infty(\Gamma)$, the exponential function is Lipschitz on the range of $(u_1,u_2)$, and hence \nc
    \begin{equation*}
        \lvert\exp(u_1)-\exp(u_2)\rvert \leq \| u_1-u_2\|_{L^\infty(\Gamma)} \cdot \exp\Bigl(\max \bigl\{\|u_1\|_{L^\infty(\Gamma)},\|u_2\|_{L^\infty(\Gamma)}\bigr\}\Bigr).
    \end{equation*}
    Using this bound, we deduce that
    \begin{align*}
        &\kappa^2 \|q\|^2_{L^2(\Gamma)}+  e^{-\|u_1\|_{L^\infty(\Gamma)}} \|\nabla q\|^2_{L^2(\Gamma)}
        \leq \ \kappa^2\int_\Gamma q^2 \,{\rm dx}+\int_\Gamma \exp(u_1)\lvert\nabla q\rvert^2\,{\rm dx}\\
        &\leq \ \left\lvert \int_{\Gamma} \bigl(\exp(u_2)-\exp(u_1)\bigr)\nabla p_2\cdot \nabla q\,{\rm dx}\right \rvert    \\
        &\leq \ \exp\Bigl(\max\{\|u_1\|_{L^\infty(\Gamma)},\|u_2\|_{L^\infty(\Gamma)}\}\Bigr)\|u_1-u_2\|_{L^\infty(\Gamma)}\|\nabla p_2\|_{L^2(\Gamma)}\|\nabla q\|_{L^2(\Gamma)} \\ 
        &  \le \frac{1}{2} \exp\bigl(\|u_1\|_{L^\infty(\Gamma)}\bigr)   \exp\Bigl(2\max\bigl\{\|u_1\|_{L^\infty(\Gamma)},\|u_2\|_{L^\infty(\Gamma)}\bigr\}\Bigr) \|u_1-u_2\|_{L^\infty(\Gamma)}^2  \|\nabla p_2\|_{L^2(\Gamma)}^2
        \\ 
         &\hspace{2cm}+\frac12 \exp\bigl(-\|u_1\|_{L^\infty(\Gamma)}\bigr) \|\nabla q\|^2_{L^2(\Gamma)}, \\
        &  \le \frac{1}{2} \exp\Bigl(3\max\bigl\{\|u_1\|_{L^\infty(\Gamma)},\|u_2\|_{L^\infty(\Gamma)}\bigr\}\Bigr) \|u_1-u_2\|_{L^\infty(\Gamma)}^2  \|\nabla p_2\|_{L^2(\Gamma)}^2  +\frac{1}{2} e^{-\|u_1\|_{L^\infty(\Gamma)}} \|\nabla q\|^2_{L^2(\Gamma)},
     \end{align*}
     where in the second-to-last step we used Young's inequality $ ab \le \frac{1}{2\varepsilon} a^2 + \frac{\varepsilon}{2} b^2$  
 with
    $a = \exp\Bigl(\max \bigl\{\|u_1\|_{L^\infty(\Gamma)},\|u_2\|_{L^\infty(\Gamma)}\bigr\}\Bigr)\|u_1-u_2\|_{L^\infty(\Gamma)}\|\nabla p_2\|_{L^2(\Gamma)},$  $
    b = \|\nabla q\|_{L^2(\Gamma)}, $ and
    $\varepsilon = e^{-\|u_1\|_{L^\infty(\Gamma)}}.$
Therefore, since
$
\|q\|_{H^1(\Gamma)}^2 = \|q\|_{L^2(\Gamma)}^2 + \|\nabla q\|_{L^2(\Gamma)}^2,
$
we obtain that
\begin{align*}
    &\|q\|_{H^1(\Gamma)}^2  \\
    & \leq \frac{1}{2\min\{\kappa^2,\frac{1}{2}e^{-\|u_1\|_{L^\infty(\Gamma)}}\}} \exp\Bigl(3 \max \bigl\{ \|u_1\|_{L^\infty(\Gamma)}, \|u_2\|_{L^\infty(\Gamma)} \bigr\} \Bigr)\|u_1 - u_2\|_{L^\infty(\Gamma)}^2 \|\nabla p_2\|_{L^2(\Gamma)}^2\\
    & \leq \frac{e^{\|u_1\|_{L^\infty(\Gamma)}}}{\min\{2\kappa^2, 1\}}\exp\Bigl(3 \max \bigl\{ \|u_1\|_{L^\infty(\Gamma)}, \|u_2\|_{L^\infty(\Gamma)} \bigr\} \Bigr)\|u_1 - u_2\|_{L^\infty(\Gamma)}^2 \|\nabla p_2\|_{L^2(\Gamma)}^2\\
    & \leq \frac{1}{\min \bigl\{2\kappa^2, 1 \bigr\}}\exp\Bigl(4 \max\{ \|u_1\|_{L^\infty(\Gamma)}, \|u_2\|_{L^\infty(\Gamma)} \} \Bigr)\|u_1 - u_2\|_{L^\infty(\Gamma)}^2 \|\nabla p_2\|_{L^2(\Gamma)}^2.
\end{align*}

Finally, applying the estimate from \eqref{equa:boundedness} to control $\|\nabla p_2\|_{L^2(\Gamma)}$ in terms of $\|f\|_{L^{2}(\Gamma)}$ and $\|u_2\|_{L^\infty(\Gamma)}$, we conclude that
\begin{equation*}
    \|q\|_{H^1(\Gamma)} 
    \leq C\, \|f\|_{L^{2}(\Gamma)}\, 
    \exp\Bigl( 3 \max \bigl\{ \|u_1\|_{L^\infty(\Gamma)}, \|u_2\|_{L^\infty(\Gamma)} \bigr\} \Bigr) 
    \|u_1 - u_2\|_{L^\infty(\Gamma)},
\end{equation*}
for a constant $C>0$ that depends only on $\kappa$ and $\Gamma$. We have hence established the bound \eqref{eq:thmbound2} in the $H^1(\Gamma)$ norm. The bound in the $L^\infty(\Gamma)$ norm then follows from Remark~\ref{rem:Linfty-stability}, concluding the proof of Theorem~\ref{thm:boundedness_solution fractional operator} for $\beta = 1$.
\end{proof}

\subsubsection{Fractional Elliptic Problem}\label{ssec:fractional}
Here, we prove Theorem~\ref{thm:boundedness_solution fractional operator} in the case $\beta>1$. The proof technique that we used in the case $\beta =1$ ---based on a weak formulation--- cannot be applied directly, since now \(\L_u^\beta\) is defined spectrally. Instead, we will leverage a series representation of the solution, along with  Weyl's law for \(\L_u\) and the stability of the elliptic forward map proved in the previous subsection.

We first introduce a new version of the Weyl's law on metric graphs established in \cite{bolin2024regularity}, in which we ensure that the constants are independent of the input parameter $u.$

\begin{lemma}\label{lem:weyl's law}
    Let $\L_u$ be the second-order elliptic operator defined in \eqref{equa:p} with $u \in \black H^1 \nc(\Gamma)$. Then, there exist constants $C_1, C_2 > 0$ independent of $u$ such that, for all $j \in \mathbb{N},$
    \begin{equation*}
         C_1\exp \bigl(-\|u\|_{L^\infty(\Gamma)} \bigr)\, j^2 \leq \lambda_j \leq C_2\exp \bigl(\|u\|_{L^\infty(\Gamma)}\bigr) \, j^2,
    \end{equation*}
    where $\{\lambda_j\}_{j \in \mathbb{N}}$ are the eigenvalues of $\L_u$ arranged in non-decreasing order.
\end{lemma}
\begin{proof}
 To clarify the dependence of the constants on $u$ using the classical Weyl's law \cite{MR3894619}, we denote $M := \|u\|_{L^\infty(\Gamma)}$ and consider two constant‑coefficient elliptic operators as extreme cases:
\begin{equation*}
\left\{
\begin{aligned}
     \L_{-M}
  &:= -\,\nabla \cdot\bigl(e^{-M}\nabla\cdot\bigr)
      + \kappa^2
      = \kappa^2 - e^{-M}\Delta,\\
     \L_{+M}
  &:= -\,\nabla \cdot\bigl(e^{M}\nabla\cdot\bigr)
      + \kappa^2
      = \kappa^2 - e^{M}\Delta,
\end{aligned}
\right.
\end{equation*}
each equipped with the same Kirchhoff vertex conditions as $\L_u$.  
For any non‑zero $v\in H^1(\Gamma)$, introduce the Rayleigh quotients (see \cite{courant-hilbert-vol1}, pp. 401-405) 
\begin{equation*}
\left\{
\begin{aligned}
  \mathcal{R}_{-M}(v)
  &= \frac{\displaystyle \int_\Gamma e^{-M}\lvert v'\rvert^2 \, {\rm dx}
            +  \kappa^2 \int_\Gamma  |v|^2 \, {\rm dx}}
           {\|v\|_{L^2(\Gamma)}^2},\\[4pt]
  \mathcal{R}_u(v)
  &= \frac{\displaystyle \int_\Gamma e^{u}\lvert v'\rvert^2 \, {\rm dx}
            + \kappa^2 \int_\Gamma  |v|^2 \, {\rm dx}}
           {\|v\|_{L^2(\Gamma)}^2},\\[4pt]
  \mathcal{R}_{+M}(v)
  &= \frac{\displaystyle \int_\Gamma e^{M}\lvert v'\rvert^2 \, {\rm dx}
            + \kappa^2 \int_\Gamma  |v|^2 \, {\rm dx}}
           {\|v\|_{L^2(\Gamma)}^2}.
\end{aligned}
\right.
\end{equation*}
Since $e^{-M}\le e^{u}\le e^{M}$ on $\Gamma$, it follows that 
\[
  \mathcal{R}_{-M}(v)\le \mathcal{R}_u(v)\le \mathcal{R}_{+M}(v),
  \quad\forall\,v\in H^1(\Gamma)\setminus\{0\}.
\]

The paper \cite{bolin2024regularity} establishes that \(\L_u\) is self-adjoint with compact inverse on the compact, locally one-dimensional graph \(\Gamma\). Hence, the Courant–Fischer (min–max) characterization (see \cite{courant-hilbert-vol1}, pp.~405–408) applies verbatim as in the classical Euclidean or manifold setting. Applying the Courant–Fischer min–max principle, we have
\begin{equation}
\label{equa:rayleigh}
  \lambda_j(\L_{-M})
  = \min_{\substack{V\subset H^1(\Gamma)\\\dim V=j}}
    \max_{v\in V} \mathcal{R}_{-M}(v)
  \;\le\;
  \lambda_j(\L_u)
  \;\le\;
  \min_{\substack{V\subset H^1(\Gamma)\\\dim V=j}}
    \max_{v\in V} \mathcal{R}_{+M}(v)
  = \lambda_j(\L_{+M}).
\end{equation}

Classical Weyl's law \cite{MR3894619} for the constant‑coefficient operator $\L=-\Delta$
(with Kirchhoff vertex conditions) on a compact metric graph gives
constants $C'_1,C'_2>0$ (depending only on~$\Gamma$) such that
\[
  C'_1 j^2 \le \lambda_j(\L) \le C'_2 j^2.
\]
Hence
\begin{equation}
\label{equa:constant-bound}
  C'_1 e^{-M} j^2 \;\le\; \lambda_j(\L_{-M}),
  \qquad
  \lambda_j(\L_{+M}) \;\le\; C'_2 e^{M} j^2.
\end{equation}
Combining \eqref{equa:constant-bound} with \eqref{equa:rayleigh} yields
\[
  C_1 e^{-\|u\|_{L^\infty(\Gamma)}} j^2
  \;\le\;
  \lambda_j(\L_u)
  \;\le\;
  C_2 e^{\|u\|_{L^\infty(\Gamma)}} j^2,
\]
where $C_1$ and $C_2$ depend only on $\Gamma$ and $\kappa$. 

\end{proof}

 Lemma \ref{lem:weyl's law}  can be used to establish the following lemma, which controls the difference between powers of the eigenvalues of $\L_{u_1}$ and $\L_{u_2},$ for given $u_1,u_2 \in \black H^1 \nc(\Gamma),$ in terms of the $L^\infty(\Gamma)$ deviation between $u_1$ and $u_2.$

\begin{lemma}
\label{lemma:eigenvale boundedness}
    Let $\L_{u_1}$ and $\L_{u_2}$ be defined as in \eqref{equa:p} with $u_1,u_2 \in \black H^1 \nc(\Gamma)$. Then, there exists a constant $C > 0$ independent of $u_1, u_2$ such that, for $\black t\nc \ge 0,$ 
    \begin{equation*}
        \bigl|(\lambda_j^{(1)})^{-\black t\nc} - (\lambda_j^{(2)})^{-\black t\nc}\bigr| \leq C\exp\Big((\black t\nc+2)\max \bigl\{\|u_1\|_{L^\infty(\Gamma)},\|u_2\|_{L^\infty(\Gamma)} \bigr\}\Big) \,  \|u_1 - u_2\|_{L^\infty(\Gamma)},
    \end{equation*}
    where $\{\lambda_j^{(i)}\}_{j \in \mathbb{N}}$ are the eigenvalues of $\L_{u_i}$ arranged in non-decreasing order, $(i=1,2)$.
\end{lemma}

\begin{proof}
We first estimate $\lvert\lambda_j^{(1)}-\lambda_j^{(2)}\rvert$. Setting the Rayleigh quotient (see \cite{courant-hilbert-vol1}, pp.~401–404) of $\L_{u_i},\,(i=1,2)$ to be the following form
\begin{align*}
    \mathcal{R}_{u_i}(v)
  = \frac{\displaystyle \int_\Gamma e^{u_i}|v'|^2 \, {\rm dx}
          +\kappa^2\int_\Gamma |v|^2 \, {\rm dx}}
         {\|v\|_{L^2(\Gamma)}^2},
  \quad (i=1,2)
\end{align*}
then by the Courant–Fischer min–max principle (see \cite{courant-hilbert-vol1}, pp.~405–408), it holds that
\begin{align*}
      \lambda_j^{(i)}
  = \min_{\substack{V\subset H^1(\Gamma)\\\dim V=j}}
    \max_{0\neq v\in V}
    \mathcal{R}_{u_i}(v).
\end{align*}
Hence, the eigenvalues can be expressed in the min-max form
\begin{align*}
     \lambda_j^{(1)}
  = \min_{\substack{V\subset H^1(\Gamma)\\\dim V=j}}
    \max_{0\neq v\in V}\mathcal{R}_{u_1}(v),
  \quad
  \lambda_j^{(2)}
  = \min_{\substack{V\subset H^1(\Gamma)\\\dim V=j}}
    \max_{0\neq v\in V}\mathcal{R}_{u_2}(v).
\end{align*}
Let \(V_2\) be a \(j\)-dimensional subspace realizing the minimum for \(\lambda_j^{(2)}\).  Then,
\begin{align*}
     \lambda_j^{(1)} - \lambda_j^{(2)}
  \;\le\;
  \max_{0\neq v\in V_2} \mathcal{R}_{u_1}(v)
  \;-\;
  \max_{0\neq v\in V_2} \mathcal{R}_{u_2}(v)
  \;\le\;
  \max_{0\neq v\in V_2}\bigl|\mathcal{R}_{u_1}(v)-\mathcal{R}_{u_2}(v)\bigr|.
\end{align*}
Similarly, letting \(V_1\) realize the minimum for \(\lambda_j^{(1)}\) gives
\begin{align*}
     \lambda_j^{(2)} - \lambda_j^{(1)}
  \;\le\;
  \max_{0\neq v\in V_1}\bigl|\mathcal{R}_{u_1}(v)-\mathcal{R}_{u_2}(v)\bigr|.
\end{align*}

Combining these two estimates yields
\begin{align*}
      \bigl|\lambda_j^{(1)}-\lambda_j^{(2)}\bigr|
  \;\le\;
  \max\Bigl\{
    \max_{0\neq v\in V_1},\,
    \max_{0\neq v\in V_2}
  \Bigr\}
  \bigl|\mathcal{R}_{u_1}(v)-\mathcal{R}_{u_2}(v)\bigr|
  \;\le\;
  \sup_{0\neq v\in H^1(\Gamma)}
  \bigl|\mathcal{R}_{u_1}(v)-\mathcal{R}_{u_2}(v)\bigr|.
\end{align*}
Notice that, for any \(v\neq0\),
\begin{align}
  \bigl|\mathcal{R}_{u_1}(v)-\mathcal{R}_{u_2}(v)\bigr|
  &= \frac{\displaystyle\Bigl|\int_\Gamma (e^{u_1}-e^{u_2})|v'|^2 \, {\rm dx}\Bigr| }
         {\|v\|_{L^2(\Gamma)}^2} \nonumber\\[4pt]
  &\le \|e^{u_1}-e^{u_2}\|_{L^\infty(\Gamma)}\,
     \frac{\|v'\|_{L^2(\Gamma)}^2}{\|v\|_{L^2(\Gamma)}^2}
  \;\le\;
     e^M\,\|u_1-u_2\|_{L^\infty(\Gamma)}
     \,\frac{\|v'\|_{L^2(\Gamma)}^2}{\|v\|_{L^2(\Gamma)}^2},
     \label{equa: rayleigh quotient for difference}
\end{align}
where \(M:=\max\bigl\{\|u_1\|_{L^\infty(\Gamma)},\|u_2\|_{L^\infty(\Gamma)}\bigr\}\).  Taking the infimum over
all \(j\)-dimensional subspaces in the last ratio gives the \(j\)th eigenvalue
of \(-\Delta\) with Kirchhoff conditions. By the classical Weyl's law on the compact metric graph (see \cite{MR3894619}), there exists a constant $c(\Gamma)$, independent of $u$, such that
\(\lambda_j(-\Delta)\le c(\Gamma)\,j^2\), for all  $j\in\N$.  Thus, \eqref{equa: rayleigh quotient for difference} implies that 
\begin{align*}
    \bigl|\lambda_j^{(1)}-\lambda_j^{(2)}\bigr|
  \;\le\;
  e^M\,\|u_1-u_2\|_{L^\infty(\Gamma)}\;\lambda_j(-\Delta)
  \;\le\;
  c(\Gamma)\,e^M\,j^2\,\|u_1-u_2\|_{L^\infty(\Gamma)}.
\end{align*}

For general \(\black t\nc\geq0\), the mean‐value theorem applied to \(x\mapsto x^{-\black t\nc}\) yields
\begin{align*}
      \bigl|(\lambda_j^{(1)})^{-\black t\nc} - (\lambda_j^{(2)})^{-\black t\nc}\bigr|
  = \black t\nc\,\xi_j^{-\black t\nc-1}\,\bigl|\lambda_j^{(1)}-\lambda_j^{(2)}\bigr|,
\end{align*}
where \(\xi_j\) lies between \(\lambda_j^{(1)}\) and \(\lambda_j^{(2)}\).  Using
\(\xi_j^{-\black t\nc-1}\le C\,e^{(\black t\nc+1)M}j^{-2(\black t\nc+1)}\) from
Lemma~\ref{lem:weyl's law}, in which $C$ is independent of $u$ and the above \(O(j^2)\) bound, we obtain
\begin{align*}
     \bigl|(\lambda_j^{(1)})^{-\black t\nc} - (\lambda_j^{(2)})^{-\black t\nc}\bigr|
  &\;\le\;
  C\,e^{(\black t\nc+2)M}\,j^{-2\black t\nc}\,\|u_1-u_2\|_{L^\infty(\Gamma)}\\
  &\;\le\;
  C\,\exp\Bigl((\black t\nc+2)\max\bigl\{\|u_1\|_{L^\infty(\Gamma)},\|u_2\|_{L^\infty(\Gamma)}\bigr\}\Bigr)\,
  \|u_1-u_2\|_{L^\infty(\Gamma)},
\end{align*}
since \(j^{-2\black t\nc}\le1\) for \(j\ge1\) and \(\black t\nc\geq 0\).  This completes the proof.
\end{proof}

We are now ready to prove the stability of the forward map for the fractional elliptic problem.

\begin{proof}[Proof of Theorem~\ref{thm:boundedness_solution fractional operator} for $\beta>1$] 
    For $\beta > 1$, since the operator is defined in the spectral sense, we cannot directly apply the weak formulation used in the proof for the case $\beta = 1$. Instead, we begin with the series representation of the solution:
    \begin{equation}
        \label{equa: series solution}
        p = \sum_{j=1}^\infty f_j \lambda_j^{-\beta} e_j,
    \end{equation}
    where $\{\lambda_j\}$ are the non-decreasingly ordered eigenvalues of $\L_u$, $f_j = \langle f,e_j\rangle_{L^2(\Gamma)}$, and $\{e_j\}$ are the corresponding $L^2(\Gamma)$-normalized eigenfunctions.
    
    For the case $\beta = 1$, we have already established that there exists $C_0>0$ independent of $u$, such that  the solution $p_1$ satisfies:
    \begin{equation*}
        \|p_1\|^2_{L^2(\Gamma)} =\sum_{j=1}^\infty \lvert f_j \lambda_j^{-1}\rvert^2 <  C_0  \exp \bigl(2\|u\|_{L^\infty(\Gamma)}\bigr).
    \end{equation*}
    We now estimate $\|p\|_{L^2(\Gamma)}$ for the general case $\beta > 1$:
    \begin{equation}
        \label{equa: norm of p}
        \begin{aligned}
            \|p\|^2_{L^2(\Gamma)} 
            &= \sum_{j=1}^\infty \lvert f_j \lambda_j^{-\beta} \rvert^2
            = \sum_{j=1}^\infty \lvert f_j \rvert^2 \cdot \lambda_j^{-2\beta}\\
            &\overset{\hspace{-0.1cm}\text{(i)}}{\le}  C_0  \exp \bigl(2\beta\|u\|_{L^\infty(\Gamma)} \bigr)\cdot  \sum_{j=1}^\infty\lvert f_j \rvert^2j^{-4\beta} 
            \leq  c_0(\beta) \exp \bigl(2\beta\|u\|_{L^\infty(\Gamma)} \bigr),
        \end{aligned}
    \end{equation}
    where inequality (i) follows from Lemma~\ref{lem:weyl's law}, and $c_0(\beta)$ is a decreasing function of $\beta$.

    Next, we analyze $\|\nabla p\|_{L^2(\Gamma)}$. We begin with the eigenfunctions $\{e_j\}$. By the definition of the eigenfunction $e_j$, we have
    \begin{equation*}
        \bigl(\kappa^2 - \nabla\cdot(\exp(u)\nabla \cdot)\bigr) \, e_j = \lambda_j e_j.
    \end{equation*}
    Taking the inner product with $e_j$ in $L^2(\Gamma)$ and applying integration by parts and using that the eigenfunctions are in the domain of the operator, and thus satisfy the Kirchhoff vertex conditions, we obtain 
    \begin{equation*}
        \int_{\Gamma} \kappa^2 e_j^2 \, {\rm dx} +  \int_\Gamma \exp(u)\lvert \nabla e_j \rvert^2 \, {\rm dx} = \lambda_j \int_{\Gamma} e_j^2\,{\rm dx}.
    \end{equation*}
    Using the normalization condition $\|e_j\|_{L^2(\Gamma)} = 1$, it follows that
    \begin{equation*}
        \int_{\Gamma} \exp(u)\lvert \nabla e_j \rvert^2 \,{\rm dx} = \lambda_j - \kappa^2.
    \end{equation*}
    Since $u \in L^\infty(\Gamma)$,  we deduce that
    \begin{equation*}
        \exp \bigl(-\|u\|_{L^\infty(\Gamma)} \bigr) \|\nabla e_j\|_{L^2(\Gamma)}^2 \leq \lambda_j - \kappa^2 \leq \exp \bigl(\|u\|_{L^\infty(\Gamma)} \bigr) \|\nabla e_j\|_{L^2(\Gamma)}^2.
    \end{equation*}
    Consequently,
    \begin{equation*}   
        \exp \bigl(-\|u\|_{L^\infty(\Gamma)} \bigr) (\lambda_j - \kappa^2)
        \leq 
        \|\nabla e_j\|^2_{L^2(\Gamma)}\leq \exp\bigl(\|u\|_{L^\infty(\Gamma)}\bigr) (\lambda_j - \kappa^2),
    \end{equation*}
    implying that there exist constants $C_1$, $C_2>0$, independent of $u$, such that
    \begin{equation}
        \label{equa: order estimation of nabla ej}
        C_1\exp\Bigl(-\frac{1}{2}\|u\|_{L^\infty(\Gamma)} \Bigr) \lambda_j^{1/2}
        \leq 
        \|\nabla e_j\|_{L^2(\Gamma)}
        \leq 
        C_2\exp \Bigl(\frac{1}{2}\|u\|_{L^\infty(\Gamma)}\Bigr) \lambda_j^{1/2}.
    \end{equation}
    Noting that $\Gamma$ is a compact metric graph, term-by-term differentiation of the series is justified due to uniform convergence. We estimate $\|\nabla p\|_{L^2(\Gamma)}$ as follows:
 \begin{align}
  \|\nabla p\|_{L^2(\Gamma)}
  &= \Bigl\|\sum_{j=1}^\infty f_j \lambda_j^{-\beta}\,\nabla e_j \Bigr\|_{L^2(\Gamma)} \nonumber
  \le \sum_{j=1}^\infty \lvert f_j\lambda_j^{-\beta}\rvert\,\|\nabla e_j\|_{L^2(\Gamma)} \nonumber\\
&\overset{\hspace{-0.1cm}\text{(ii)}}{\le} C_2\exp \Bigl(\frac{1}{2}\|u\|_{L^\infty(\Gamma)} \Bigr)\sum_{j=1}^\infty |f_j|\,\lambda_j^{-\beta}\,\lambda_j^{1/2}  \nonumber
  = C_2\exp \Bigl(\frac{1}{2}\|u\|_{L^\infty(\Gamma)} \Bigr)\sum_{j=1}^\infty |f_j|\,\lambda_j^{-(\beta-\frac{1}{2})} \nonumber\\
  &\overset{\hspace{-0.1cm}\text{(iii)}}{\le} C_3\,\exp\biggl({ \Bigl(\beta-\frac{1}{2}+\frac{1}{2} \Bigr)\|u\|_{L^\infty(\Gamma)}}\biggr)
     \sum_{j=1}^\infty |f_j|^2 j^{-(\beta-\frac{1}{2})}
      \nonumber\\
  &\le  c_1(\beta) \,\exp\bigl(\beta\|u\|_{L^\infty(\Gamma)}\bigr).
    \label{equa: nabla p}
  \end{align}
  Here, we obtained inequality~(ii) by using \eqref{equa: order estimation of nabla ej}, and inequality~(iii) follows directly since by assumption $f\in L^2(\Gamma)$. Combining estimates \eqref{equa: norm of p} and \eqref{equa: nabla p}, we obtain a bound for the $H^1(\Gamma)$ norm of $p$ when $\beta > 1$:
    \begin{equation}
        \label{equa: sobolev norm of p}
        \|p\|_{H^1(\Gamma)} \leq c(\beta) \exp\bigl(\beta\|u\|_{L^\infty(\Gamma)}\bigr),
    \end{equation}
    where $c(\beta)$ is independent of $u$ and a decreasing function of $\beta$.

We now derive an estimate for $\|p_1-p_2\|_{H^1(\Gamma)}$. Employing the series representation of the solutions $p_1$ and $p_2$, we express their difference as
\begin{align*}
\tilde{p} &:= p_1-p_2
 = \sum_{j=1}^\infty f_j^{(1)}(\lambda_j^{(1)})^{-\beta} e_j^{(1)} - \sum_{j=1}^\infty f_j^{(2)}(\lambda_j^{(2)})^{-\beta} e_j^{(2)}\\[1mm]
& = \sum_{j=1}^\infty \Bigl[ f_j^{(1)}(\lambda_j^{(1)})^{-\beta} e_j^{(1)} - f_j^{(2)}(\lambda_j^{(2)})^{-\beta} e_j^{(2)} \Bigr]\\[1mm]
& = \sum_{j=1}^\infty \Bigl[ (\lambda_j^{(1)})^{-(\beta-1)}\, f_j^{(1)}(\lambda_j^{(1)})^{-1} e_j^{(1)} - (\lambda_j^{(2)})^{-(\beta-1)}\, f_j^{(2)}(\lambda_j^{(2)})^{-1} e_j^{(2)} \Bigr],
\end{align*}
where $\{(\lambda_j^{(1)},\,e_j^{(1)})\}_{j\in \N}$ and $\{(\lambda_j^{(2)},\,e_j^{(2)})\}_{j\in \N}$ represent the eigenpairs of the operators $\L_{u_1}$ and $\L_{u_2}$ with standard Kirchhoff vertex conditions, ordered and normalized in the standard way.

We first estimate $\|p_1-p_2\|_{L^2(\Gamma)}$. Since the eigenfunctions for different operators $\L_{u_1}$ and $\L_{u_2}$ are no longer orthogonal, we will need to employ the results for $\beta=1$, as well as Lemmas \ref{lem:weyl's law} and \ref{lemma:eigenvale boundedness}. 
As shown in \cite{bolin2024regularity}, $\L_{u_1}$ is self-adjoint, positive definite with compact inverse in $L^2(\Gamma)$. Hence, the spectral theorem yields that for $\beta>1$, $\L_{u_1}^{-(\beta-1)}$ is again self-adjoint, positive definite, and compact, with discrete spectrum.
Thus, we can estimate the operator norm of $\L_{u_1}^{-(\beta-1)}$ by using Lemma~\ref{lem:weyl's law},
\begin{align}
     \bigl\|\L_{u_1}^{-(\beta-1)}\bigr\|_{op}
  = \sup_j (\lambda_j^{(1)})^{-(\beta-1)}
  = (\lambda_1^{(1)})^{-(\beta-1)}
  \leq C\exp\big(-(\beta-1)\|u_1\|_{L^\infty(\Gamma)}\big),
  \label{equa:operator norm}
\end{align}
where $C$ is a constant independent of $u_1$. 
Therefore,
\begin{align}
  \|p_1-p_2\|_{L^2(\Gamma)}
  &=\Bigl\|\sum_{j=1}^{\infty}(\lambda_j^{(1)})^{-(\beta-1)}
          f_j^{(1)}(\lambda_j^{(1)})^{-1}e_j^{(1)}
        -(\lambda_j^{(2)})^{-(\beta-1)}
          f_j^{(2)}(\lambda_j^{(2)})^{-1}e_j^{(2)}
     \Bigr\|_{L^2(\Gamma)}\nonumber\\
  &\le \Bigl\|\sum_{j=1}^{\infty}(\lambda_j^{(1)})^{-(\beta-1)}
          \bigl[f_j^{(1)}(\lambda_j^{(1)})^{-1}e_j^{(1)}
               -f_j^{(2)}(\lambda_j^{(2)})^{-1}e_j^{(2)}\bigr]
     \Bigr\|_{L^2(\Gamma)}\nonumber\\
  &\hspace{1.5cm}
      +\Bigl\|\sum_{j=1}^{\infty}\bigl[(\lambda_j^{(1)})^{-(\beta-1)}
               -(\lambda_j^{(2)})^{-(\beta-1)}\bigr]\,
            f_j^{(2)}(\lambda_j^{(2)})^{-1}e_j^{(2)}
     \Bigr\|_{L^2(\Gamma)}\nonumber\\
  &\overset{\text{(iv)}}{\le}
     \|\L_{u_1}^{-(\beta-1)}\|_{\mathrm{op}}\,
     \|p_1^{(1)}-p_2^{(1)}\|_{L^2(\Gamma)}\nonumber\\
  &\hspace{1.5cm}
      +\Bigl\|\sum_{j=1}^{\infty}\bigl[(\lambda_j^{(1)})^{-(\beta-1)}
               -(\lambda_j^{(2)})^{-(\beta-1)}\bigr]\,
            f_j^{(2)}(\lambda_j^{(2)})^{-1}e_j^{(2)}
     \Bigr\|_{L^2(\Gamma)}\nonumber\\
  &\overset{\text{(v)}}{\le}
      c_2(\beta) \,\exp \bigl((\beta-1)\|u_1\|_{L^\infty(\Gamma)}\bigr)\,
     \|p_1^{(1)}-p_2^{(1)}\|_{L^2(\Gamma)}\nonumber\\
  &\hspace{1.5cm}
     + c_3(\beta) \,\exp \bigl((\beta+1)
            \max \bigl\{\|u_1\|_{L^\infty(\Gamma)},\|u_2\|_{L^\infty(\Gamma)}\bigr\}\bigr)\,
       \|p_2^{(1)}\|_{L^2(\Gamma)}\nonumber\\
  &\overset{\text{(vi)}}{\le}
      c_4(\beta)\,\exp \Bigl((\beta+2)
            \max \bigl\{\|u_1\|_{L^\infty(\Gamma)},\|u_2\|_{L^\infty(\Gamma)} \bigr\}\Bigr)\,
       \|u_1-u_2\|_{L^\infty(\Gamma)},
  \label{equa: norm of p1-p2}
\end{align}
where $p_1^{(1)}$ and $p_2^{(1)}$ denote the solutions of  \eqref{equa:p} with $\beta=1$ and parameters $u_1$ and $u_2,$ respectively. In {inequality~(iv)} we invoke the operator–norm estimate \eqref{equa:operator norm}. Inequality~(v) relies on Lemmas~\ref{lem:weyl's law} and \ref{lemma:eigenvale boundedness}, while inequality~(vi) is obtained by applying Theorem~\ref{thm:boundedness_solution fractional operator} with \(\beta = 1\).

Next, we estimate $\|\nabla p_1-\nabla p_2\|_{L^2(\Gamma)}$. Proceeding analogously, we have
\begin{align}
  \|\nabla p_1 - \nabla p_2\|_{L^2(\Gamma)}
  &= \Bigl\|\sum_{j=1}^\infty (\lambda_j^{(1)})^{-(\beta-1)}\,f_j^{(1)}(\lambda_j^{(1)})^{-1}\nabla e_j^{(1)}
       - (\lambda_j^{(2)})^{-(\beta-1)}\,f_j^{(2)}(\lambda_j^{(2)})^{-1}\nabla e_j^{(2)}\Bigr\|_{L^2(\Gamma)}
  \nonumber\\
  &\le \Bigl\|\sum_{j=1}^\infty (\lambda_j^{(1)})^{-(\beta-1)}
        \bigl[f_j^{(1)}(\lambda_j^{(1)})^{-1}\nabla e_j^{(1)}
             - f_j^{(2)}(\lambda_j^{(2)})^{-1}\nabla e_j^{(2)}\bigr]\Bigr\|_{L^2(\Gamma)}
  \nonumber\\
  &\hspace{1.5cm}
    + \Bigl\|\sum_{j=1}^\infty \bigl[(\lambda_j^{(1)})^{-(\beta-1)}
        - (\lambda_j^{(2)})^{-(\beta-1)}\bigr]\,
        f_j^{(2)}(\lambda_j^{(2)})^{-1}\nabla e_j^{(2)}\Bigr\|_{L^2(\Gamma)}
  \nonumber\\
  &\overset{\hspace{-0.1cm}\text{(vii)}}{\le} \max_j(\lambda_j^{(1)})^{-(\beta-1)}
         \Bigl\|\sum_{j=1}^\infty\bigl[f_j^{(1)}(\lambda_j^{(1)})^{-1}\nabla e_j^{(1)}
                                - f_j^{(2)}(\lambda_j^{(2)})^{-1}\nabla e_j^{(2)}\bigr]\Bigr\|_{L^2(\Gamma)}
  \nonumber\\
  &\hspace{1.5cm}
    + \Bigl\|\sum_{j=1}^\infty \bigl[(\lambda_j^{(1)})^{-(\beta-1)}
        - (\lambda_j^{(2)})^{-(\beta-1)}\bigr]\,
        f_j^{(2)}(\lambda_j^{(2)})^{-1}\nabla e_j^{(2)}\Bigr\|_{L^2(\Gamma)}
  \nonumber\\
  &\le  c_2(\beta)\,\exp \bigl((\beta-1)\|u_1\|_{L^\infty(\Gamma)}\bigr)\,
         \|\nabla p_1^{(1)} - \nabla p_2^{(1)}\|_{L^2(\Gamma)}
  \nonumber\\
  &\hspace{1.5cm}
    +  c_3(\beta) \,\exp \bigl((\beta+1)\|u_1\|_{L^\infty(\Gamma)}\bigr)\,
      \|u_1 - u_2\|_{L^\infty(\Gamma)}\,\|\nabla p_2\|_{L^2(\Gamma)}
  \nonumber\\
  & \overset{\hspace{-0.1cm}\text{(viii)}}{\le}   c_4(\beta) \,\exp \bigl((\beta+2)\|u_1\|_{L^\infty(\Gamma)}\bigr)\,
         \|u_1 - u_2\|_{L^\infty(\Gamma)},
  \label{equa: norm of nabla p1-p2}
\end{align}
where in (vii) we used Lemmas \ref{lem:weyl's law} and \ref{lemma:eigenvale boundedness}, and (viii) follows from Theorem~\ref{thm:boundedness_solution fractional operator} with $\beta = 1,$ proved in Subsection~\ref{ssec:elliptic}.
Combining the estimates in \eqref{equa: norm of p1-p2} and \eqref{equa: norm of nabla p1-p2}, we deduce that
\begin{equation*}
    \|p_1-p_2\|_{H^1(\Gamma)} \leq c(\beta) \, \exp \bigl((\beta+2)\|u_1\|_{L^\infty(\Gamma)} \bigr)\|u_1-u_2\|_{L^\infty(\Gamma)},
\end{equation*}
as desired.
\end{proof}

\section{Numerical Experiments}\label{sec:numerics}

In this section, we illustrate the numerical solution of elliptic and fractional elliptic Bayesian inverse problems on compact metric graphs. Subsection~\ref{ssec:implementation} details the algorithmic implementation, while Subsections~\ref{ssec:output standard case} and~\ref{ssec:output fractional case} present numerical results for the elliptic and fractional elliptic problems, respectively.
\black We compare the posterior mean estimator with a \emph{maximum a posteriori} (MAP) estimator computed using Markov chain Monte Carlo (MCMC) output. Specifically, in the high-dimensional discretized setting considered here, the MAP-type estimator we report is a \emph{pointwise marginal MAP} (posterior marginal mode), obtained by maximizing the one-dimensional marginal posteriors at each mesh node; see Equation~\ref{eq:u_mmap} below for the precise definition and computation.
In addition, we explore the quantification of uncertainty through the posterior marginal standard deviation along points in the graph.\nc

\subsection{Implementation}
\label{ssec:implementation}
We work on a letter-shaped metric graph naturally embedded in $\R^2$ (see e.g. Figure \ref{fig:mu-p-comparison}) and adopt the Bayesian  framework introduced in Section~\ref{sec:problemformulation}. For the forward map induced by Equation~\eqref{equa:p}, we set $\kappa=1$ and $f(x)=z_1^2(x)-z_2^2(x)$, where \( \bigl(z_1(x),z_2(x)\bigr)\) denotes the Cartesian coordinates of the point \(x\in \Gamma\subset \mathbb R^{2}\). We examine both the elliptic case with 
\(\beta = 1\) (Subsection~\ref{ssec:output standard case}) and the fractional elliptic case with \(\beta = 3/2\) (Subsection~\ref{ssec:output fractional case}). The prior for $u$ is a Whittle–Matérn Gaussian field as specified in Equation~\eqref{equa:prior gaussian}, with smoothness parameter $\alpha = 1$. 
\black The validity of this choice is ensured by Proposition~\ref{prop:forward_under_H1}, 
and it facilitates the use of existing sampling and discretization strategies. \nc 
We further take \(\black \kappa_{\rm ref}\nc = \sqrt{0.2} \cdot {2}/{3}\) and \({\color{black} a_{\rm ref}} = 0.2\) \black in Equation~\eqref{equa:prior gaussian}\nc, which leads to a correlation range of $3$ in the setting of the \textsf{MetricGraph}  R package \cite{MetricGraph}.

We discretize the letter-shaped metric graph in a fine mesh and approximate the forward model and the prior using the finite-element method introduced in \cite{bolin2019rational} and implemented in the \textsf{rSPDE}  {\cite{rSPDE}} and \textsf{MetricGraph}  {\cite{MetricGraph}} R packages. \black A brief recap of the finite-element discretization (including the assembly of the discrete operators) and of the quadrature-based rational approximation used for the fractional case is provided in Appendix~\ref{app:fem-rational}. \nc
As we choose $\alpha = 1$, the random field $u$ can be simulated exactly using the methods in \cite{bolin2023statistical}; however, as we discretize the forward problem, we also use a finite-element approximation of $u$ as introduced in \cite{bolin2024regularity} and further refined in \cite{bolin2025new}. 
To balance accuracy and computational efficiency, we use a mesh size of \(h = 0.05\) for the elliptic problem and \(h = 0.2\) for the fractional elliptic case, thus reducing the computational cost in the fractional elliptic case. For the latter case, we employ the operator-based approach from \cite{bolin2019rational}.

\black
Let $\{x_i\}_{i=1}^{N_h}$ denote the mesh nodes on $\Gamma$. In the finite-element discretization, we represent the log-conductivity and the state by their nodal degrees of freedom, storing the vectors
\(
u=\bigl(u(x_1),\dots,u(x_{N_h})\bigr)^\top\in\mathbb R^{N_h},\,
p(u)=\bigl(p(x_1),\dots,p(x_{N_h})\bigr)^\top\in\mathbb R^{N_h},
\)
or, equivalently, the corresponding coefficient vectors in the nodal finite-element basis.
The numbers of mesh nodes used for the elliptic and fractional elliptic problems are $N_h=1{,}998$ and $N_h=507$, respectively, chosen to balance accuracy and computational cost. \nc In each case, we draw a (discretized) ground-truth parameter \(u_0\sim\mu_0\) from the prior and compute \(p_0:=p(u_0)\). 
We then generate point-wise observations at all mesh points (i.e.\ a full observation setting on the chosen discretization) by adding spatially-varying Gaussian noise:
\begin{equation}
    \label{eq:noise model}
    y_i \;=\; p_0[ i  ] \;+\;\bigl(\mathsf{n}_{\rm rel}\,\lvert p_0[i ]\rvert \;+\;\mathsf{n}_{\rm abs}\bigr)\,\varepsilon_i, 
  \qquad \varepsilon_i\iid \mathcal{N}(0,1),\quad i=1,\dots,{N_h},
\end{equation}
where \black \(\mathsf{n}_{\rm rel}=0.03\) \nc controls the relative noise level and \black \(\mathsf{n}_{\rm abs}=0.40\) \nc ensures a nonnegligible noise floor even when \(p_0[ \, i \, ]\) is small. 
This mixed noise model yields meaningful uncertainty at both large and small scales, and helps avoid the ill-conditioning of the likelihood that would emerge when using a single fixed noise level over the wide range of values of \(p_0\). 
\black
Under our full pointwise observation design on the chosen discretization, the data dimension also equals $N_h$. At the same time, $N_h$ is the dimension of the discretized unknown $u$ and of the discretized state $p(u)$. Therefore, changing $N_h$ between the elliptic and fractional elliptic experiments changes \emph{both} the spatial resolution and the dimensionality of the discretized inverse problem. Consequently, reconstruction errors should not be interpreted as depending monotonically on $N_h$ alone.
\nc

For posterior sampling, we use the preconditioned Crank–Nicolson (pCN) MCMC method \cite{cotter2013}, which preserves the Gaussian prior in the proposal and reduces the Metropolis–Hastings ratio to the posterior-to-prior term alone, thereby avoiding explicit covariance inversions. The discretization-independent rate of convergence (uniform spectral gap) of pCN was established in \cite{hairer2014spectral}; see also \cite{trillos2017consistency} for a study of pCN on (non-metric) graphs.    
Here, we further incorporate a simple self‐adaptive update for the step size \(\tau\) to balance acceptance and exploration, and apply temperature annealing to accelerate convergence of the chain; we refer to \cite[Chapter 7]{sanz2024first} for an introduction to annealing strategies for Monte Carlo methods. The procedure is outlined in Algorithm~\ref{alg:adaptive-pcn} below.

\begin{minipage}{0.95\linewidth} 
\centering
\begin{algorithm}[H]
\footnotesize
\caption{\small Adaptive pCN MCMC with Temperature Annealing}
\label{alg:adaptive-pcn}
\begin{algorithmic}[1]
  \Require Prior \(\mu_0\), data \(y\) and 
  potential function $\Phi(\cdot; y),$
  initial step size \(\tau\), minimum step size \(\tau_{\min}\), annealing start \(T_0\), cooling factor \(\zeta\), sample size \(N\), adaptation interval \(N_{\mathrm{adapt}}\), target acceptance rate \(r_{\mathrm{target}}\), burn-in \(B\).
  \State {\bf Initialization:} Sample \(u^{(0)}\sim\mu_0\). 
  \For{\(n=1,\dots,N\)}
    \State {\bf Proposal step:} Set 
      \[
        v := \sqrt{1-\tau^2}\;u^{(n-1)} \;+\;\tau\,\xi, \qquad \xi\sim\mu_0.
      \]
    \State {\bf Accept/reject step:} Set temperature $T_n :=\max\bigl(1,\;T_0\,\zeta^{\lfloor n/N_{\mathrm{adapt}}\rfloor}\bigr),$ and set 
    \[
    u^{(n)} := 
    \begin{cases}
      v, & \text{with probability  
      \(\min \Bigl\{1, \exp \Bigl( \frac{\Phi(u^{(n-1)}\,;y) - \Phi(v; y )}{T_n}\Bigr) \Bigr\}\)}, \\
      u^{(n-1)}, & \text{with probability  \( 1 - \min \Bigl\{1,\exp \bigl(\frac{ \Phi(u^{(n-1)}\,;y) - \Phi(v; y )}{T_n} \bigr) \Bigr\}\)}.
    \end{cases}
    \]
    \State Record acceptance indicator.
    \If{\(n \bmod N_{\mathrm{adapt}} = 0\)}
      \State Compute recent acceptance rate \(\bar r\).
      \If{\(\bar r < 0.9\,r_{\mathrm{target}}\)} 
         \(\tau \gets \max(0.9\,\tau,\tau_{\min})\).
      \ElsIf{\(\bar r > 1.1\,r_{\mathrm{target}}\)}
         \(\tau \gets 1.2\,\tau\).
      \EndIf
    \EndIf
  \EndFor
  \State \Return \(\{u^{(n)}\}_{n=B+1}^N\) (discard first \(B\) burn-in samples) and acceptance history.
\end{algorithmic}
\end{algorithm}
\end{minipage}

\bigskip

\black
We next define the two point estimators used throughout the numerical results, together with the corresponding plug-in reconstructions and nodewise uncertainty quantifiers computed from the retained MCMC samples.
Let $\{u^{(n)}\}_{n=B+1}^{N}$ denote the post burn-in MCMC samples. The posterior mean estimator is computed nodewise as
\begin{equation}
    \label{eq:u_postmean}
    \hat u_{\rm mean}[i] \;=\; \frac{1}{N-B}\sum_{n=B+1}^{N} u^{(n)}[i], \qquad i=1,\dots,N_h.
\end{equation}

\black The MAP-type point estimator we report is the \emph{pointwise marginal MAP} (posterior marginal mode): for each node $i$, we fit a one-dimensional kernel density estimator to $\{u^{(n)}[i]\}_{n=B+1}^{N}$ and define
\begin{equation}
    \label{eq:u_mmap}
    \hat u_{\rm mode}[i] \;=\; \arg\max_{z\in\mathbb R}\ \widehat{\pi}_i(z), \qquad i=1,\dots,N_h,
\end{equation}
where $\widehat{\pi}_i$ is the estimated marginal posterior density of $u[i]\mid y$. This marginal-MAP summary is computed componentwise and, in general, differs from the joint MAP in $\mathbb R^{N_h}$. Indeed, because the forward map $u\mapsto \mathcal F(u)$ is nonlinear, the negative log-likelihood is typically nonconvex in $u$; consequently, the posterior need not be log-concave and may exhibit multiple local modes. 
We therefore use \eqref{eq:u_mmap} as a computationally convenient MAP-type summary extracted from the MCMC output.
We then form plug-in reconstructions of the state by solving the forward problem at these point estimates,
\[
\hat p_{\rm mean} := p(\hat u_{\rm mean}),
\qquad 
\hat p_{\rm mode} := p(\hat u_{\rm mode}).
\]
To quantify uncertainty, we compute the posterior marginal standard deviation of $u$ at each node from the retained samples,
\begin{equation}
    \label{eq:post_sd}
     \widehat{\mathrm{sd}}_{u\mid y}[i]
  \;=\;
  \Bigg(\frac{1}{N-B-1}\sum_{n=B+1}^{N}\Big(u^{(n)}[i]-\hat u_{\rm mean}[i]\Big)^2\Bigg)^{1/2},
  \qquad i=1,\dots,N_h.
\end{equation}

\nc

\black
To provide quantitative comparisons between the two point estimators $(\hat u_{\rm mean},\hat p_{\rm mean})$ and $(\hat u_{\rm mode},\hat p_{\rm mode})$, we report both absolute and scale-free relative errors for the log-conductivity $u$, the state $p$, and the physical conductivity $a=\exp(u)$.
At the continuum level, for any field $v$ on $\Gamma$, we consider the $L^{2}(\Gamma)$ error $\|v_{\rm est}-v_0\|_{L^{2}(\Gamma)}$ and its relative version $\|v_{\rm est}-v_0\|_{L^{2}(\Gamma)}/\|v_0\|_{L^{2}(\Gamma)}$.\nc
\black
In our finite-element implementation, $v_0$ and $v_{\rm est}$ are represented by nodal vectors in $\mathbb R^{N_h}$.
We compute the following discrete metrics:
\[
  \mathrm{RMSE}(v_{\rm est},v_0)
  \;:=\;
  \Bigg(\frac{1}{N_h}\sum_{i=1}^{N_h}\bigl(v_{\rm est}[i]-v_0[i]\bigr)^2\Bigg)^{1/2},
  \qquad
  \mathrm{RelErr}_{L^2(\Gamma)}(v_{\rm est},v_0)
  \;:=\;
  \frac{\|v_{\rm est}-v_0\|_{L_h^2(\Gamma)}}{\|v_0\|_{L^2_h(\Gamma)}},
\]
where $\|w\|_{L^2_h(\Gamma)} := (w^\top C w)^{1/2}$ denotes the discrete $L^2(\Gamma)$ norm induced by the finite-element mass matrix $C$.
The mass matrix $C$ is assembled in the standard way in the finite-element discretization; see, e.g., \cite{bolin2024regularity}.
We apply these definitions with $(v_{\rm est},v_0)=(\hat u_{\rm mean},u_0)$ and $(\hat u_{\rm mode},u_0)$ for the parameter, with $(v_{\rm est},v_0)=(\hat p_{\rm mean},p_0)$ and $(\hat p_{\rm mode},p_0)$ for the state, and with $(v_{\rm est},v_0)=(\exp(\hat u_{\rm mean}),\exp(u_0))$ and $(\exp(\hat u_{\rm mode}),\exp(u_0))$ for the conductivity. \nc
\black
For the log-conductivity $u$, we additionally report a range-normalized error
\[
  \mathrm{NRMSE}_{\rm range}(u_{\rm est},u_0)
  \;:=\;
  \frac{\mathrm{RMSE}(u_{\rm est},u_0)}{\max_{1\le i\le N_h} u_0[i]-\min_{1\le i\le N_h} u_0[i]},
\]
which remains informative across repetitions in which $\|u_0\|_{L^2(\Gamma)}$ (and hence $\|u_0\|_2$) may be small.\nc

\black \black
In our parameterization, $u$ is the \emph{log}-conductivity and the physical coefficient entering the PDE is $a=\exp(u)$. For some randomly drawn ground truths $u_0\sim\mu_0$, the norm $\|u_0\|_{L^2(\Gamma)}$ can be small, so the standard relative error $\|u_{\rm est}-u_0\|_{L^2(\Gamma)}/\|u_0\|_{L^2(\Gamma)}$ may be artificially inflated and thus not fully representative of reconstruction quality.
We therefore complement $\mathrm{RMSE}(u_{\rm est},u_0)$ with the scale-free range-normalized error $\mathrm{NRMSE}_{\rm range}(u_{\rm est},u_0)$. This choice is also consistent with our prior specification, which is calibrated to avoid unrealistically large amplitudes or excessively rapid spatial variations in $u$ that would lead to highly heterogeneous conductivities $a=\exp(u)$ and, in turn, to sharply peaked posteriors and poor MCMC mixing.
Moreover, since the forward model is driven by the physical conductivity, we additionally report the conductivity error $\mathrm{RelErr}_{L^2(\Gamma)}(a_{\rm est},a_0)$ with $a_{\rm est}=\exp(u_{\rm est})$ and $a_0=\exp(u_0)$. This yields a physically meaningful notion of “relative error” while remaining comparable across repetitions.\nc

In the experiments described in the following subsections, we carefully select hyperparameters such as $B$, $\tau$, and $\tau_{\rm min}$ to balance exploration and acceptance rate.
\black
In the elliptic experiment we use $N = 60{,}000$ iterations with burn-in $B = 30{,}000$, while in the fractional elliptic experiment we use $N = 70{,}000$ iterations with the same burn-in $B = 30{,}000$ (thus retaining $30{,}000$ and $40{,}000$ samples, respectively). We use the same adaptive strategy in both cases.
To assess robustness across different realizations, we repeat the experiment for multiple independently drawn ground truths (i.e.\ different prior sample paths $u_0$ and the corresponding $p_0$). In the elliptic case we report results over $30$ independent trials, whereas in the fractional elliptic case we report results over $10$ independent trials.
The slightly longer chain in the fractional elliptic case is motivated by the fact that the fractional forward model (and the induced likelihood) typically yields a more challenging posterior, with stronger correlations and slower mixing in practice. Empirically, we therefore retain more post burn-in samples in the fractional elliptic experiment to stabilize Monte Carlo estimates of the posterior mean, the pointwise marginal mode, and the pointwise posterior standard deviation curves.
\nc
In both cases we initialize $\tau=0.3$, impose $\tau_{\min}=0.01$, set $N_{\rm adapt}=500$, and target acceptance rate $r_{\rm target}=0.40$. For annealing we take $T_0=20$ and $\zeta=0.95$, so that $T_n=\max(1,T_0\zeta^{\lfloor n/N_{\rm adapt}\rfloor})$.
The typical acceptance rate during the stable period ranges from $38\%$ to $42\%$, indicating effective exploration while maintaining computational efficiency.

\subsection{Elliptic Problem on a Letter-shaped Metric Graph}
\label{ssec:output standard case} 

\begin{figure}[htbp]
  \centering
  \textbf{Elliptic problem: estimates for parameter and PDE solution}\\[3ex]
  \begin{subfigure}[b]{0.5\textwidth}
    \includegraphics[width=\textwidth]{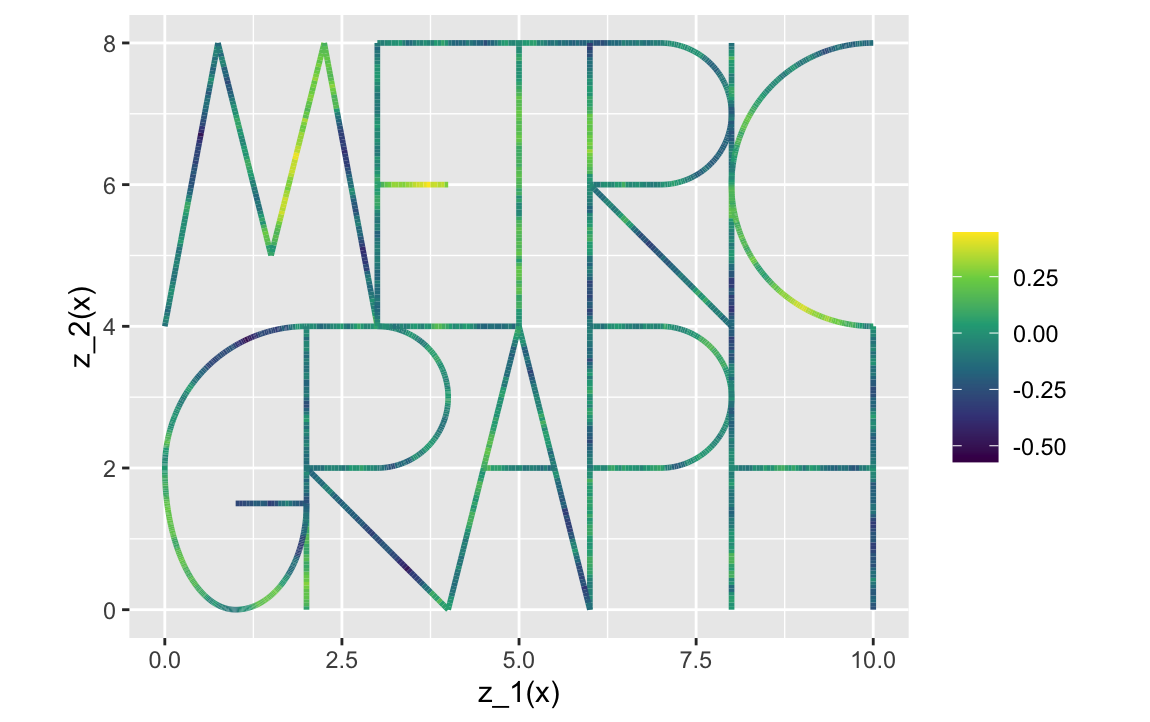}
    \caption{Ground truth $u_0$}
    \label{fig:mu-true}
  \end{subfigure}\hfill
  \begin{subfigure}[b]{0.5\textwidth}
    \includegraphics[width=\textwidth]{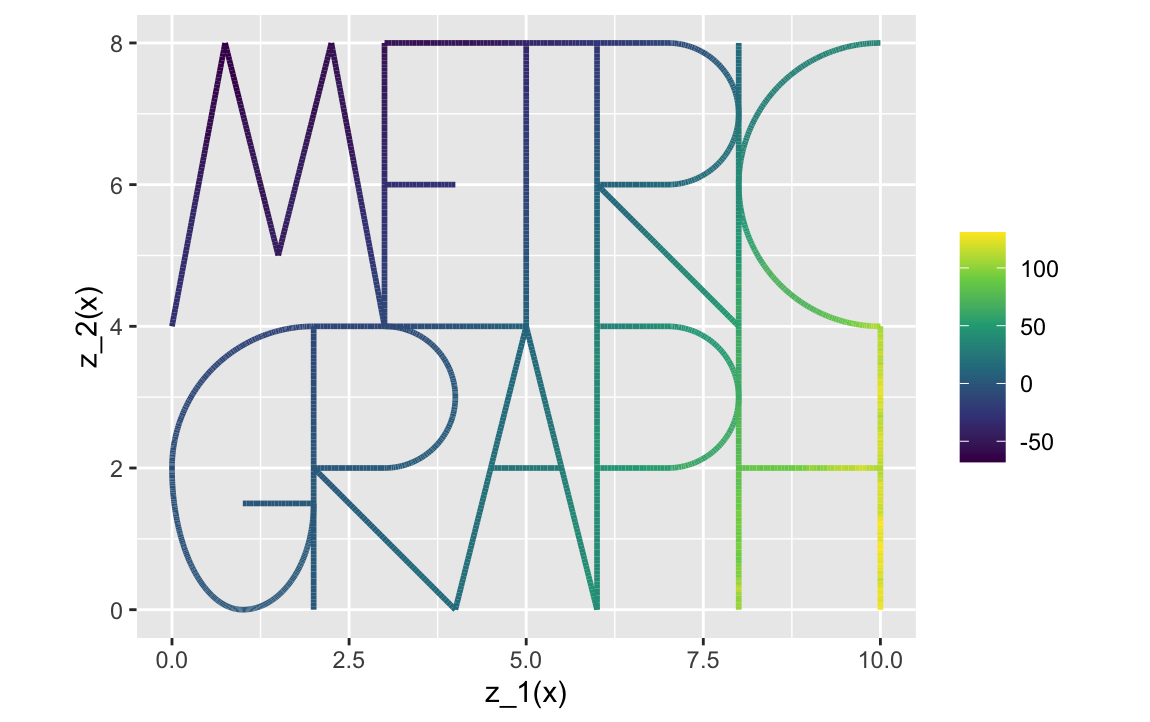}
    \caption{True solution $p_0$}
    \label{fig:p-true}
  \end{subfigure}
  
  
  \begin{subfigure}[b]{0.5\textwidth}
    \includegraphics[width=\textwidth]{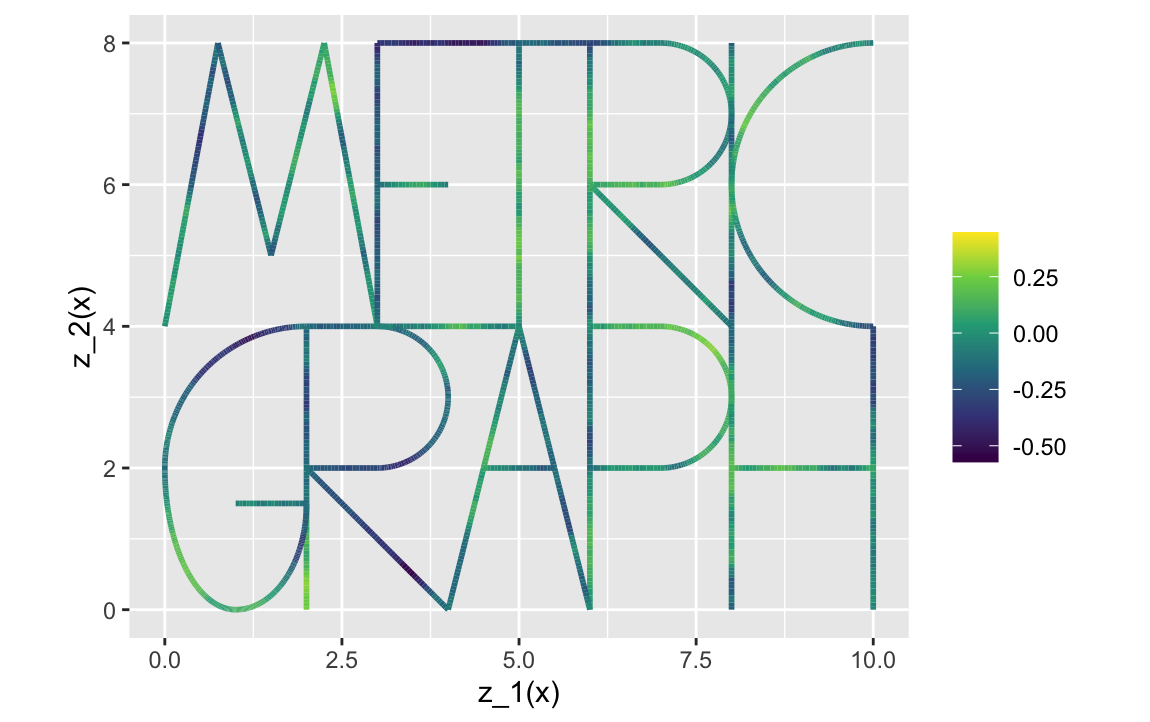}
    \caption{Posterior mean estimate of $u_0$}
    \label{fig:mu-postmean}
  \end{subfigure}\hfill
  \begin{subfigure}[b]{0.5\textwidth}
    \includegraphics[width=\textwidth]{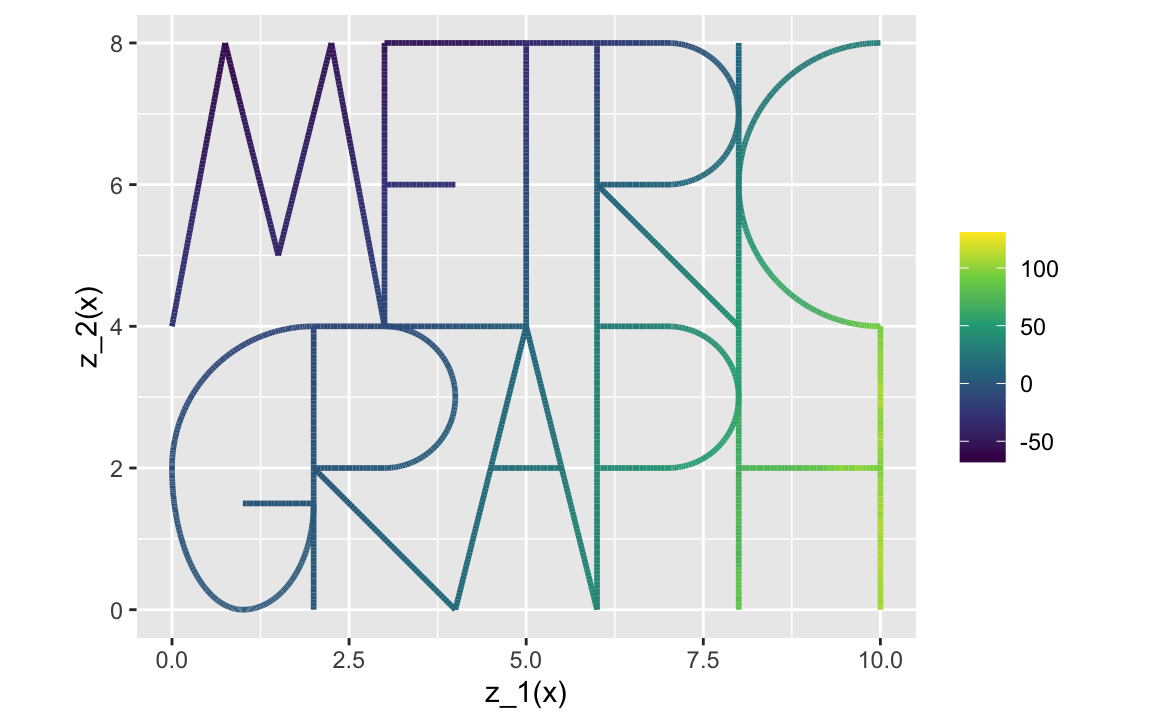}
    \caption{Posterior estimate of $p_0$}
    \label{fig:p-postmean}
  \end{subfigure}


  \begin{subfigure}[b]{0.5\textwidth}
    \includegraphics[width=\textwidth]{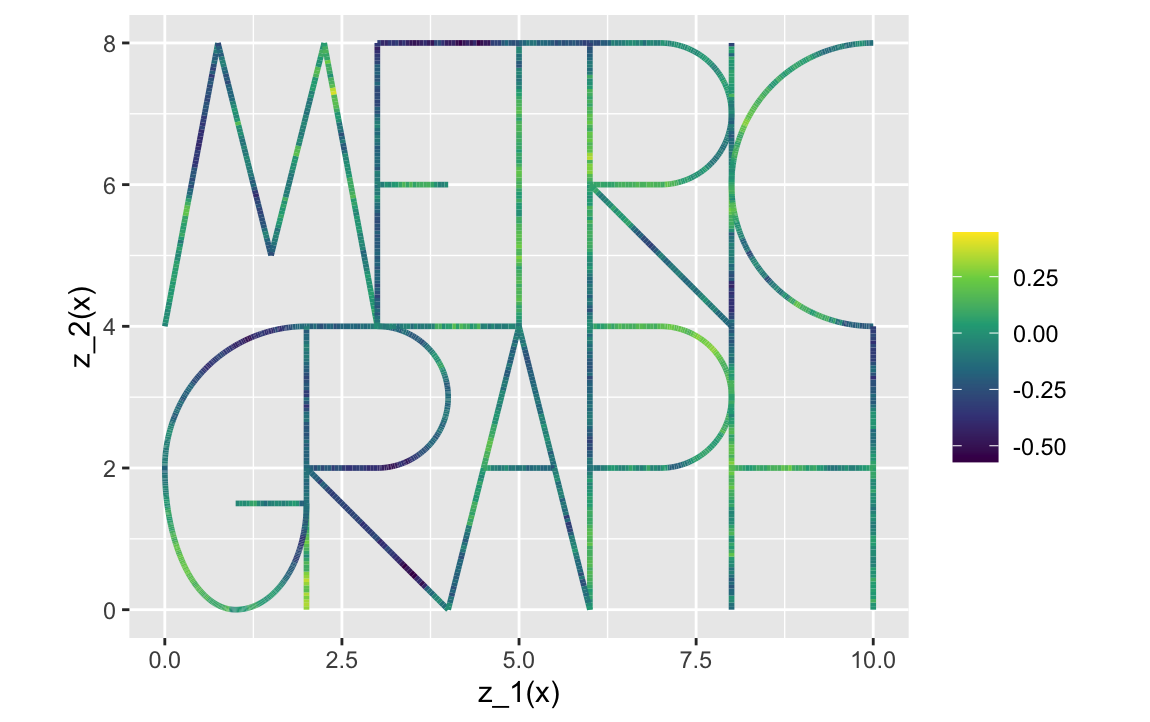}
    \caption{MAP estimate of $u_0$}
    \label{fig:mu-map}
  \end{subfigure}\hfill
  \begin{subfigure}[b]{0.5\textwidth}
    \includegraphics[width=\textwidth]{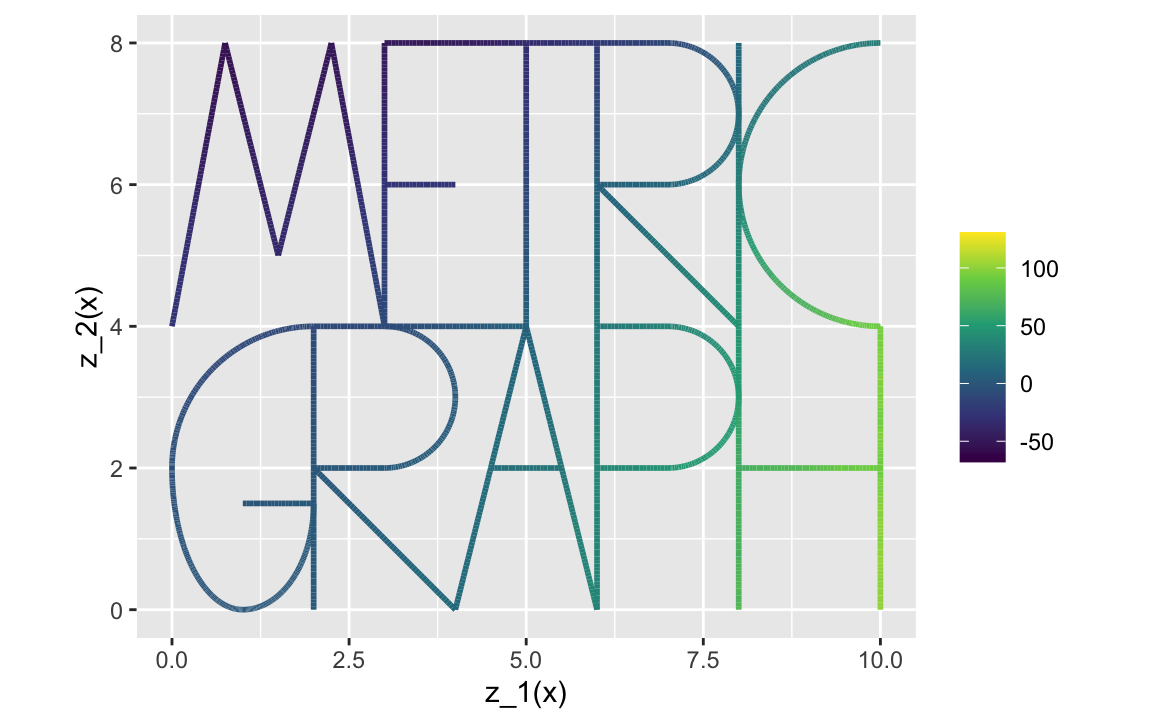}
    \caption{MAP estimate of $p_0$}
    \label{fig:p-map}
  \end{subfigure}

  \caption{Elliptic problem ($\beta =1).$ Comparison of ground truth, posterior mean, and MAP estimate for the parameter (left column) and for the corresponding PDE solution (right column).}
  \label{fig:mu-p-comparison}
\end{figure}

\begin{figure}[htbp]
  \centering
  \textbf{Elliptic problem: posterior marginal standard deviation}\\[3ex]
  \begin{subfigure}[b]{0.5\textwidth}
    \includegraphics[width=\textwidth]{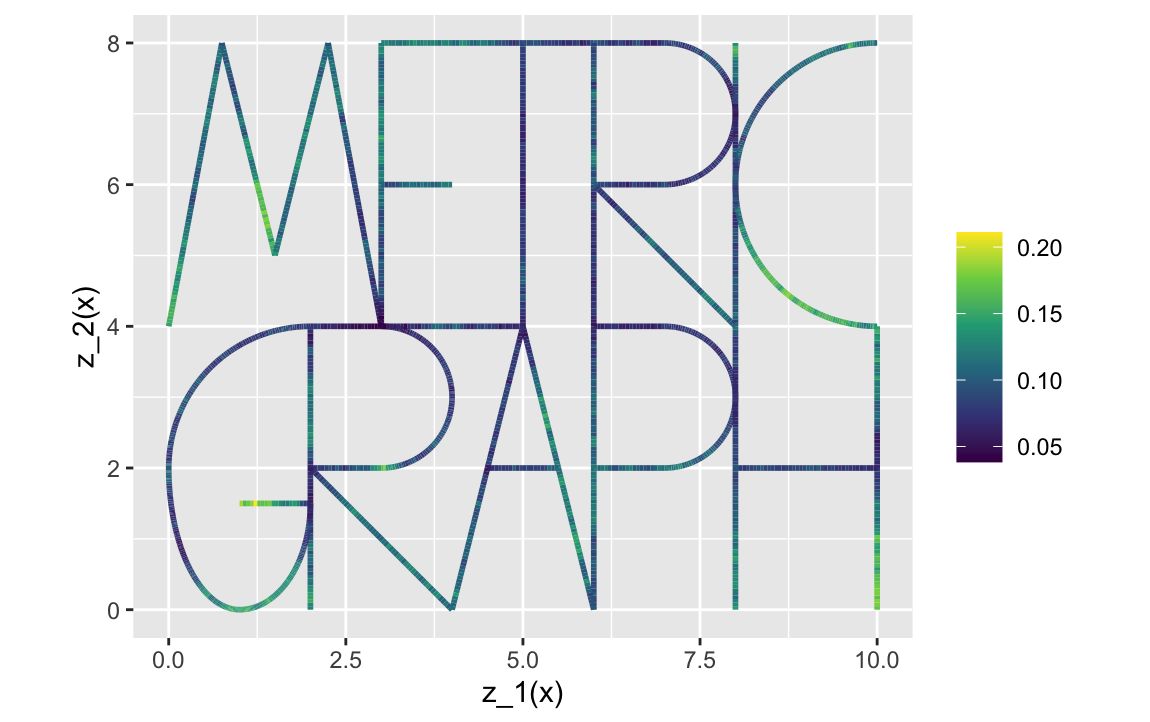}
    \caption{Standard deviation for parameter}
    \label{fig:var-u}
  \end{subfigure}\hfill
  \begin{subfigure}[b]{0.5\textwidth}
    \includegraphics[width=\textwidth]{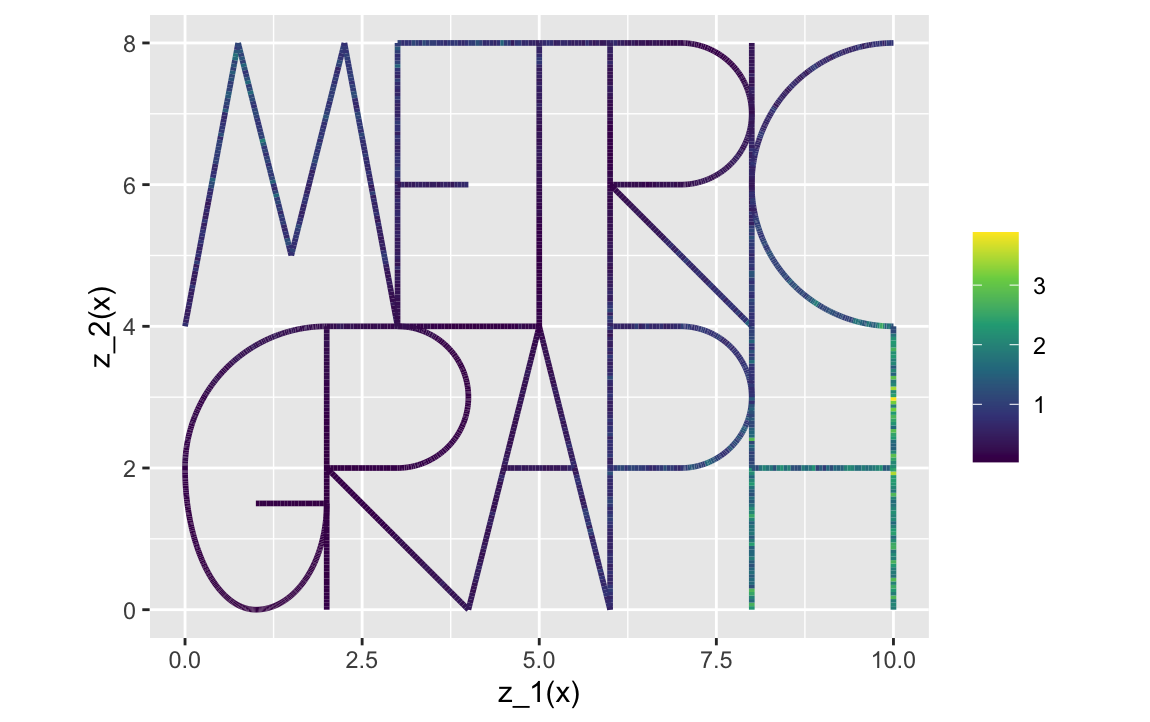}
    \caption{Standard deviation for PDE solution}
    \label{fig:var-p}
  \end{subfigure}

  \vspace{6ex}

  \textbf{Elliptic problem: absolute difference between truth and posterior mean}\\[3ex]
  \begin{subfigure}[b]{0.5\textwidth}
    \includegraphics[width=\textwidth]{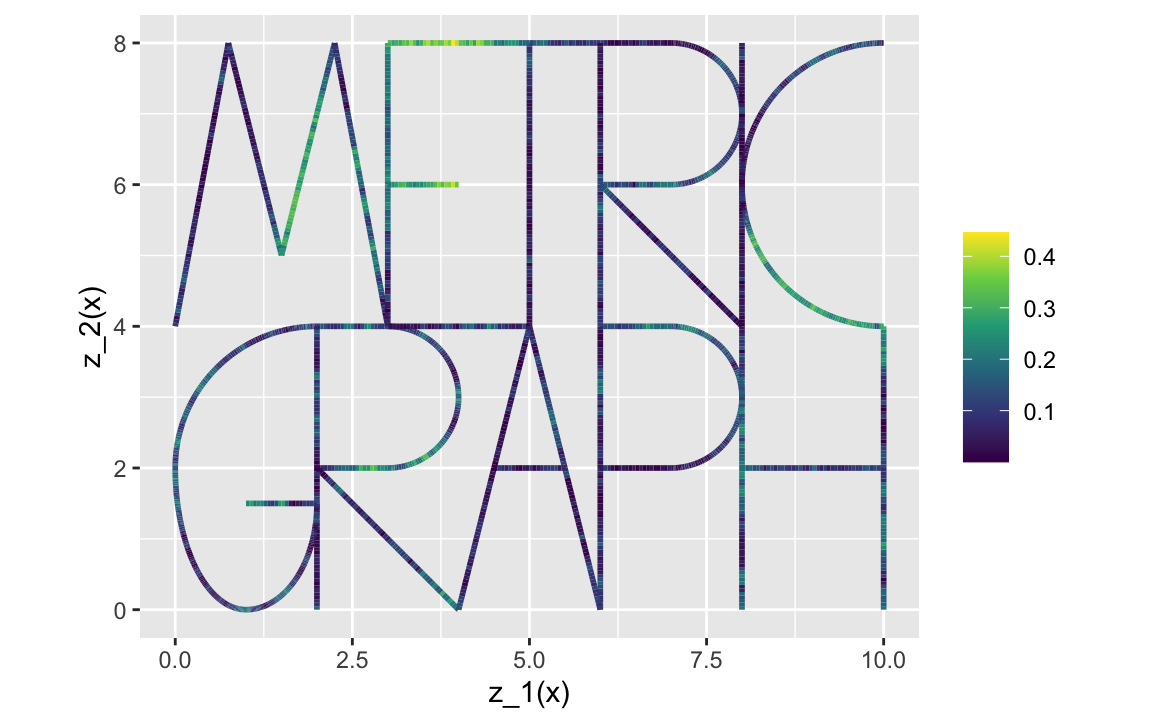}
    \caption{Difference between truth $u_0$ and posterior mean}
    \label{fig:u-postmean-difference}
  \end{subfigure}\hfill
  \begin{subfigure}[b]{0.5\textwidth}
    \includegraphics[width=\textwidth]{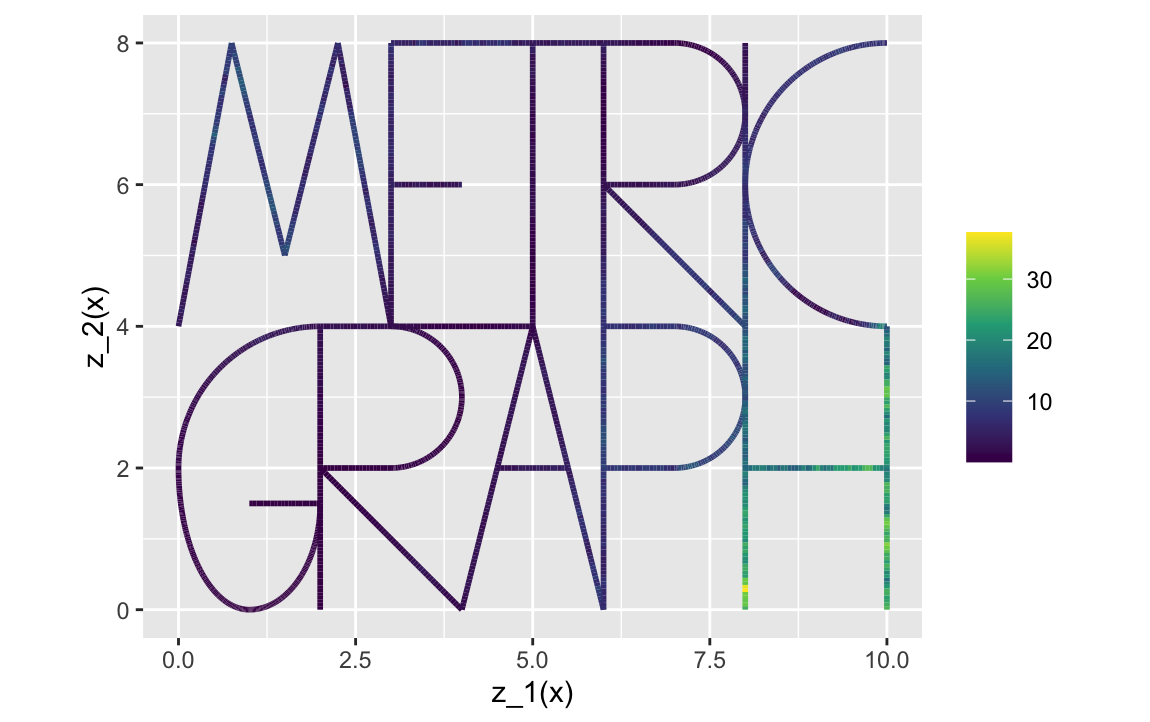}
    \caption{Difference between truth $p_0$ and posterior mean}
    \label{fig:p-postmean-difference}
  \end{subfigure}\hfill

  \caption{Elliptic case ($\beta =1)$. Top row:  posterior marginal standard deviation computed from MCMC samples for the parameter (left), \black defined as in Equation~\eqref{eq:post_sd}\nc, and for the corresponding PDE solution (right). Bottom row: difference between the truth and the posterior mean for the parameter (left) and the PDE solution (right).}
  \label{fig:variance-up}
\end{figure}

We begin by illustrating the performance of the posterior mean and posterior marginal mode (MAP-type) estimators \black on a \emph{single} representative realization of the inverse problem, corresponding to one draw of the ground-truth parameter $u_0\sim\mu_0$ and the associated state $p_0=p(u_0)$. \nc Figure~\ref{fig:mu-p-comparison} displays the posterior mean estimator \black $\hat u_{\rm mean}$ \nc and the posterior marginal mode estimator \black $\hat u_{\rm mode}$\nc alongside the ground-truth parameter $u_0$. Both estimators successfully recover the global structure of $u_0$, with only minor discrepancies at a few localized points. We also evaluate the corresponding plug-in state reconstructions \black $\hat p_{\rm mean}$ and $\hat p_{\rm mode}$\nc , which exhibit close agreement with the ground-truth solution $p_0$, except near the edges where $p_0$ attains large values. To assess uncertainty, Figure~\ref{fig:variance-up} shows the posterior marginal standard deviation \black $\widehat{\mathrm{sd}}_{u\mid y}$\nc. Regions with larger reconstruction errors coincide with elevated marginal standard deviations, confirming that the Bayesian approach correctly identifies areas of increased uncertainty.

\black
To assess average performance beyond a single realization, we repeat the entire experiment over multiple independently drawn ground truths $u_0\sim\mu_0$ and summarize trial-averaged reconstruction metrics (mean $\pm$ sd) for both point estimators. For the elliptic experiment ($\beta=1$ in Equation~\ref{equa:p}), we draw $30$ independent ground truths to illustrate the typical behavior of the reconstructions; the resulting metrics are reported in Table~\ref{tab:recon-metrics-standard}. Complementary box plots of these metrics, showing their distribution across trials, are provided in Appendix~\ref{app:boxplots-metrics}.\nc

\black
\begin{table}[h]
\centering
\small
\setlength{\tabcolsep}{6pt}
\renewcommand{\arraystretch}{1.15}
\begin{tabular}{@{}lcc@{}}
\toprule
\textbf{Metric} & \textbf{Posterior mean} & \textbf{Posterior marginal mode} \\
\midrule
$\mathrm{RMSE}(u_{\rm est}, u_0)$
& $0.0414 \pm 0.0107$ & $0.0430 \pm 0.0106$ \\

$\mathrm{RMSE}(p_{\rm est},p_0)$
& $0.402 \pm 0.124$ & $0.474 \pm 0.205$ \\

$\mathrm{NRMSE}_{\rm range}(u_{\rm est},u_0)$
& $0.183 \pm 0.0427$ & $0.190 \pm 0.0416$ \\

$\mathrm{RelErr}_{L^2(\Gamma)}(p_{\rm est},p_0)$
& $0.0104 \pm 0.00230$ & $0.0123 \pm 0.00469$ \\

$\mathrm{RelErr}_{L^2(\Gamma)}(a_{\rm est}, a_0)$
& $0.0412 \pm 0.0107$ & $0.0427 \pm 0.0106$ \\
\bottomrule
\end{tabular}
\caption{\black Elliptic experiment ($\beta=1$): reconstruction metrics (mean $\pm$ sd across $30$ independent trials) comparing the posterior mean estimator with the posterior marginal mode estimator. Here $a=\exp(u)$ and $\mathrm{RelErr}_{L^2(\Gamma)}$ denotes the discrete $L^2(\Gamma)$ relative error.\nc}
\label{tab:recon-metrics-standard}
\end{table}
\nc


\begin{table}[htp]
\centering
\small
\setlength{\tabcolsep}{6pt}
\renewcommand{\arraystretch}{1.15}
\begin{tabular}{@{}lcc@{}}
\toprule
\textbf{Metric} & \textbf{Posterior mean} & \textbf{Posterior marginal mode} \\
\midrule
$\mathrm{RMSE}(u_{\rm est},u_0)$ & $0.0293 \pm 0.00355$ & $0.0309 \pm 0.00385$ \\
$\mathrm{RMSE}(p_{\rm est},p_0)$ & $0.855 \pm 0.160$ & $0.882 \pm 0.160$ \\
$\mathrm{NRMSE}_{\rm range}(u_{\rm est},u_0)$ & $0.135 \pm 0.0202$ & $0.143 \pm 0.0226$ \\
$\mathrm{RelErr}_{L^2(\Gamma)}(p_{\rm est},p_0)$ & $0.0100 \pm 0.00168$ & $0.0103 \pm 0.00160$ \\
$\mathrm{RelErr}_{L^2(\Gamma)}(a_{\rm est},a_0)$ & $0.0285 \pm 0.00351$ & $0.0298 \pm 0.00379$ \\
\bottomrule
\end{tabular}
\caption{\black Fractional elliptic experiment ($\beta=3/2$): reconstruction metrics (mean $\pm$ sd across $10$ independent trials) comparing the posterior mean estimator with the posterior marginal mode estimator. Here $a=\exp(u)$ and $\mathrm{RelErr}_{L^2(\Gamma)}$ denotes the discrete $L^2(\Gamma)$ relative error.\nc}
\label{tab:recon-metrics-fractional}
\end{table}

\subsection{Fractional Elliptic Problem on a Letter-shaped Metric Graph}
\label{ssec:output fractional case}

\begin{figure}[htbp]
  \centering
  \textbf{Fractional elliptic problem: estimates for parameter and PDE solution}\\[3ex]
  \begin{subfigure}[b]{0.5\textwidth}
    \includegraphics[width=\textwidth]{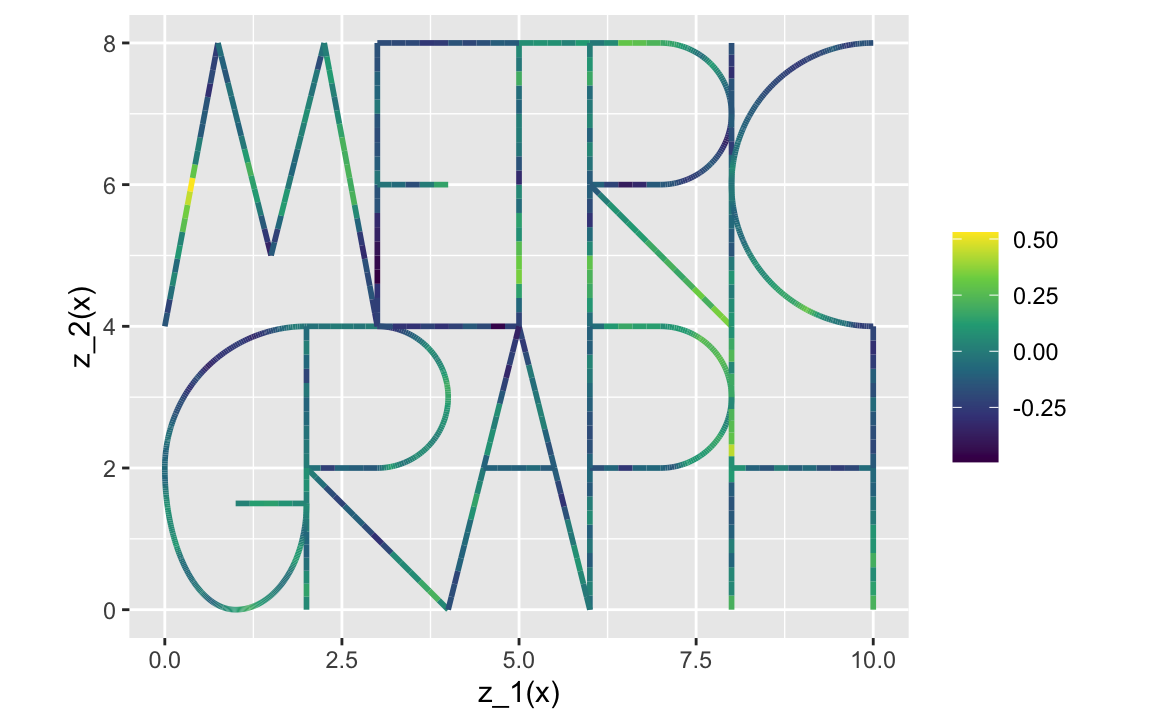}
    \caption{True parameter $u_0$}
    \label{fig:frac mu-true}
  \end{subfigure}\hfill
  \begin{subfigure}[b]{0.5\textwidth}
    \includegraphics[width=\textwidth]{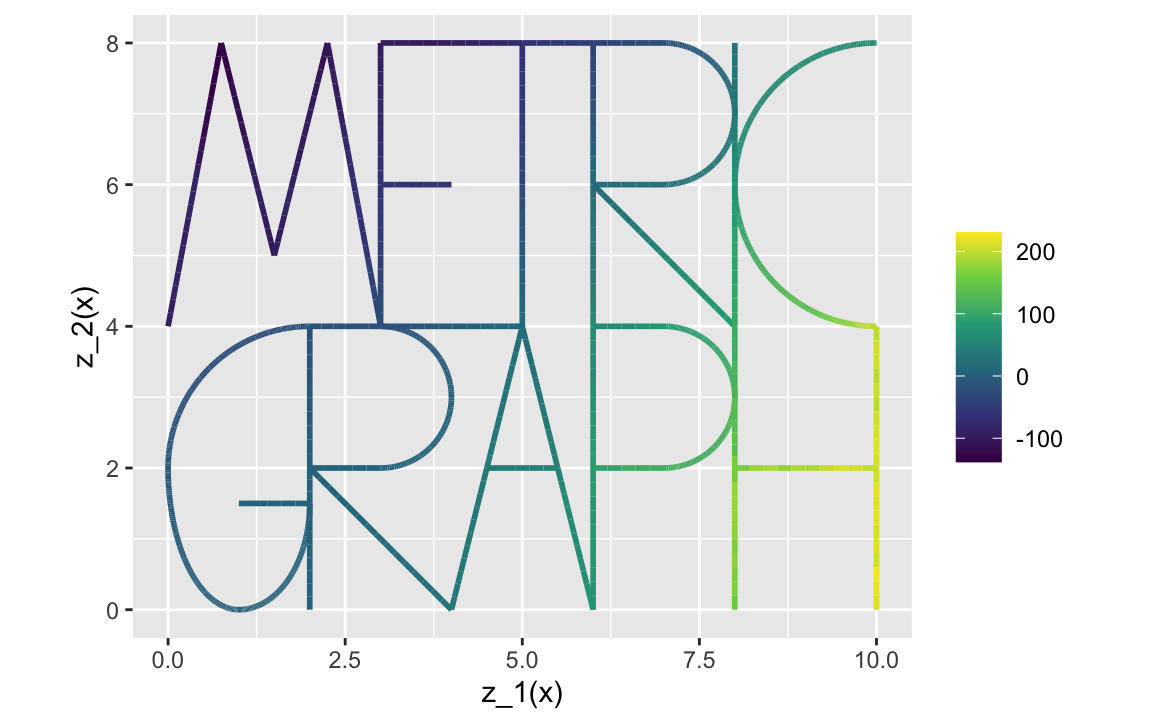}
    \caption{True solution $p_0$}
    \label{fig:frac p-true}
  \end{subfigure}

  \vspace{4ex}

  \begin{subfigure}[b]{0.5\textwidth}
    \includegraphics[width=\textwidth]{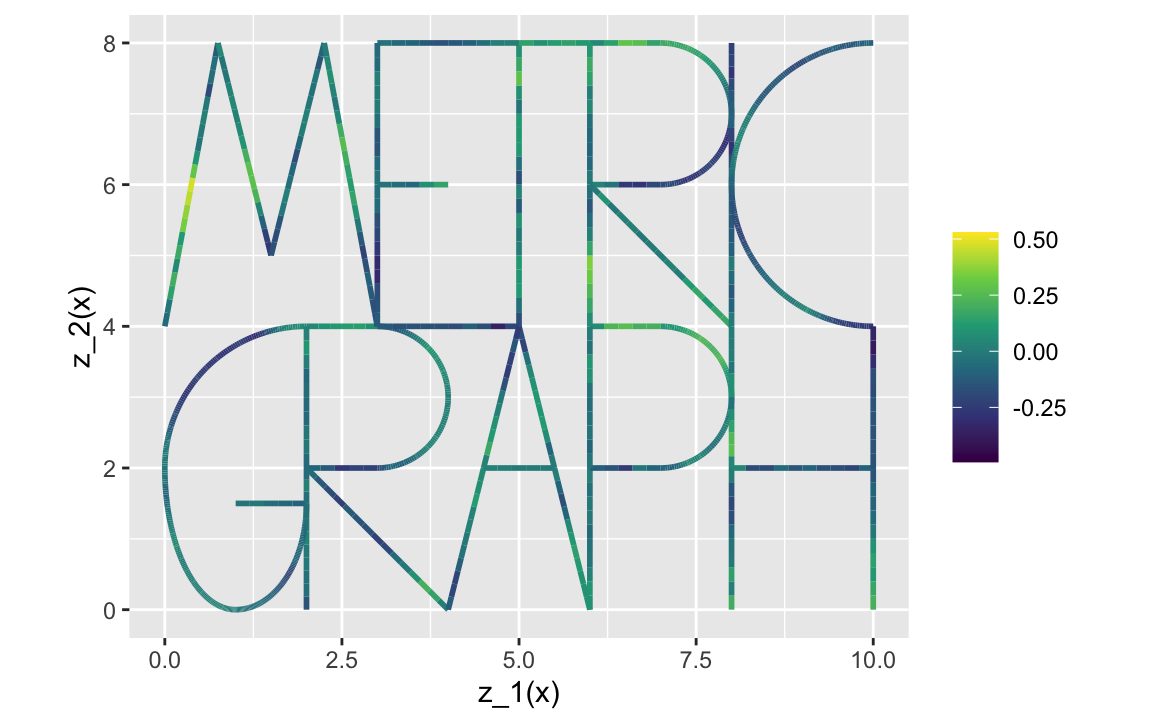}
    \caption{Posterior mean estimate of $u_0$}
    \label{fig:frac mu-postmean}
  \end{subfigure}\hfill
  \begin{subfigure}[b]{0.5\textwidth}
    \includegraphics[width=\textwidth]{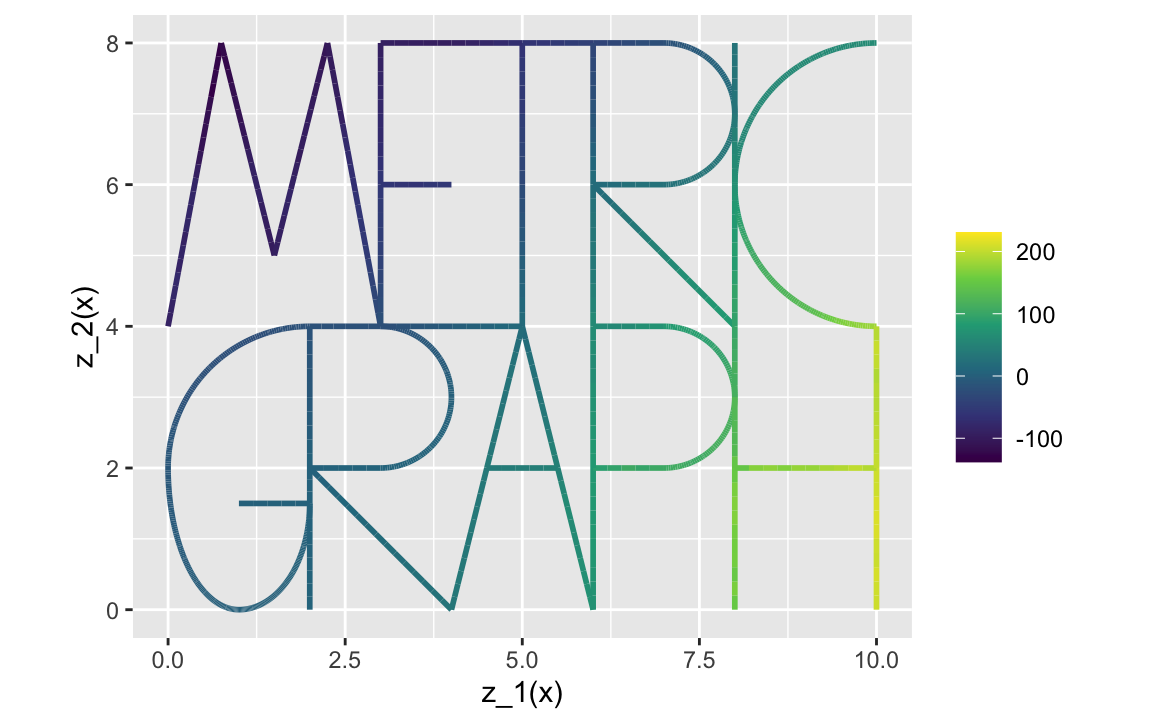}
    \caption{Posterior mean estimate of $p_0$}
    \label{fig:frac p-postmean}
  \end{subfigure}

  \vspace{4ex}

  \begin{subfigure}[b]{0.5\textwidth}
    \includegraphics[width=\textwidth]{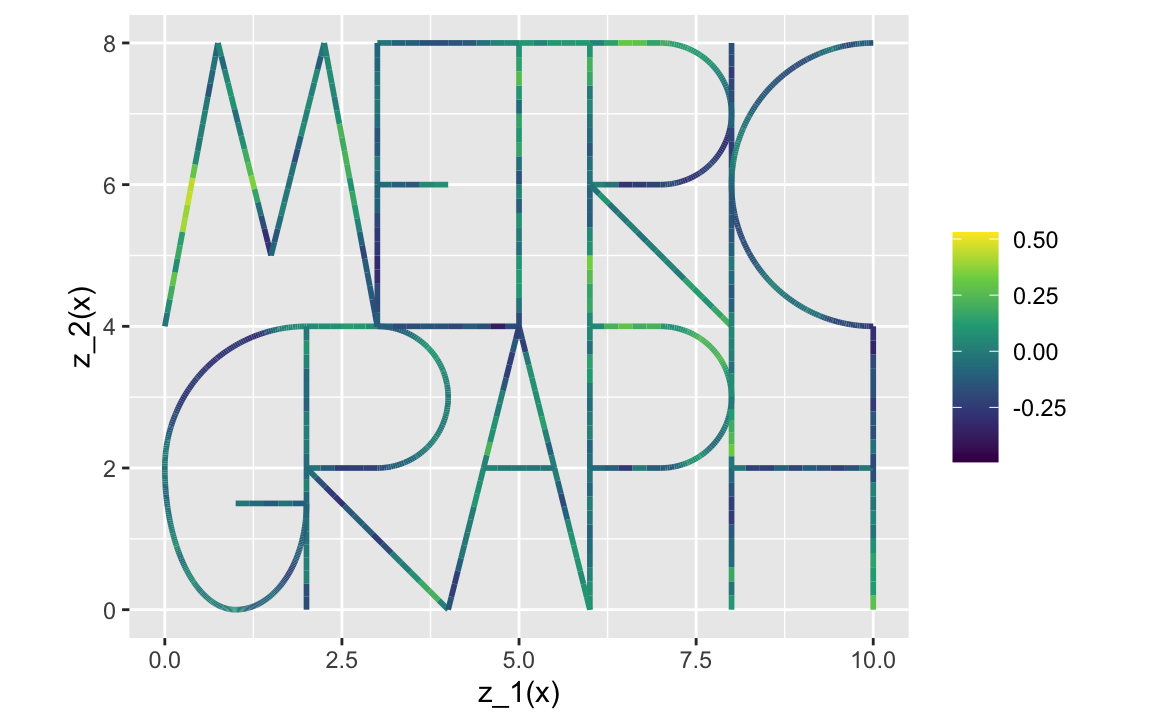}
    \caption{MAP estimate of $u_0$}
    \label{fig:frac mu-map}
  \end{subfigure}\hfill
  \begin{subfigure}[b]{0.5\textwidth}
    \includegraphics[width=\textwidth]{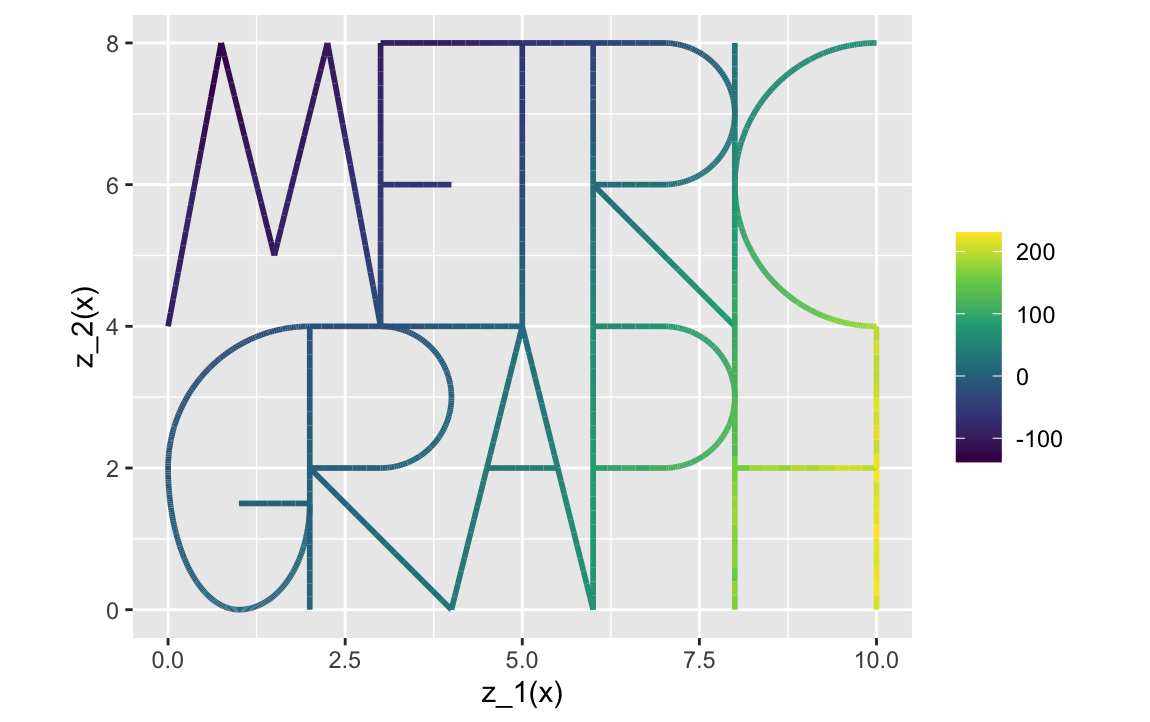}
    \caption{MAP estimate of $p_0$}
    \label{fig:frac p-map}
  \end{subfigure}

  \caption{Fractional elliptic problem ($\beta = 3/2).$ Comparison of ground truth, posterior mean, and MAP estimate for the parameter (left column) and for the corresponding PDE solution  (right column).}
  \label{fig: fractional mu-p-comparison}
\end{figure}

\begin{figure}[htbp]
  \centering
  \textbf{Fractional elliptic problem: posterior marginal standard deviation}\\[3ex]
  \begin{subfigure}[b]{0.5\textwidth}
    \includegraphics[width=\textwidth]{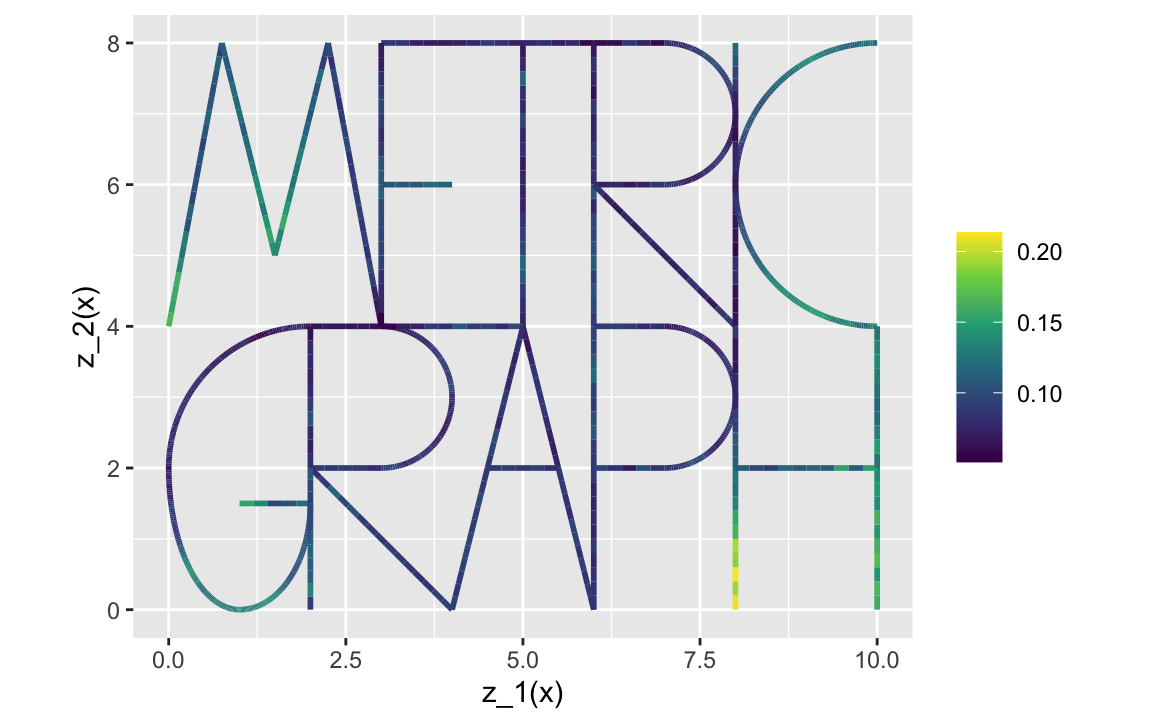}
    \caption{Standard deviation for parameter}
    \label{fig:frac var-u}
  \end{subfigure}\hfill
  \begin{subfigure}[b]{0.50\textwidth}
   \includegraphics[width=\textwidth]{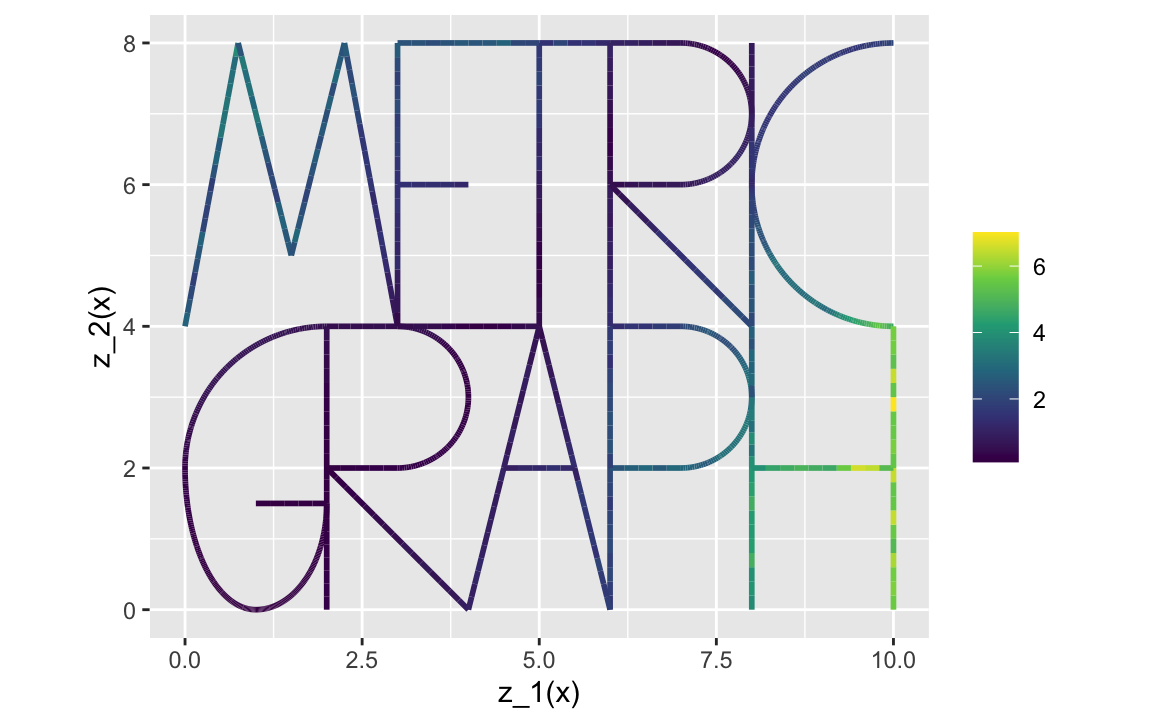}
    \caption{Standard deviation for PDE solution}
    \label{fig:frac var-p}
  \end{subfigure}

  \vspace{6ex}

  \textbf{Fractional elliptic problem: absolute difference between truth and posterior mean}\\[3ex]
  \begin{subfigure}[b]{0.5\textwidth}
    \includegraphics[width=\textwidth]{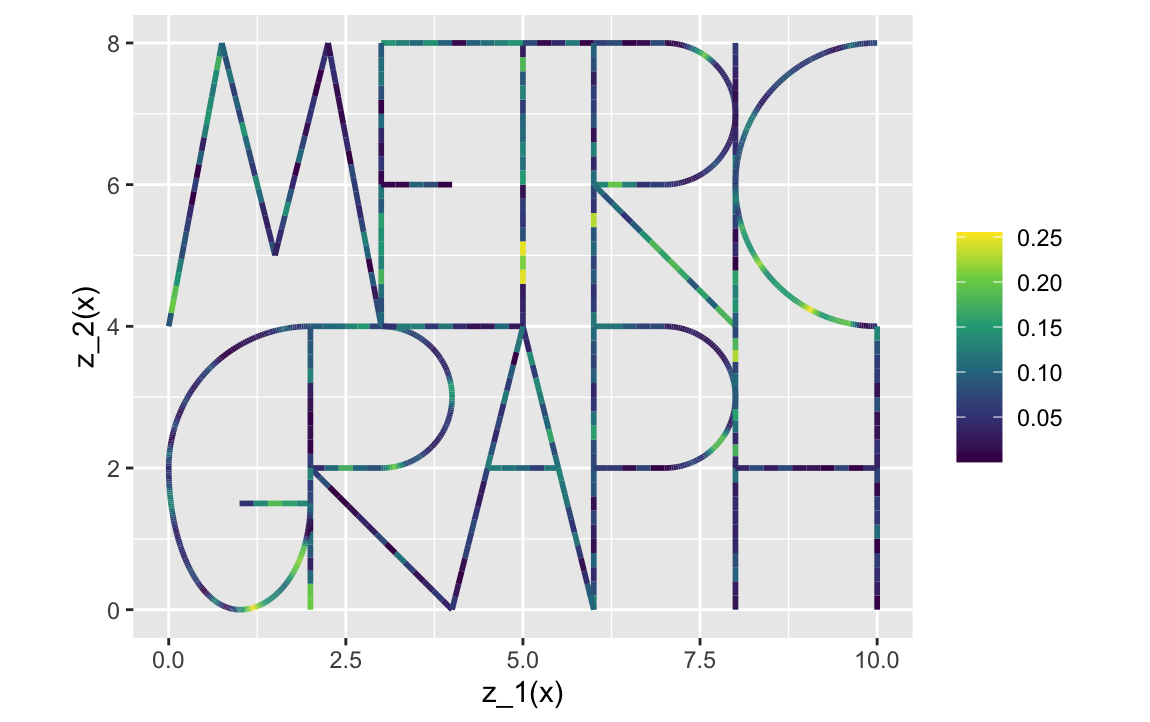}
    \caption{Difference between truth $u_0$ and posterior mean}
    \label{fig:u-postmean-difference-frac}
  \end{subfigure}\hfill
  \begin{subfigure}[b]{0.5\textwidth}
    \includegraphics[width=\textwidth]{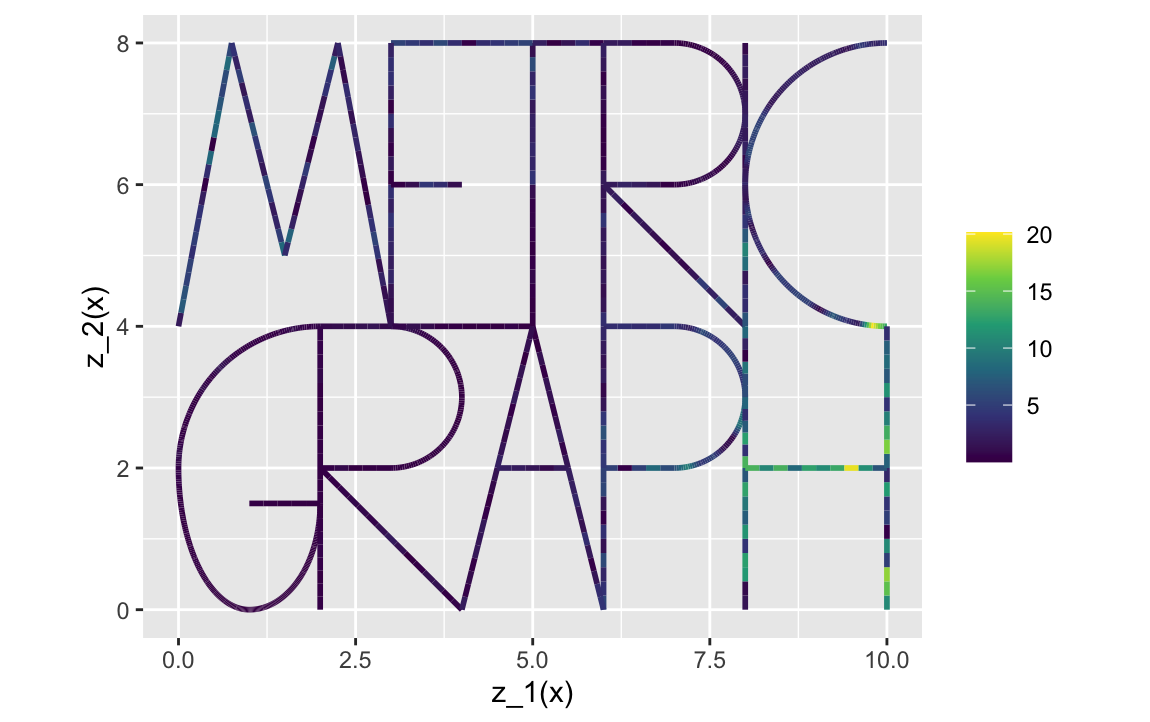}
    \caption{Difference between truth $p_0$ and posterior mean}
    \label{fig:p-postmean-difference-frac}
  \end{subfigure}\hfill

  \caption{Fractional elliptic problem ($\beta = 3/2)$. Top row: posterior marginal standard deviation  computed from MCMC samples for the parameter $u$ (left) and for the corresponding PDE solution $p$ (right). Bottom row: difference between the truth and the posterior mean for the parameter (left) and the PDE solution (right).}
  
  \label{fig:frac variance-up}
\end{figure}

For the fractional elliptic case, we set the order parameter $\beta = 3/2$ in Equation~\ref{equa:p} and follow the same workflow as in Subsection~\ref{ssec:output standard case}. Figure~\ref{fig: fractional mu-p-comparison} shows strong agreement between the posterior mean estimator and \black the posterior marginal mode \nc estimator and the ground-truth parameter, with the $\mathrm{RMSE}(u_{\rm est},u_0)$ falling below $0.08$ \black in this representative realization\nc. Figure~\ref{fig:frac variance-up} displays the posterior marginal standard deviation computed from the MCMC samples; as in the elliptic case, regions with higher marginal standard deviations align with regions of larger reconstruction errors.

\black
Here, we draw $10$ independent ground truths and report the corresponding metrics in Table~\ref{tab:recon-metrics-fractional}. We use fewer trials in this case because each fractional forward solve (and hence each MCMC run) is considerably more time consuming, and because the resulting reconstruction metrics exhibit consistently small standard deviations, suggesting that $10$ trials are sufficient to characterize the typical performance in this setting. Complementary box plots illustrating the distribution of these metrics across trials are provided in Appendix~\ref{app:boxplots-metrics}.
\nc

\section{Conclusions}\label{sec:conclusions}
This paper has studied the formulation, well-posedness, and numerical solution of Bayesian inverse problems on metric graphs, focusing on elliptic and fractional elliptic inverse problems. We have leveraged recent Gaussian process models on metric graphs to specify the prior, and we have built upon recent regularity theory for PDEs on metric graphs to establish the stability of the forward model. Numerical results demonstrate accurate reconstruction and effective uncertainty quantification. 

Several research directions stem from this work, including: (1) Analyzing Bayesian inversion for other PDEs on metric graphs, which will likely require new regularity theory and stability estimates. (2) Developing prior models on metric graphs  \cite{bolin2025log}, including for instance level-set techniques for geometric inverse problems \cite{iglesias2016bayesian}. (3) Studying the effect of numerical discretization of the forward model and the prior in the posterior \cite{sanz2024analysis,sanz2022finite}; doing so will motivate the design and analysis of new solvers for PDEs on metric graphs.     

\black
\section*{Data Availability}
The code that supports the findings of this study is openly available on Zenodo at \url{https://doi.org/10.5281/zenodo.18714824}. The development repository is hosted at \url{https://github.com/WenwenLi2002/Bayesian-Inverse-Problems-on-Metric-Graphs}.
\nc

\section*{Acknowledgments}
DB was partially funded by King Abdullah University of Science and Technology (KAUST) under Award No.~ORFS-CRG12-2024-6399. DSA was partly funded by NSF CAREER DMS-2237628.



\bibliographystyle{siam}
\bibliography{references}

\newpage
\appendix

\begingroup
\color{black}

\section{Appendix}
\label{app:additional-results}
\subsection{Box Plots of Reconstruction Metrics}
\label{app:boxplots-metrics}
\begin{figure}[H]
    \centering

    \begin{subfigure}[t]{0.32\textwidth}
        \centering
        \includegraphics[width=\linewidth]{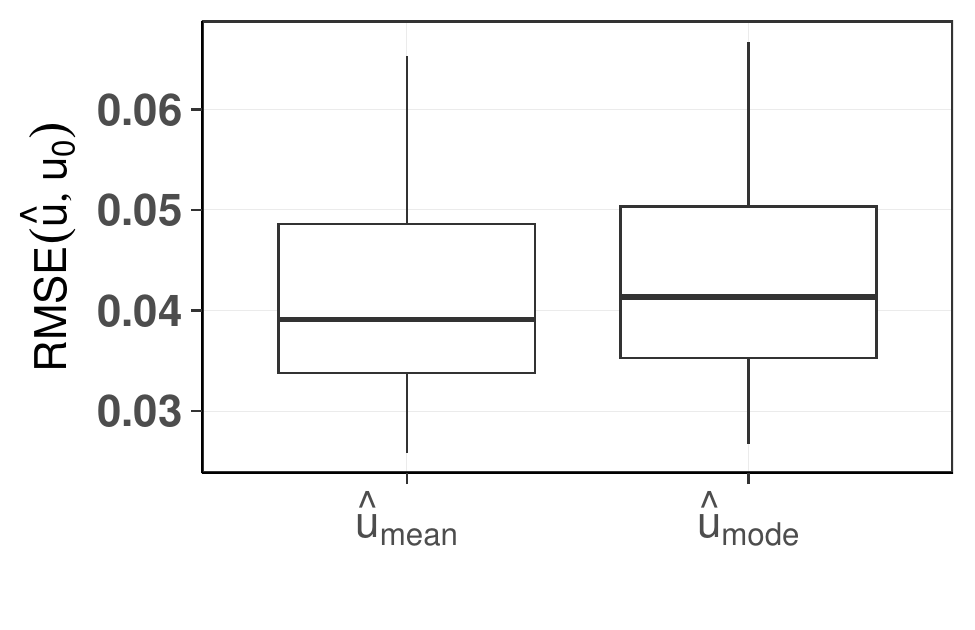}
        \caption{$\mathrm{RMSE}(\hat u,u_0)$}
        \label{fig:box-rmse-u}
    \end{subfigure}
    \hfill
    \begin{subfigure}[t]{0.32\textwidth}
        \centering
        \includegraphics[width=\linewidth]{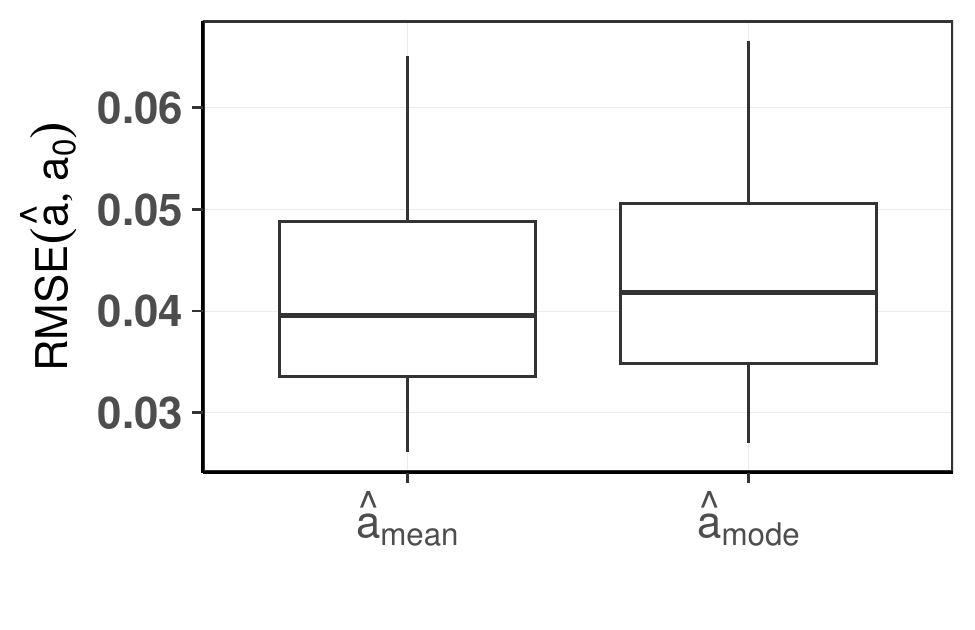}
        \caption{$\mathrm{RMSE}(\hat a,a_0)$}
        \label{fig:box-rmse-a}
    \end{subfigure}
    \hfill
    \begin{subfigure}[t]{0.32\textwidth}
        \centering
        \includegraphics[width=\linewidth]{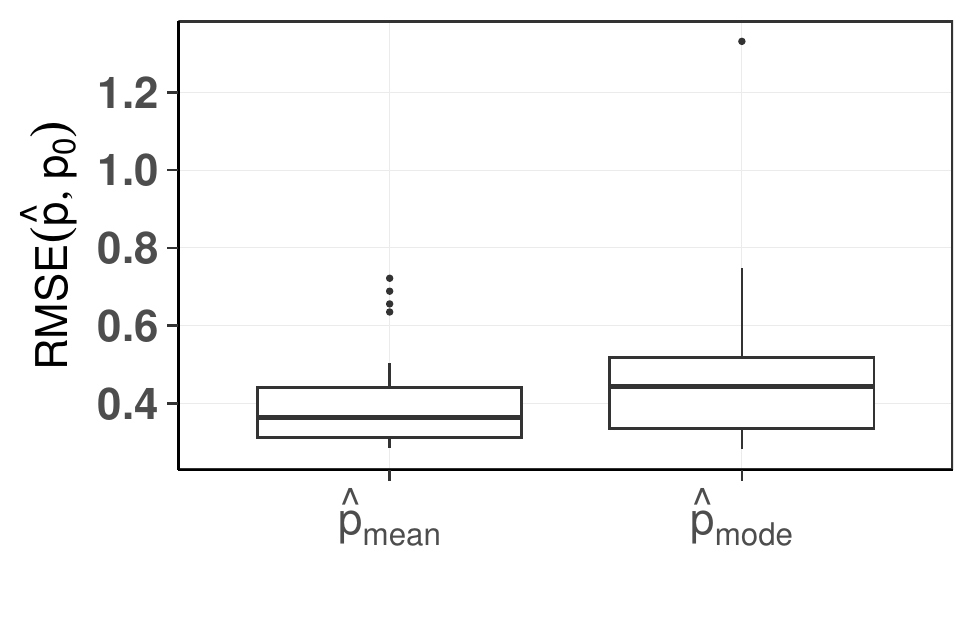}
        \caption{$\mathrm{RMSE}(\hat p,p_0)$}
        \label{fig:box-rmse-p}
    \end{subfigure}
    \vspace{0.25cm}

    \begin{subfigure}[t]{0.32\textwidth}
        \centering
        \includegraphics[width=\linewidth]{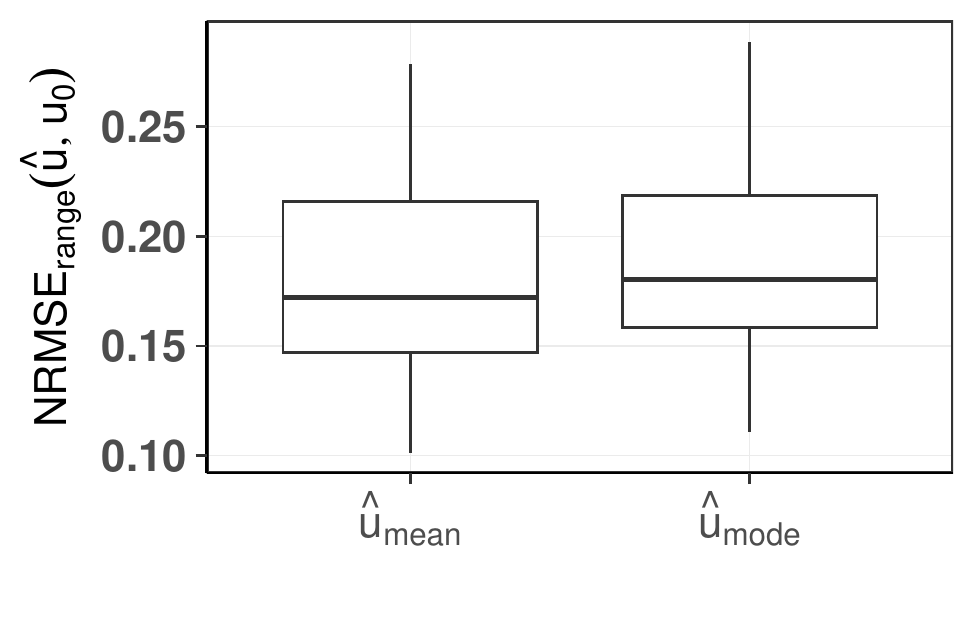}
        \caption{$\mathrm{NRMSE}_{\rm range}(\hat u,u_0)$}
        \label{fig:box-nrmse-u-range}
    \end{subfigure}
    \hfill
    \begin{subfigure}[t]{0.32\textwidth}
        \centering
        \includegraphics[width=\linewidth]{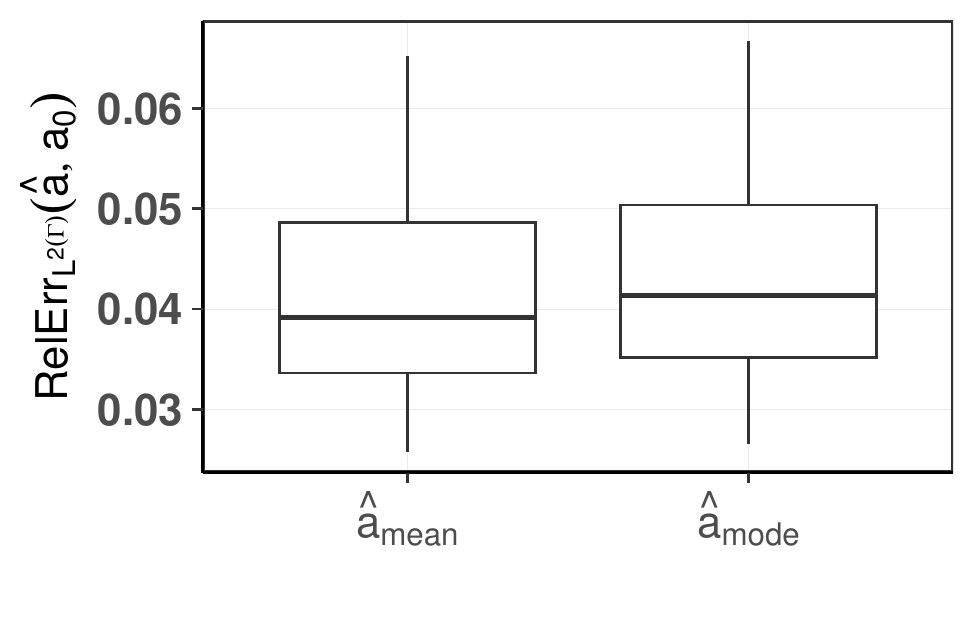}
        \caption{$\mathrm{RelErr}_{L^2(\Gamma)}(\hat a,a_0)$}
        \label{fig:box-relerr-a}
    \end{subfigure}
    \hfill
    \begin{subfigure}[t]{0.32\textwidth}
        \centering
        \includegraphics[width=\linewidth]{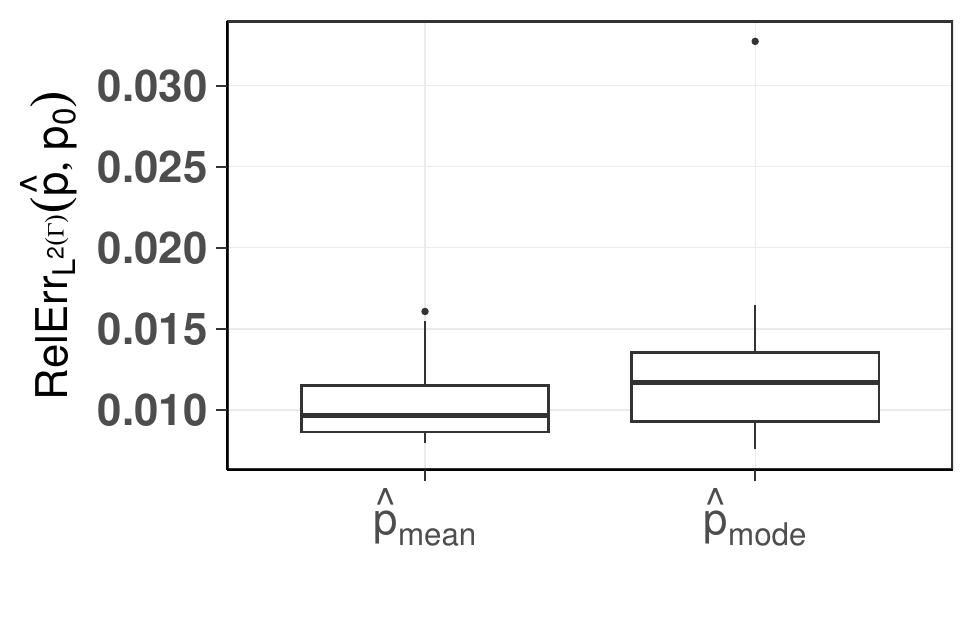}
        \caption{$\mathrm{RelErr}_{L^2(\Gamma)}(\hat p,p_0)$}
        \label{fig:box-relerr-p}
    \end{subfigure}

    \caption{ \black
    Boxplots of reconstruction errors over repeated trials.
    Top row: RMSE for the log-conductivity estimator $\hat u$ and the corresponding plug-in state estimator $\hat p=p(\hat u)$.
    Bottom row: range-normalized error for $u$ and relative $L^2(\Gamma)$ errors (mass-matrix based) for the conductivity $\hat a=\exp(\hat u)$ and the state $\hat p$.\nc
    }
    \label{fig:boxplots-2plus3-errors}
\end{figure}

\begin{figure}[H]
    \centering

    \begin{subfigure}[t]{0.32\textwidth}
        \centering
        \includegraphics[width=\linewidth]{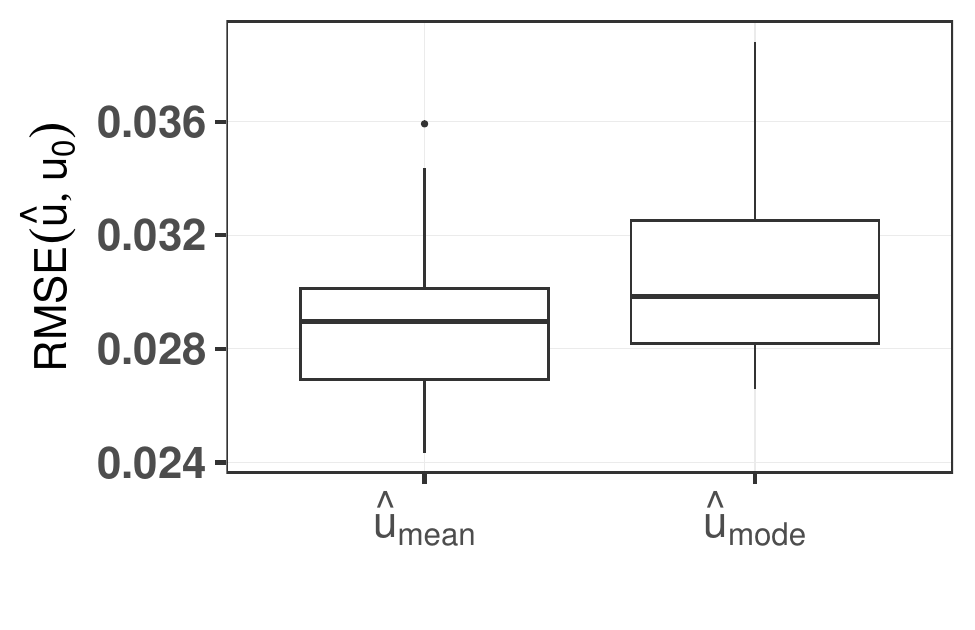}
        \caption{$\mathrm{RMSE}(\hat u,u_0)$}
        \label{fig:box-rmse-u-frac}
    \end{subfigure}
    \hfill
    \begin{subfigure}[t]{0.32\textwidth}
        \centering
        \includegraphics[width=\linewidth]{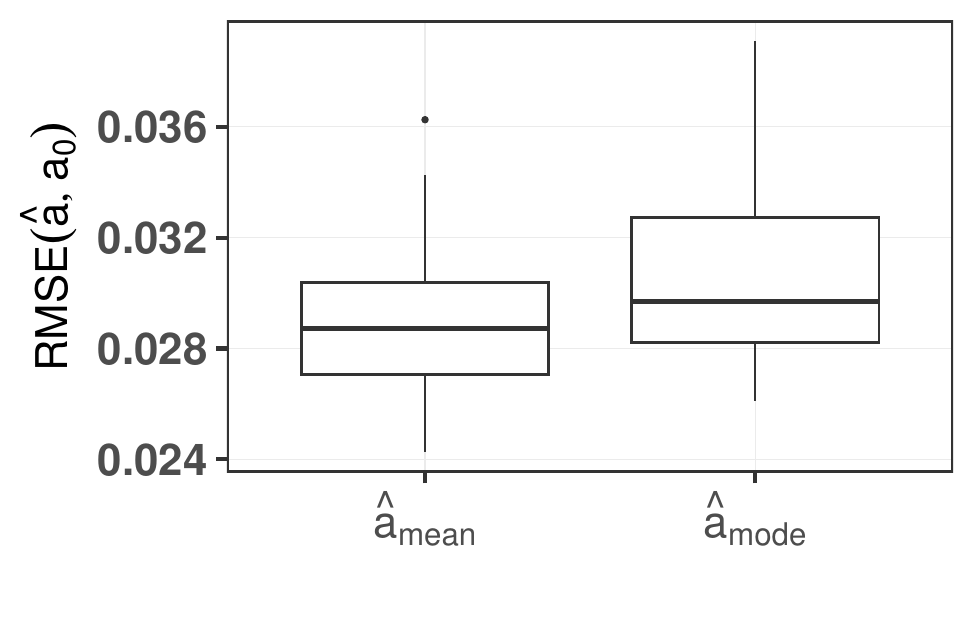}
        \caption{$\mathrm{RMSE}(\hat p,p_0)$}
        \label{fig:box-rmse-a-frac}
    \end{subfigure}
    \hfill
    \begin{subfigure}[t]{0.32\textwidth}
        \centering
        \includegraphics[width=\linewidth]{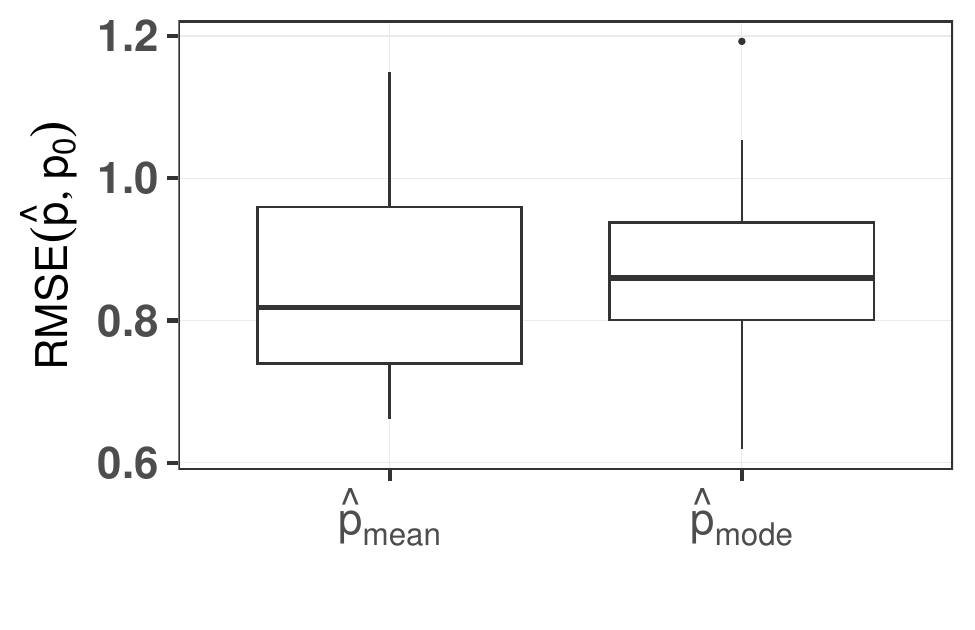}
        \caption{$\mathrm{RMSE}(\hat p,p_0)$}
        \label{fig:box-rmse-p-frac}
    \end{subfigure}


    \begin{subfigure}[t]{0.32\textwidth}
        \centering
        \includegraphics[width=\linewidth]{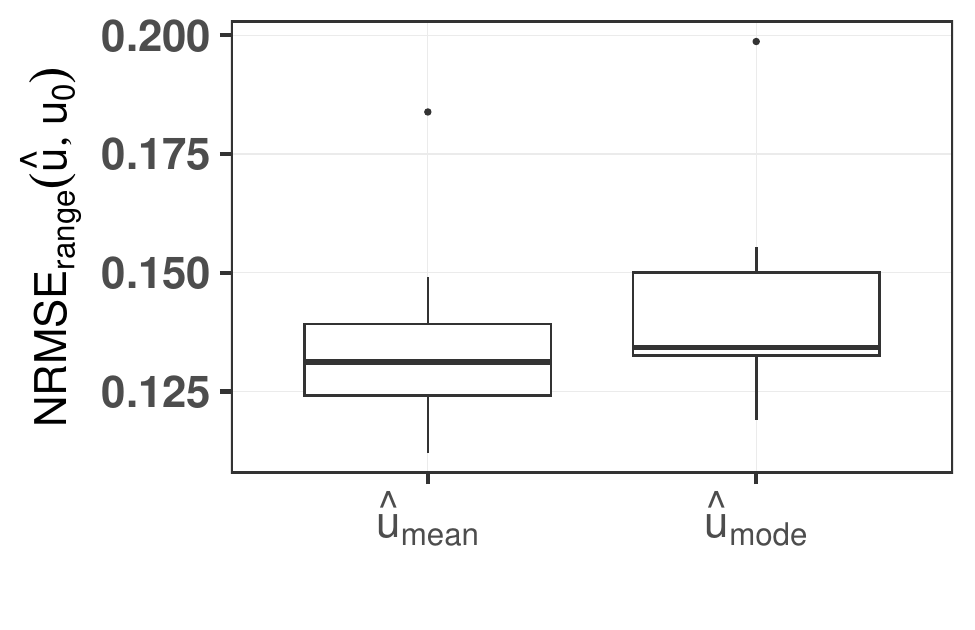}
        \caption{$\mathrm{NRMSE}_{\rm range}(\hat u,u_0)$}
        \label{fig:box-nrmse-u-range-frac}
    \end{subfigure}
    \hfill
    \begin{subfigure}[t]{0.32\textwidth}
        \centering
        \includegraphics[width=\linewidth]{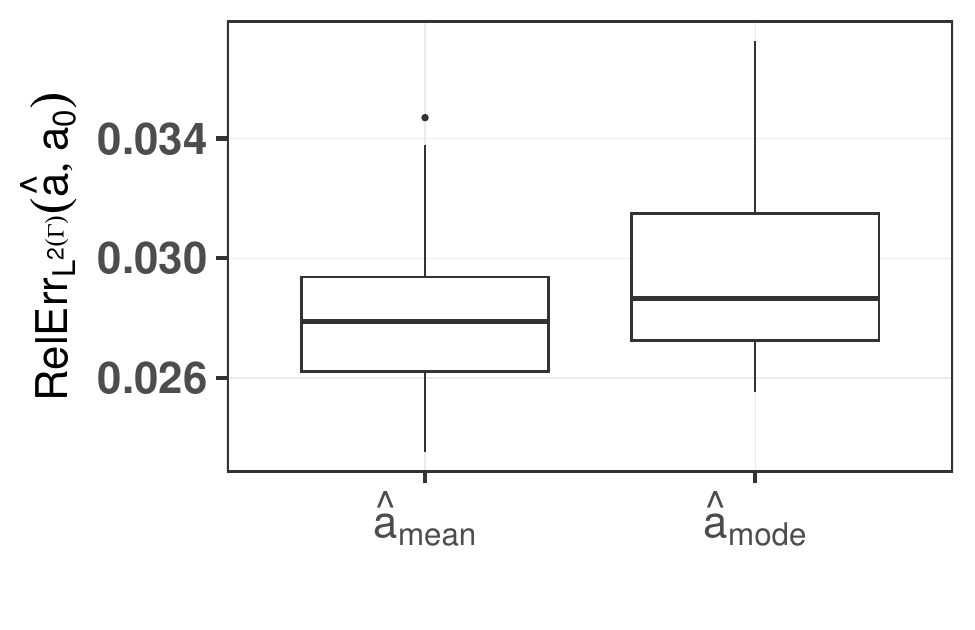}
        \caption{$\mathrm{RelErr}_{L^2(\Gamma)}(\hat a,a_0)$}
        \label{fig:box-relerr-a-frac}
    \end{subfigure}
    \hfill
    \begin{subfigure}[t]{0.32\textwidth}
        \centering
        \includegraphics[width=\linewidth]{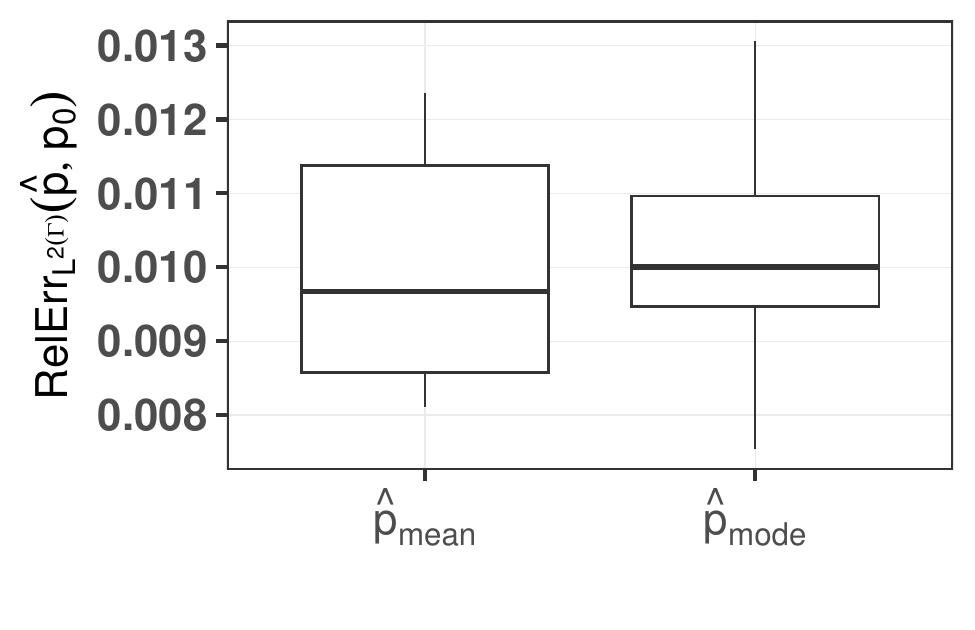}
        \caption{$\mathrm{RelErr}_{L^2(\Gamma)}(\hat p,p_0)$}
        \label{fig:box-relerr-p-frac}
    \end{subfigure}

    \caption{
    \black Fractional elliptic case: Boxplots of reconstruction errors over repeated trials. The reported metrics and the display convention are the same as in Figure~\ref{fig:boxplots-2plus3-errors}.\nc
    }
    \label{fig:boxplots-2plus3-errors-frac}
\end{figure}

Figures~\ref{fig:boxplots-2plus3-errors}--\ref{fig:boxplots-2plus3-errors-frac} visualize the trial-to-trial variability of the reconstruction errors via box plots, and complement the mean$\pm$sd summaries reported earlier (cf. Tables~\ref{tab:recon-metrics-standard} and \ref{tab:recon-metrics-fractional}).  
Overall, all metrics are well-behaved: the box plots are relatively tight with moderate whiskers, indicating that the reconstructions are consistently accurate and do not exhibit systematic instability across trials. 
Within this generally strong performance, a clear and consistent pattern emerges: posterior mean summaries ($\hat{u}_{\rm mean}$) tend to be slightly more stable than posterior marginal modes ($\hat{u}_{\rm mode}$). 
Across all panels, the ``mean'' estimators typically have comparable (often slightly smaller) medians and noticeably smaller interquartile ranges, together with shorter upper whiskers, suggesting reduced trial-to-trial variability. 
This effect is most visible in the state reconstruction, where the marginal-mode box plots show heavier upper tails and occasional high-error outliers, particularly for $\mathrm{RMSE}(\hat p_{\star},p_0)$, whereas the posterior-mean box plots remain more concentrated. 
For the coefficient reconstructions, $\mathrm{RMSE}(\hat u_\star,u_0)$ and the plug-in conductivity error $\mathrm{RMSE}(\hat a,a_0)$ are uniformly small and tightly clustered in both cases, with the posterior mean again exhibiting slightly tighter dispersion. 
The bottom-row diagnostics convey the same message: the range-normalized error for $u$ and the relative $L^2(\Gamma)$ errors for $a$ and $p$ are tightly distributed, and the posterior mean summaries are consistently less variable. 
A plausible explanation is that $\hat u_{\rm mode}$ is defined as a posterior \emph{marginal} mode, which may be less smooth and can exhibit localized irregularities. Such features may only have a limited impact on $\mathrm{RMSE}(\hat u_\star,u_0)$, yet they can propagate through the plug-in map $p(\hat u)$ and lead to noticeably larger errors in the recovered state.
Taken together, these box plots confirm that the reconstructions perform well across trials, while also highlighting a modest but systematic robustness advantage of posterior mean summaries over marginal modes.

\subsection{Numerical Verification of Well-posedness}
\label{app:wp-numerical}

We complement the baseline reconstruction metrics with a \emph{data-perturbation sensitivity test},
aimed at empirically probing the stability of posterior summaries with respect to perturbations of the observed data.
For each trial we start from a baseline dataset $y^0$ and generate perturbed observations $y^\delta$,
rerun the full inference pipeline, and compute posterior point summaries
$\hat u^\delta_\star$ and the associated fields
$\hat p^\delta_\star := p(\hat u^\delta_\star)$ and $\hat a^\delta_\star := \exp(\hat u^\delta_\star)$,
where $\star\in\{\mathrm{mean},\mathrm{mode}\}$ denotes the posterior mean and the posterior marginal mode
(pointwise marginal MAP) summaries.

Throughout this test we \emph{fix} the observation noise model (cf. Equation~\ref{eq:noise model}) used in the likelihood.
Specifically, given a parameter field $u$ and the corresponding solution $p(u)$, we set the
data-dependent standard deviation at an observation location $x_i$ to be
\[
\sigma_i(u)\;:=\;\alpha_{\mathrm{noise}}\bigl|p(u)(x_i)\bigr|+\beta_{\mathrm{noise}},
\]
with $\alpha_{\mathrm{noise}},\beta_{\mathrm{noise}}>0$ held fixed across all perturbation levels.
The perturbation level $\delta\ge 0$ is introduced \emph{only} through the data-generation step:
letting $p^\dagger$ denote the forward solution under the ground truth $u^\dagger$ and defining
$\sigma_i^\dagger := \alpha_{\mathrm{noise}}|p^\dagger(x_i)|+\beta_{\mathrm{noise}}$, we generate
\[
y_i^\delta \;=\; y_i^0 + \delta\,\sigma_i^\dagger\,\xi_i,
\qquad \xi_i\overset{\text{i.i.d.}}{\sim}N(0,1),
\]
using the same ground truth $(u^\dagger,p^\dagger)$ for all $\delta$ within each trial.
Importantly, $\delta$ does \emph{not} enter the likelihood (and hence does not redefine the posterior);
for each $\delta$ we rerun inference under the \emph{same} Bayesian model, with the only change being the
realized dataset $y^\delta$.
We then quantify sensitivity by comparing the resulting summaries
$\hat u_\star^\delta,\hat p_\star^\delta,\hat a_\star^\delta$ against their baseline counterparts
$\hat u_\star^0,\hat p_\star^0,\hat a_\star^0$ using discrete $L^2(\Gamma)$-type norms.

\paragraph{Metrics.}
For fields on the mesh we use the discrete $L^2(\Gamma)$ norm induced by the finite-element mass matrix,
\[
\|w\|_{L^2_h(\Gamma)} := (w^\top C w)^{1/2}.
\]
We quantify the perturbation size by the relative data change
\[
  \mathrm{RelPert}_\delta(y)
  :=
  \frac{\|y^\delta-y^0\|_{L^2_h(\Gamma)}}{\|y^0\|_{L^2_h(\Gamma)}}.
\]
For a generic posterior summary $v^\delta_\star$ with baseline counterpart $v^0_\star$, define the relative change
\[
  \mathrm{RelChg}_\delta(v_\star)
  :=
  \frac{\|v^\delta_\star-v^0_\star\|_{L^2_h(\Gamma)}}{\|v^0_\star\|_{L^2_h(\Gamma)}},
\]
and the empirical amplification factor
\[
  A_\delta(v_\star)
  :=
  \frac{\|v^\delta_\star-v^0_\star\|_{L^2_h(\Gamma)}}{\|y^\delta-y^0\|_{L^2_h(\Gamma)}}.
\]
Here, we report $\mathrm{RelPert}_\delta(y)$, the state-level relative change $\mathrm{RelChg}_\delta(\hat p_\star)$,
and amplification factors $A_\delta(\hat u_\star)$ and $A_\delta(\hat a_\star)$.
We focus on amplification in the inferred coefficient summaries, since this directly probes parameter stability
in the inverse problem; for the state $\hat p_\star$ we report the scale-free relative change $\mathrm{RelChg}_\delta(\hat p_\star)$.

\paragraph{Perturbation levels.}
In the elliptic experiments we take $\delta\in\{0.1,0.2,0.5\}$.
In the fractional elliptic experiments we take $\delta\in\{1.1,1.3,1.5\}$ to probe a stronger perturbation regime.

\paragraph{Results and interpretation.}
Tables~\ref{tab:wp-standard-compact}--\ref{tab:wp-fractional-compact} report the perturbation-to-output response
(mean $\pm$ sd across trials) and allow us to separate (i) the imposed perturbation size from (ii) the response of posterior summaries.

\smallskip
\noindent\textbf{(i) Perturbation size and calibration.}
In the elliptic case (Table~\ref{tab:wp-standard-compact}, first row), the relative data perturbation
$\mathrm{RelPert}_\delta(y)$ increases in a nearly linear fashion with $\delta$:
it is about $0.32\%$ at $\delta=0.1$, $0.63\%$ at $\delta=0.2$, and $1.58\%$ at $\delta=0.5$.
In the fractional elliptic case (Table~\ref{tab:wp-fractional-compact}, first row), the imposed perturbations are intentionally larger,
and $\mathrm{RelPert}_\delta(y)$ remains monotone and approximately linear in $\delta$,
rising from $2.86\%$ at $\delta=1.1$ to $3.82\%$ at $\delta=1.5$ (with slightly larger trial-to-trial variability).
Overall, the  amplification factors remain small at all perturbation levels.

\smallskip
\noindent\textbf{(ii) Response of posterior summaries: state-level relative changes.}
A robust, scale-free indicator of stability is the relative change in the state summary $\mathrm{RelChg}_\delta(\hat p_\star)$.
In the elliptic case, $\mathrm{RelChg}_\delta(\hat p_{\mathrm{mean}})\approx 3.1\%$ and
$\mathrm{RelChg}_\delta(\hat p_{\mathrm{mode}})\approx 4.5$--$4.8\%$, with no systematic deterioration as $\delta$ increases.
In the fractional elliptic case, despite data perturbations at the few-percent level,
the state summaries vary at the one-percent level:
$\mathrm{RelChg}_\delta(\hat p_{\mathrm{mean}})\approx 0.78\%$--$1.06\%$ and
$\mathrm{RelChg}_\delta(\hat p_{\mathrm{mode}})\approx 1.04\%$--$1.19\%$ across $\delta\in\{1.1,1.3,1.5\}$.
Thus, in both regimes the state summaries remain stable in relative terms over the tested perturbation ranges.

\smallskip
\noindent\textbf{(iii) Amplification factors for coefficient summaries and stable input--output response.}
The amplification factors $A_\delta(\hat u_\star)$ and $A_\delta(\hat a_\star)$ compare the absolute change in inferred coefficients
to the absolute size of the imposed data perturbation.
In the elliptic case, the mean amplification decreases as $\delta$ increases
(e.g., $A_\delta(\hat u_{\mathrm{mean}})$ drops from $0.367$ at $\delta=0.1$ to $0.0798$ at $\delta=0.5$,
and similarly $A_\delta(\hat a_{\mathrm{mean}})$ from $0.298$ to $0.0651$).
This pattern is consistent with a benign \emph{ratio effect}:
$\|y^\delta-y^0\|_{L^2_h(\Gamma)}$ grows approximately linearly in $\delta$ by construction,
while the induced coefficient changes $\|\hat u_\star^\delta-\hat u_\star^0\|_{L^2_h(\Gamma)}$ and
$\|\hat a_\star^\delta-\hat a_\star^0\|_{L^2_h(\Gamma)}$ remain of comparable magnitude over the tested range.
Accordingly, smaller perturbations naturally yield more variable ratios (a smaller denominator),
without indicating any deterioration of the underlying mapping from data to posterior summaries.
More importantly for our purposes, across the tested levels we do not see the amplification
\emph{escalate} as the perturbation size is reduced, which would be a typical numerical signature of unstable dependence.

In the fractional elliptic case, we probe a substantially larger perturbation regime with $\mathrm{RelPert}_\delta(y)\approx 2.86\%$--$3.82\%$,
yet the coefficient amplification factors remain uniformly small:
$A_\delta(\hat u_\star)$ and $A_\delta(\hat a_\star)$ stay on the order of $10^{-2}$
(about $7\times 10^{-3}$--$9\times 10^{-3}$ across $\delta$ and across $\star\in\{\mathrm{mean},\mathrm{mode}\}$).
Thus, even when the data are perturbed at the few-percent level, the induced coefficient changes are only a tiny fraction of the input change.
Taken together, the elliptic and fractional elliptic results provide empirical support for a stable perturbation-to-solution response
of the posterior coefficient summaries over the ranges of $\delta$ considered, rather than only in an infinitesimal-perturbation limit.

\smallskip
\noindent\textbf{(iv) Variability across trials.}
Several metrics exhibit standard deviations comparable to or larger than the corresponding means, especially in the elliptic regime
at the smallest perturbation level.
This is consistent with the ratio-based nature of $A_\delta(\cdot)$: when $\|y^\delta-y^0\|_{L^2_h(\Gamma)}$ is tiny,
moderate trial-to-trial fluctuations in the numerator translate into noticeable variability in the ratio.
Importantly, this variability does not coincide with a systematic deterioration in the state-level relative changes
$\mathrm{RelChg}_\delta(\hat p_\star)$, nor with growth of the coefficient amplification factors as $\delta$ decreases.
Overall, the numerical evidence in Tables~\ref{tab:wp-standard-compact}--\ref{tab:wp-fractional-compact}
is consistent with stable dependence of posterior coefficient summaries on the observed data within the tested perturbation ranges.

\begin{table}[htp]
\centering
\scriptsize
\setlength{\tabcolsep}{5pt}
\renewcommand{\arraystretch}{1.15}
\caption{\black Elliptic experiment ($\beta=1$): perturbation-to-output response (mean $\pm$ sd across $30$ trials). We report the relative data perturbation $\mathrm{RelPert}_\delta(y)$, the relative change $\mathrm{RelChg}_\delta(\hat p_\star)$, and amplification factors $A_\delta(\hat u_\star)$ and $A_\delta(\hat a_\star)$ for $\star\in\{\mathrm{mean},\mathrm{mode}\}$, where ``mode'' denotes the posterior marginal mode summary.\nc}
\label{tab:wp-standard-compact}
\begin{tabular}{@{}lccc@{}}
\toprule
 & $\delta=0.1$ & $\delta=0.2$ & $\delta=0.5$ \\
\midrule
$\mathrm{RelPert}_\delta(y)$ 
& $3.21\times 10^{-3}\pm 7.94\times 10^{-5}$ 
& $6.34\times 10^{-3}\pm 1.69\times 10^{-4}$ 
& $1.58\times 10^{-2}\pm 3.91\times 10^{-4}$ \\
\addlinespace[2pt]

\multicolumn{4}{@{}l}{\emph{Posterior mean summary} ($\star=\mathrm{mean}$)}\\
$\mathrm{RelChg}_\delta(\hat p_{\mathrm{mean}})$ 
& $0.0307\pm 0.0970$ & $0.0310\pm 0.0970$ & $0.0314\pm 0.0968$ \\
$A_\delta(\hat u_{\mathrm{mean}})$ 
& $0.367\pm 1.02$ & $0.193\pm 0.535$ & $0.0798\pm 0.213$ \\
$A_\delta(\hat a_{\mathrm{mean}})$ 
& $0.298\pm 0.635$ & $0.156\pm 0.334$ & $0.0651\pm 0.133$ \\
\addlinespace[2pt]

\multicolumn{4}{@{}l}{\emph{Posterior marginal mode summary} ($\star=\mathrm{mode}\,$)}\\
$\mathrm{RelChg}_\delta(\hat p_{\mathrm{mode}})$ 
& $0.0481\pm 0.106$ & $0.0454\pm 0.107$ & $0.0461\pm 0.106$ \\
$A_\delta(\hat u_{\mathrm{mode}})$ 
& $0.422\pm 1.08$ & $0.221\pm 0.570$ & $0.0899\pm 0.227$ \\
$A_\delta(\hat a_{\mathrm{mode}})$ 
& $0.344\pm 0.655$ & $0.180\pm 0.345$ & $0.0735\pm 0.138$ \\
\bottomrule
\end{tabular}
\end{table}

\begin{table}[htp]
\centering
\scriptsize
\setlength{\tabcolsep}{5pt}
\renewcommand{\arraystretch}{1.15}
\caption{\black Fractional elliptic experiment ($\beta=3/2$): perturbation-to-output response (mean $\pm$ sd across $10$ trials). The reported metrics and the display convention are the same as in Table~\ref{tab:wp-standard-compact}.\nc}
\label{tab:wp-fractional-compact}
\begin{tabular}{@{}lccc@{}}
\toprule
 & $\delta=1.1$ & $\delta=1.3$ & $\delta=1.5$ \\
\midrule
$\mathrm{RelPert}_\delta(y)$ 
& $0.0286\pm 1.69\times 10^{-4}$ 
& $0.0337\pm 2.48\times 10^{-3}$ 
& $0.0382\pm 3.08\times 10^{-3}$ \\
\addlinespace[2pt]

\multicolumn{4}{@{}l}{\emph{Posterior mean summary} ($\star=\mathrm{mean}$)}\\
$\mathrm{RelChg}_\delta(\hat p_{\mathrm{mean}})$ 
& $7.78\times 10^{-3}\pm 2.68\times 10^{-3}$ 
& $9.55\times 10^{-3}\pm 8.23\times 10^{-4}$ 
& $1.06\times 10^{-2}\pm 9.26\times 10^{-4}$ \\
$A_\delta(\hat u_{\mathrm{mean}})$ 
& $7.32\times 10^{-3}\pm 8.01\times 10^{-4}$ 
& $7.23\times 10^{-3}\pm 1.01\times 10^{-3}$ 
& $7.17\times 10^{-3}\pm 1.32\times 10^{-3}$ \\
$A_\delta(\hat a_{\mathrm{mean}})$ 
& $7.32\times 10^{-3}\pm 7.90\times 10^{-4}$ 
& $7.24\times 10^{-3}\pm 1.02\times 10^{-3}$ 
& $7.16\times 10^{-3}\pm 1.31\times 10^{-3}$ \\
\addlinespace[2pt]

\multicolumn{4}{@{}l}{\emph{Posterior marginal mode summary} ($\star=\mathrm{mode}$)}\\
$\mathrm{RelChg}_\delta(\hat p_{\mathrm{mode}})$ 
& $1.04\times 10^{-2}\pm 4.17\times 10^{-3}$ 
& $1.11\times 10^{-2}\pm 1.52\times 10^{-3}$ 
& $1.19\times 10^{-2}\pm 1.23\times 10^{-3}$ \\
$A_\delta(\hat u_{\mathrm{mode}})$ 
& $8.61\times 10^{-3}\pm 8.64\times 10^{-4}$ 
& $8.11\times 10^{-3}\pm 1.59\times 10^{-3}$ 
& $8.15\times 10^{-3}\pm 9.90\times 10^{-4}$ \\
$A_\delta(\hat a_{\mathrm{mode}})$ 
& $8.60\times 10^{-3}\pm 8.74\times 10^{-4}$ 
& $8.11\times 10^{-3}\pm 1.59\times 10^{-3}$ 
& $8.14\times 10^{-3}\pm 9.84\times 10^{-4}$ \\
\bottomrule
\end{tabular}
\end{table}

\newpage
\subsection{Reconstruction Metrics Under an Oscillatory Forcing}
\label{app:metrics-oscillatory-rhs}
In the numerical experiments in Section~\ref{ssec:implementation}, the forcing is chosen as
$f(x)=z_1^2(x)-z_2^2(x)$, which is smooth and effectively low-frequency in the embedded coordinates.
To assess the robustness of the inference--reconstruction pipeline with respect to the forcing,
we repeat the same set of trials with a more oscillatory right-hand side.
Specifically, let $z(x)=(z_1(x),z_2(x))$ denote the embedded coordinates of a point $x\in\Gamma$,
and define their rescaled versions $\tilde z_1(x),\tilde z_2(x)\in[0,1]$ by
\[
\tilde z_i(x) := \frac{z_i(x)-\min_{\Gamma} z_i}{\max_{\Gamma} z_i-\min_{\Gamma} z_i},
\qquad i=1,2.
\]
We then consider the oscillatory forcing
\begin{equation}
\label{eq:osc_rhs}
    f_{\mathrm{osc}}(x)
= 10\,\sin \bigl(4\pi \,\tilde z_1(x)\bigr)\,
      \cos \bigl(4\pi \,\tilde z_2(x)\bigr).
\end{equation}
In the implementation, we post-process $f_{\mathrm{osc}}$ at the finite-element level:
we remove its $L_h^2(\Gamma)$-weighted mean (so that the discrete forcing has zero average with respect to the mass matrix),
and we rescale it by its $L^2_h(\Gamma)$ norm (as defined in Section~\ref{ssec:implementation}), $\|g\|_{L^2_h(\Gamma)}^2:=g^\top C g$, to keep the forcing on a numerically stable scale.

The oscillatory forcing $f_{\mathrm{osc}}$ (cf. Equation~\eqref{eq:osc_rhs}) induces higher spatial variation in the state $p$ along the network,
and therefore puts more weight on the ability of the forward solver and the inverse procedure to
resolve comparatively fine-scale features.
In particular, compared with the baseline low-frequency forcing (see Section~\ref{ssec:implementation}), the induced observations may be more sensitive
to localized discrepancies in the reconstructed coefficient, while the posterior may also exhibit slightly larger
trial-to-trial variability due to the interaction between oscillations and observational noise.

\begin{figure}[H]
    \centering

    \begin{subfigure}[t]{0.32\textwidth}
        \centering
        \includegraphics[width=\linewidth]{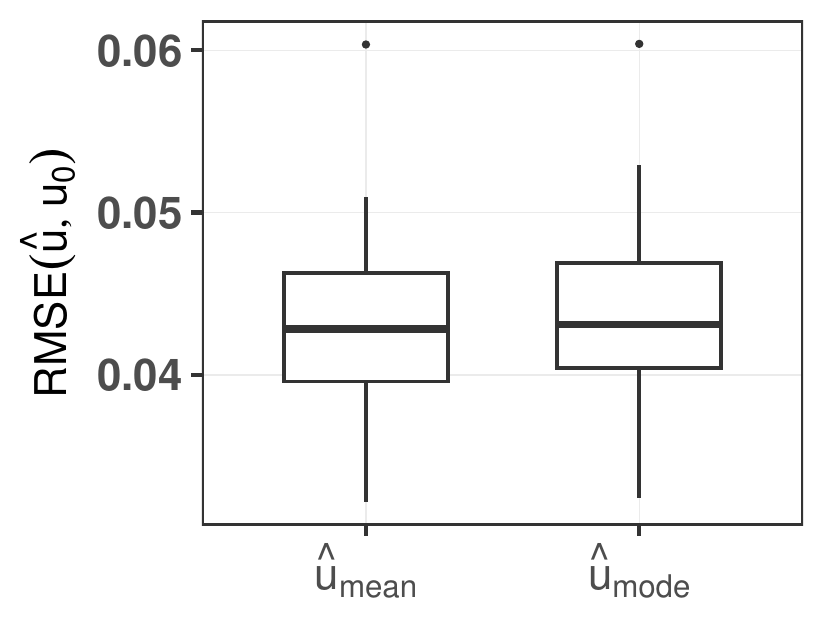}
        \caption{$\mathrm{RMSE}(\hat u,u_0)$}
        \label{fig:box-rmse-u-rhs}
    \end{subfigure}
    \hfill
    \begin{subfigure}[t]{0.32\textwidth}
        \centering
        \includegraphics[width=\linewidth]{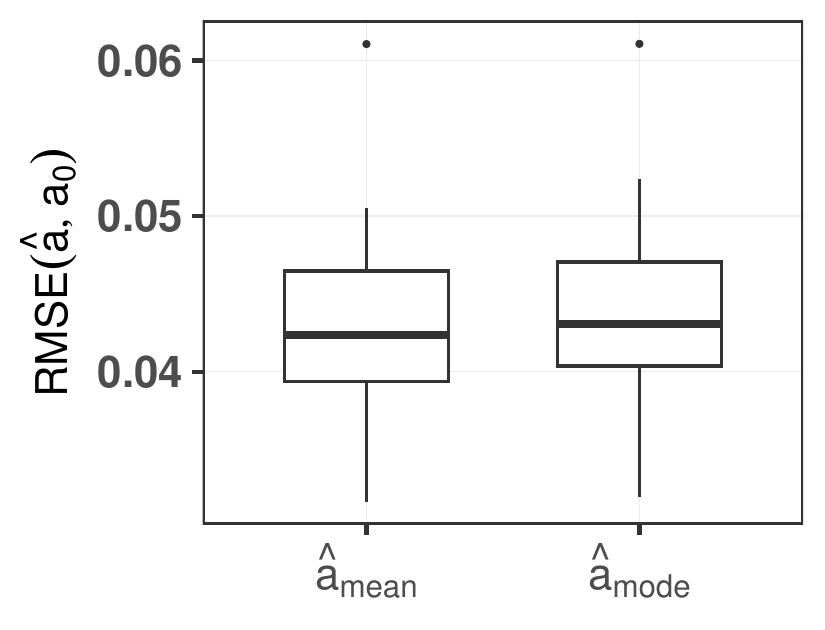}
        \caption{$\mathrm{RMSE}(\hat a,a_0)$}
        \label{fig:box-rmse-a-rhs}
    \end{subfigure}
    \hfill
    \begin{subfigure}[t]{0.32\textwidth}
        \centering
        \includegraphics[width=\linewidth]{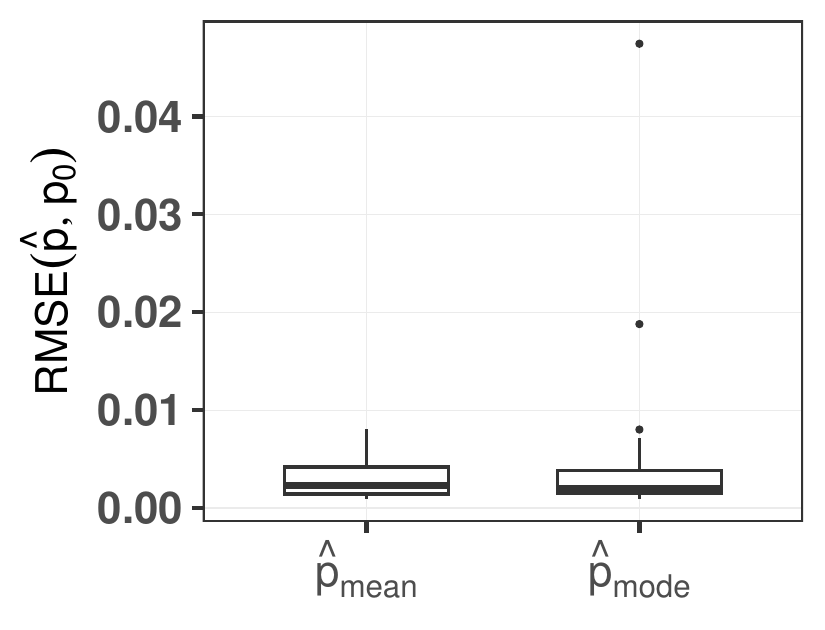}
        \caption{$\mathrm{RMSE}(\hat p,p_0)$}
        \label{fig:box-rmse-p-rhs}
    \end{subfigure}
    \vspace{0.25cm}

    \begin{subfigure}[t]{0.32\textwidth}
        \centering
        \includegraphics[width=\linewidth]{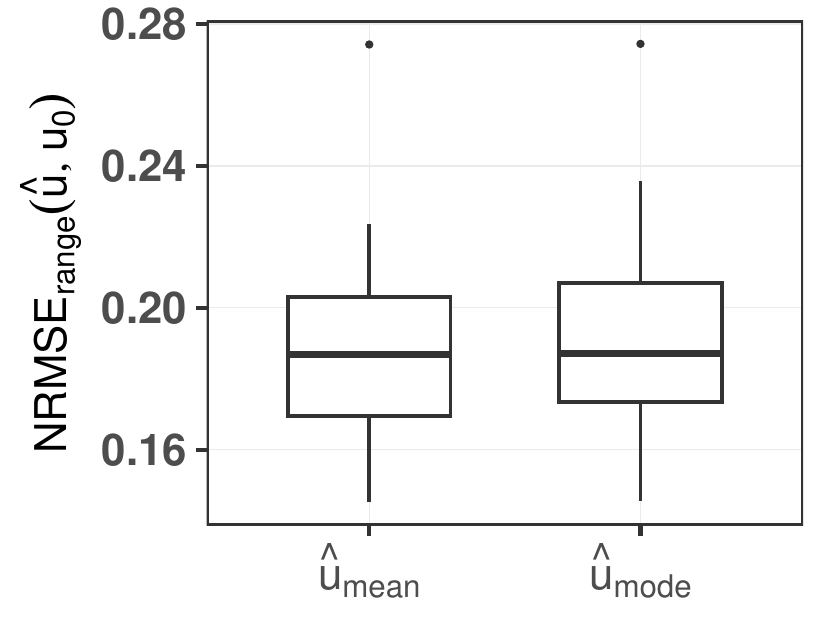}
        \caption{$\mathrm{NRMSE}_{\rm range}(\hat u,u_0)$}
        \label{fig:box-nrmse-u-range-rhs}
    \end{subfigure}
    \hfill
    \begin{subfigure}[t]{0.32\textwidth}
        \centering
        \includegraphics[width=\linewidth]{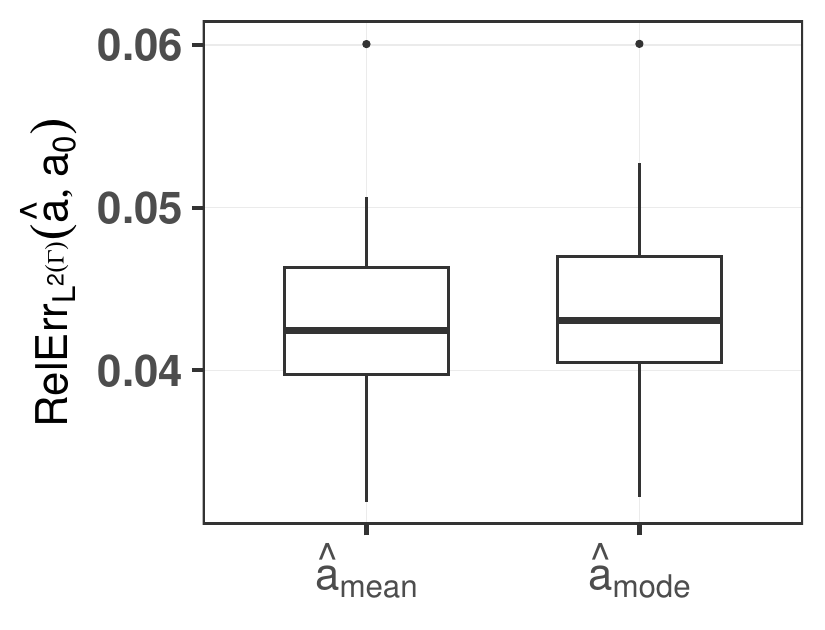}
        \caption{$\mathrm{RelErr}_{L^2(\Gamma)}(\hat a,a_0)$}
        \label{fig:box-relerr-a-rhs}
    \end{subfigure}
    \hfill
    \begin{subfigure}[t]{0.32\textwidth}
        \centering
        \includegraphics[width=\linewidth]{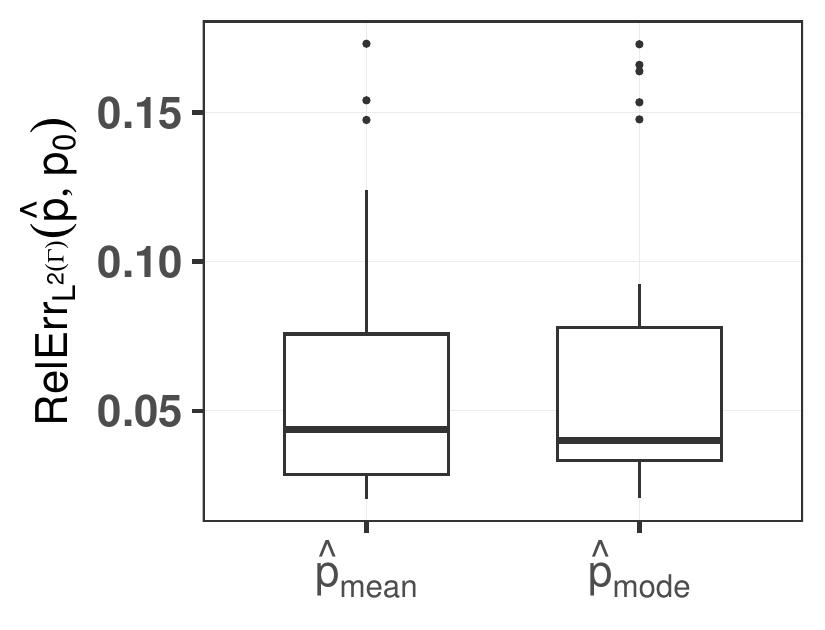}
        \caption{$\mathrm{RelErr}_{L^2(\Gamma)}(\hat p,p_0)$}
        \label{fig:box-relerr-p-rhs}
    \end{subfigure}

    \caption{
    \black Boxplots of reconstruction errors over repeated trials for the \emph{oscillatory forcing term} experiment described above.
    The reported metrics and the display convention are the same as in Figure~\ref{fig:boxplots-2plus3-errors}.\nc
    }
    \label{fig:boxplots-2plus3-errors-rhs}
\end{figure}

Figure~\ref{fig:boxplots-2plus3-errors-rhs} summarizes the reconstruction quality under the oscillatory forcing described above.
Overall, the distributions of the parameter errors \(\mathrm{RMSE}(\hat u,u_0)\) and \(\mathrm{NRMSE}_{\mathrm{range}}(\hat u,u_0)\) in Figures~\ref{fig:box-rmse-u-rhs} and~\ref{fig:box-nrmse-u-range-rhs} remain tightly concentrated and of comparable magnitude to those in the baseline experiment of Figure~\ref{fig:boxplots-2plus3-errors}.
This indicates that, at the tested noise level and sensor configuration, replacing a smooth low-frequency forcing by a more oscillatory input does not significantly impact the accuracy of the recovered \(u\).

Figures~\ref{fig:box-rmse-a-rhs} and~\ref{fig:box-relerr-a-rhs} show that the same conclusion carries over to the conductivity \(a=\exp(u)\).
In particular, the relative \(L^2(\Gamma)\) errors for \(\hat a\) stay small and exhibit limited trial-to-trial variability, suggesting that the exponential nonlinearity does not amplify the posterior uncertainty in \(u\) into a substantially less stable reconstruction of \(a\) in this setting.
This is consistent with the fact that the posterior dispersion of \(u\) (as reflected indirectly by the spread of the error boxplots) remains controlled under the oscillatory forcing.

For the state variable \(p\), Figures~\ref{fig:box-rmse-p-rhs} and~\ref{fig:box-relerr-p-rhs} exhibit a slightly different pattern: while the median errors and interquartile ranges are still comparable to the baseline, a small number of trials produce visibly larger errors (outliers), especially in the relative \(L^2(\Gamma)\) metric.
This behavior is expected for a higher-frequency forcing, since oscillatory inputs typically excite more spatial modes of the PDE solution and can increase sensitivity of the state reconstruction to localized mismatches in \(\hat u\) (and hence in \(\hat a\)).
Importantly, these larger errors remain isolated to a few trials rather than shifting the bulk of the distribution, so the overall reconstruction quality remains stable.

Finally, across all six plots in Figure \ref{fig:boxplots-2plus3-errors-rhs}, the posterior mean and posterior mode yield nearly indistinguishable error distributions.
Beyond indicating practical robustness to the choice of point estimator, this similarity suggests that the posterior is not exhibiting strong asymmetry or multimodality at the level relevant for reconstruction: the dominant source of variability appears to be trial-to-trial fluctuations induced by the data realization and the inherent posterior uncertainty, rather than a systematic difference between \(\hat u_{\mathrm{mean}}\) and \(\hat u_{\mathrm{mode}}\).
\endgroup

\begingroup
\color{black}
\subsection{Finite-element Discretization and Rational Approximation}\label{app:fem-rational}

We briefly summarize the numerical discretization used in the experiments.
Our starting point is the general elliptic differential operator framework in Subsection~\ref{ssec:quantum_graphs}.
In particular, for a strictly positive diffusion coefficient $a\in C(\Gamma)$ bounded away from zero
and $\kappa>0$ in \eqref{eq:edge_ODE}, the second-order elliptic operator $\mathcal{L}$ equipped with the Kirchhoff
vertex conditions \eqref{eq:kirchhoff} is self-adjoint and strictly positive on $L^2(\Gamma)$.
Consequently, its fractional powers $\mathcal{L}^\beta$ are well defined via functional calculus.
We consider the (fractional) elliptic differential equation
\[
\mathcal{L}^{\beta}p = f,\qquad f\in L^2(\Gamma),\quad \beta\ge 1,
\]
and describe its finite-element approximation below.

\paragraph{Finite-element mesh and continuous piecewise linear space on $\Gamma$.}
Let $\{x_i\}_{i=1}^{N_h}$ be the mesh nodes on $\Gamma$, obtained by subdividing each edge into uniform subintervals,
as described in Subsection~\ref{ssec:implementation}.
Following \cite[Section~6.2]{bolin2024regularity}, we use continuous piecewise linear finite elements on $\Gamma$.
Concretely, for each edge $e$ (identified with an interval $(0,\ell_e)$), we introduce a uniform partition
$0=t_{e,0}<t_{e,1}<\cdots<t_{e,M_e}=\ell_e$ with step size $h_e=\ell_e/M_e$.
On the interior nodes $\{t_{e,j}\}_{j=1}^{M_e-1}$, we use the standard edgewise hat basis functions
\[
\varphi^{h}_{e,j}(t):=
\begin{cases}
\dfrac{t-t_{e,j-1}}{h_e}, & t\in[t_{e,j-1},t_{e,j}],\\[4pt]
\dfrac{t_{e,j+1}-t}{h_e}, & t\in[t_{e,j},t_{e,j+1}],\\[4pt]
0, & \text{otherwise},
\end{cases}
\qquad j=1,\dots,M_e-1,
\]
extended by zero to all other edges.
In addition, for each vertex $v\in V$ we include a vertex-centered hat function $\varphi^{h}_{v}$ defined by
$\varphi^{h}_{v}(v)=1$, $\varphi^{h}_{v}(w)=0$ for all other vertices $w\neq v$, and whose restriction to each incident edge
$e\in E_v$ is linear in the arc-length coordinate:
if $t=0$ corresponds to $v$ on edge $e$, then
\[
\varphi^{h}_{v}\big|_{e}(t)=\max\Bigl\{1-\frac{t}{h_e},\,0\Bigr\},\qquad t\in[0,\ell_e],
\]
and similarly if $t=\ell_e$ corresponds to $v$.
These vertex basis functions ``glue'' the edgewise bases and enforce global continuity at vertices.
We denote the resulting nodal finite element space by
\[
V_h
:= \Bigl(\bigoplus_{e\in E}\mathrm{span}\{\varphi^{h}_{e,j}\}_{j=1}^{M_e-1}\Bigr)\;\oplus\;\mathrm{span}\{\varphi^{h}_{v}\}_{v\in V}
\;=\;\mathrm{span}\{\varphi_i^{h}\}_{i=1}^{N_h}\subset H^1(\Gamma),
\]
where $\{\varphi_i^{h}\}_{i=1}^{N_h}$ is the global nodal basis satisfying $\varphi_i^{h}(x_j)=\delta_{ij}$.
In particular,
\[
N_h = |V| + \sum_{e\in E}(M_e-1).
\]

\paragraph{Discrete representation of coefficient and solution.}
We represent the diffusion coefficient and the solution by their nodal degrees of freedom in $V_h$.
For a diffusion coefficient $a$ on $\Gamma$, we denote by $a_h\in V_h$ its finite-element representation (e.g.\ the nodal interpolant), and, by a slight abuse of notation, also write $a := (a_h(x_1),\dots,a_h(x_{N_h}))^\top\in\mathbb{R}^{N_h}$ for its nodal coefficient vector.
Similarly, for $p_h=\sum_{i=1}^{N_h}p_i\varphi_i^{h}\in V_h$ we store
$p=(p_1,\dots,p_{N_h})^\top\in\mathbb{R}^{N_h}$.
In the numerical experiments we generate a ground-truth diffusion coefficient (or parameter) and compute the corresponding state;
see Subsection~\ref{ssec:implementation} for details and the chosen mesh sizes.

\paragraph{Assembly of discrete matrices.}
Let $a_h\in V_h$ be the strictly positive finite-element interpolant of the diffusion coefficient $a$ on $\Gamma$.
We define the mass matrix $C\in\mathbb{R}^{N_h\times N_h}$ and the weighted stiffness matrix
$G(a)\in\mathbb{R}^{N_h\times N_h}$ by
\begin{equation}\label{eq:mass-stiff-general}
C_{ij}=\int_\Gamma \varphi_i^{h}(x)\varphi_j^{h}(x)\,{\rm{dx}},\qquad
G_{ij}(a)=\int_\Gamma a_h(x)\,\bigl(\varphi_i^{h}\bigr)'(x)\,\bigl(\varphi_j^{h}\bigr)'(x)\, {\rm{dx}},
\end{equation}
where $\bigl(\varphi_i^{h}\bigr)'$ denotes the edgewise derivative and the integral over $\Gamma$ is the sum of edge integrals.
We then define the discrete elliptic operator
\begin{equation}\label{eq:Lh-general}
L_h(a)=\kappa^2 C + G(a).
\end{equation}
These matrices arise from the Galerkin discretization of the bilinear form
\(
\mathsf{B}(p,v):=\int_\Gamma a\,p'v' + \kappa^2\int_\Gamma pv
\)
associated with $\mathcal{L}$.
In particular, seeking $p_h=\sum_{j=1}^{N_h}p_j\varphi_j^h\in V_h$ such that
$\mathsf{B}(p_h,v_h)=(f,v_h)_{L^2(\Gamma)}$ for all $v_h\in V_h$, and testing with $v_h=\varphi_i^h$, yields
$\sum_{j=1}^{N_h} (L_h(a))_{ij}p_j=(f,\varphi_i^h)_{L^2(\Gamma)}$.
Moreover, for any $v_h=\sum_i v_i\varphi_i^{h}\in V_h$, the discrete $L^2(\Gamma)$-norm is computed via the mass matrix:
\begin{equation}\label{eq:disc-L2}
\|v_h\|_{L^2(\Gamma)}^2 = v^\top C v .
\end{equation}

\paragraph{Elliptic case $\beta=1$.}
Let $(f_h)_i=\int_\Gamma f(x)\varphi_i^{h}(x)\,dx$ be the finite-element load vector.
Then the Galerkin solution $p_h\in V_h$ is obtained by solving for $p=(p_1,\ldots,p_{N_h})^\top$ the linear system
\begin{equation}\label{eq:elliptic-linear-system}
L_h(a)\,p = f_h
\end{equation}
and setting $p_h = \sum_{i=1}^{N_h} p_i \varphi_i^h$.

\paragraph{Fractional elliptic case $\beta>1$.}
For $\beta>1$, discretization of the fractional operator via spectral decomposition is computationally prohibitive for large $N_h$. We therefore adopt the operator-based rational approximation of \cite{bolin2019rational},
as implemented in the \texttt{rSPDE} and \texttt{MetricGraph} packages, which retains sparsity and avoids computing eigenpairs.

Introduce the discrete operator $A_h(a):=C^{-1}L_h(a)$, and notice that it represents $\mathcal{L}$ on $V_h$ with respect to the
$L^2(\Gamma)$ inner product. 
Indeed, letting $u_h=\sum_i u_i\varphi_i^h$ and $v_h=\sum_i v_i\varphi_i^h$, we have
$(u_h,v_h)_{L^2(\Gamma)}=u^\top C v$ and $\mathsf{B}(u_h,v_h)=u^\top L_h(a)v$, so $A_h(a)$ is characterized by
\[
(A_h(a)u)^\top C v \;=\; u^\top L_h(a)\,v \qquad \forall\,u,v\in\mathbb{R}^{N_h},
\]
which yields $C A_h(a)=L_h(a)$, i.e.\ $A_h(a)=C^{-1}L_h(a)$.
Similarly, introduce $b_h:=C^{-1}f_h$ and notice that $b_h$ is the $L^2$--Riesz representative of the load functional: $C b_h=f_h$, so that
$(b_h,v_h)_{L^2(\Gamma)}=f_h^\top v$ for all $v\in\mathbb{R}^{N_h}$.

Write $\beta=m+\beta_0$ with $m=\lfloor \beta\rfloor\in\mathbb{N}$ and $\beta_0\in(0,1)$.
Then, 
\[
A_h(a)^{-\beta}b_h \;=\; A_h(a)^{-m}\,A_h(a)^{-\beta_0}b_h,
\]
so it suffices to approximate $A_h(a)^{-\beta_0}b_h$ for $\beta_0\in(0,1)$.

We now detail the rational approximation of $A_h(a)^{-\beta_0}$ for $\beta_0\in(0,1)$, following \cite{bolin2019rational}.
Let $0<\lambda_{h,1}\le \cdots \le \lambda_{h,N_h}$ denote the generalized eigenvalues of the pair $(L_h(a),C)$,
i.e., $L_h(a)v=\lambda\,Cv$.
Equivalently, these are the positive eigenvalues of $A_h(a)=C^{-1}L_h(a)$.
Set
\[
J_h := [\lambda_{h,N_h}^{-1},\,\lambda_{h,1}^{-1}],
\]
so that the spectrum of $A_h(a)^{-1}$ is contained in the bounded interval $J_h$.
Writing
\[
A_h(a)^{-\beta_0} = f\bigl(A_h(a)^{-1}\bigr),\qquad f(z):=z^{\beta_0},
\]
motivates constructing a rational approximation of $f$ over $J_h$.

Fix an integer approximation order $M\ge 1$.
Following \cite[Section~3]{bolin2019rational}, we decompose
\[
f(z)=z^{\beta_0}= z\,\widehat f(z),\qquad \widehat f(z):=z^{\beta_0-1}.
\]
This factorization separates off one power of $z$ and localizes the non-integer behavior at $z=0$
to the single factor $\widehat f(z)=z^{\beta_0-1}$, which has only a weak algebraic singularity.
As a result, the ensuing best uniform rational approximation on $J_h$ is numerically more stable.
We then approximate $\widehat f$ uniformly on $J_h$ by a rational function
$\widehat r_M(z)=q_1(z)/q_2(z)$ with $\deg(q_1)=M$ and $\deg(q_2)=M+1$, i.e.,
\begin{equation}\label{eq:best-rational-approx}
\widehat r_M \in 
\underset{\substack{\deg(q_1)=M\\ \deg(q_2)=M+1}}{\operatorname{arg\,min}}
\left\| z^{\beta_0-1}-\frac{q_1(z)}{q_2(z)}\right\|_{L^\infty(J_h)} .
\end{equation}
Setting $r_M(z):=z\,\widehat r_M(z)$ yields $r_M(z)\approx z^{\beta_0}$ on $J_h$. 
Using $z^{-\beta_0}=f(z^{-1})$, we obtain the induced rational approximation
\[
z^{-\beta_0} \;\approx\; r_M(z^{-1})
= \widehat r_M(z^{-1})\,z^{-1}
=\frac{q_1(z^{-1})}{q_2(z^{-1})\,z}
=\frac{p_r(z)}{p_\ell(z)},
\]
where $p_r$ and $p_\ell$ are obtained by expanding the numerator and denominator as polynomials in $z$:
\[
p_r(z):=\sum_{i=0}^{M} c_{r,i}\,z^{\,M-i},
\qquad
p_\ell(z):=\sum_{j=0}^{M+1} c_{\ell,j}\,z^{\,M+1-j}.
\]
Consequently, we obtain the operator approximation 
\[
A_h(a)^{-\beta_0}\;\approx\; \frac{p_r\bigl(A_h(a)\bigr)}{p_\ell\bigl(A_h(a)\bigr)}
= p_\ell \,\bigl(A_h(a)\bigr)^{-1} \, p_r \bigl(A_h(a)\bigr).
\]
Therefore, we approximate $A_h(a)^{-\beta_0}b_h$ by solving the non-fractional equation
\begin{equation}\label{eq:frac-rational}
P_{\ell,h}\,q_h^{(\beta_0)} \;=\; P_{r,h}\,b_h,
\qquad
P_{\ell,h}:=p_\ell \bigl(A_h(a)\bigr),\quad P_{r,h}:=p_r \bigl(A_h(a)\bigr),
\end{equation}
and setting $q_h^{(\beta_0)}\approx A_h(a)^{-\beta_0}b_h$.
As shown in \cite{bolin2019rational}, applying $P_{\ell,h}^{-1}$ amounts to solving a small number of sparse systems.

Finally, writing $\beta=m+\beta_0$ with $m=\lfloor \beta\rfloor\in\mathbb{N}$, we approximate the fractional solution vector by
\[
p \;\approx\; A_h(a)^{-m}\,q_h^{(\beta_0)},
\]
where $A_h(a)^{-m}=(L_h(a)^{-1}C)^m$ can be applied by $m$ successive solves with $L_h(a)$.

\endgroup

\end{document}